\documentclass{article}
\usepackage{amsmath,amssymb,amsfonts}
\usepackage[default]{fontsetup}
\usepackage{array}
\usepackage{supertabular}
\newcommand{\cellwidth}{0.235\textwidth}
\newcommand{\dcellwidth}{0.460\textwidth}
\newcommand\cellcenter[1]{\begin{minipage}{\cellwidth}\begin{center}\ \\ #1\\[-1.4ex]\ \end{center}\end{minipage}}
\newcommand\dcellcenter[1]{\begin{minipage}{\dcellwidth}\begin{center}\ \\ #1\\[-1.4ex]\ \end{center}\end{minipage}}
\newcommand\cc{\cellcenter}
\newcommand\dcc{\dcellcenter}
\usepackage{enumitem}
\usepackage{hhline}
\usepackage{geometry}

\makeatletter
\newcommand\arraybslash{\let\\\@arraycr}
\makeatother
\usepackage[polutonikogreek, english]{babel}
\usepackage{epigraph,calc}

\title{The Mysteries of Plato's \textit{Parmenides} Dispersed: From musical intervals to contacts/``hapseis'' to ``logoi''}
\author{Stelios Negrepontis}
\date{}

\begin{document}
\maketitle

% \begin{epigraph}
\setlength{\epigraphwidth}{\widthof{Goddess [Athena] dispersed the mist,}}
\epigraph{θεὰ σκέδασ' ἠέρα, εἴσατο δὲ χθών \\
Goddess [Athena] dispersed the mist, \\
and the land was recognized}{Odusseia, Book XIII, line 352}
%\end{epigraph}

\setcounter{section}{-1}
\section{Introduction}

Plato’s \textit{Parmenides} is the absolutely most difficult and baffling of Platonic dialogues, full of seeming contradictions (Section 1). 
We will offer the interpretation that we believe is the definitive one. 

Before approaching \textit{Parmenides}, it is helpful to obtain a
basic understanding of the basics of Plato’s philosophy obtained by
reading and studying some other, more accessible and less forbidding,
Platonic dialogues. The most helpful of these dialogues are
\textit{Theaetetus}, \textit{Sophist}, and \textit{Meno}. Our study of
these dialogues (2012~\cite{Negrepontis2012},
2018~\cite{Negrepontis2018}, 2024~\cite{Negrepontis2024c}) has led us
to the fascinating, novel and revealing conclusion that

\begin{center}
  \begin{minipage}{0.8\textwidth}
  an intelligible Being is the philosophical analogue of a dyad of geometrical lines in periodic anthyphairesis.
\end{minipage}
\end{center}

Armed with this knowledge, it should not take us long to realize that
the best way to approach \textit{Parmenides} is not to start stydying it from
the beginning, but to go directly to \textit{the Second Hypothesis} spanning
the 142b1--155e3 pages. Indeed, the principal part of the \textit{Parmenides}
consists of the Introduction 126a1--137c3, the First Hypothesis
137c4--142a8, and the Second Hypothesis.

The First Hypothesis deals with a One that is literally indivisible,
with no parts, and, as expressly stated in its final section
141e3--142a8, fails to be an intelligible Being, the Second Hypothesis
with a One that is, on the contrary, definitely an intelligible Being,
and indeed one of the form \textit{Name and Logos}, as expressly stated in its
final section 155d3--e3, a description identical with the One in the
\textit{Theaetetus} and the \textit{Sophist}.

Thus, we expect that the One of the Second Hypothesis is a philosophical analogue of a geometric dyad in periodic anthyphaireis (Section 2).

The first part of the second hypothesis 142b1--143a3 is a rather
difficult text, but with the help from Proclus’ \textit{Platonic
  Theology} 3,89,21--26 and 4,94,11--16, and with careful linguistic
analysis, indeed leads to the interpretation that the One of the
second hypothesis with its part the Being form a dyad in the
philosophical analogue of infinite anthyphairesis; and indeed, the
anthyphairesis appears in its complete form, not in the abbreviated
one observed in the \textit{Sophist}’s divisions of the Angler and the \textit{Sophist}
(Section~3).

We appear to be on the right track, and we expect that anthyphairetic
periodicity will next make its appearance. However, in 143c1--144c2,
Plato instead introduces \textit{numbers}, with units the infinite sequence of
anthyphairetic remainders of the dyad One Being. However, with the
help of the \textit{Philebus} 56d--e passage, this coarse approach of relying
only on the infinity of the anthyphairesis of the philosophical fyad
One Being, eventually is understood to lead to the introduction of the
unwanted, \textit{eristic} numbers, with unequal units, and not to the wanted
\textit{dialectic} numbers of equal/equalized unis (Section 4).

The passage 144c2--d4 marks a turning point in that oneness of the arithmetical units takes a concrete strong form: 
\begin{center}
  \begin{minipage}{0.8\textwidth}
each part of the Being must exhibit its oneness by containing a suitable part of the One. 
\end{minipage}
\end{center}
But how is this to be achieved? At the position 142d4 the reader must
realize that there is no other way to understand but to discover that
Plato hides the crucial method in a seemingly irrelevant earlier
passage 138a3--7 in the one of the first hypothesis; there it is
stated that that the One will be contained in a part of the Being
\textit{cyclically} 138a3--4, and this cyclicity is in fact a
\textit{cycle of contacts} (\textit{Hapseis}) 138a5--7. Thus, we must
understand what is meant by contact (hapsis) and how contacts form a
cycle. Ontacts is in the passage 148d--149d, where “the plus one
rule”:
$$
\text{number of terms} = \text{(number of) contacts (hapseis)} +1 $$
is enunciated.

Here we come to the most delicate and essential step of our
interpretation. It is the point where \textit{Pythagorean music} comes
to Plato’s help in introducing philosophical periodicity and
dialectical number in every intelligible Being and to our help in
order to understand Plato’s description of the intelligible Being in
the Second Hypothesis of the \textit{Parmenides}.  Indeed, the plus one rule of
the Second Hypothesis is a rephrasing of the old and fundamental rule
of \textit{Pythagorean music}:
$$
\text{the number of chords is equal to the (numbers of) musical intervals}+1. 
$$
Thus, the octave consisting of eight chords/terms produce seven
musical intervals thus: tone, tone, diesis, tone, tone, tone, diesis;
and in the \textit{Timaeus}, the universe is generated by the circularity of
thirty-four terms and thirty-three musical intervals.  But a musical
interval is essentially the ratio of two successive terms/chords;
thus, thanks to music we now realize that
\begin{center}
  \begin{minipage}{0.8\textwidth}
\textit{a contact (hapsis)} is the ratio of two successive terms in the infinite sequence of the anthyphairetic remainders of the dyad One, Being,
\end{minipage}
\end{center}
and we are smoothly led to the revelation that
\begin{center}
  \begin{minipage}{0.8\textwidth}
    \textit{the circularity of contacts} (\textit{kukloi haptesthai},
    138a5--7) is exactly the \textit{Logos Criterion} for the anthyphairetic
    periodicity of the dyad One and Being (Section 5).
\end{minipage}
\end{center}
The definition of dialectical number in an intelligible Being, 
in terms of the Logos Criterion (138a3--7) and the plus one rule (148d–149d), 
the finiteness of the dialectical number of the parts of the Being,
the equality of the dialectical number of the parts of the One and of the dialectical number of the
parts of the Being, and

 the rejection of the eristic numbers (144d4--e3)

 \noindent Once anthyphairetic periodicity is achieved, the real nature of the dialectic numbers is understood in terms of the plus one rule:
 
\begin{center}
  \begin{minipage}{0.8\textwidth}
The \textit{dialectic number} of an interval of terms in the infinite sequence of anthyphairetic remainders of the dyad One, Being is equal to the (number of) contacts plus one.
\end{minipage}
\end{center}
It follows that 
\begin{center}
  \begin{minipage}{0.8\textwidth}
the dialectic number of an infinite or finite interval is always \textit{finite}, 
\end{minipage}
\end{center}
not exceeding the length of the anthyphairetic period plus one (Section 6).
 
Now we have at our disposal all the concepts and tools that make understandable the most fundamental of Plato’s notions, notions that have appeared to verge on the contradictory and which modern Platonists, without the help of the mathematics of periodic anthyphairesis, were unable to unravel. Thus:

\medskip

\noindent--An intelligible Being is described as 

\medskip

an \textit{Indivisible} Line, 

\medskip

\noindent certainly not in the sense of indivisibility of magnitudes, that would be ridiculous, 
but in the precise sense that the anthyphairetic division of an intelligible Being, being periodic, does not produce, after the completion of a period, any new contacts (hapseis), and thus no new units for number.  

Plato’s Indivisible Line, considered mysterious by modern Platonists, is identified with Plato’s Intelligible Being (Section 7).

\medskip

\noindent--Plato’s scandalous statement in the \textit{Sophist} 257d7--258c5, 258d5--e5 that

\medskip

 not Being is a Being, 

\medskip
 
\noindent constituting another mystery for the modern Platonists, is perfectly
 understandable as follows: the two initial parts, the One and the
 Being, of the intelligible Being One in the Second Hypothesis of the
 \textit{Parmenides} are \textit{formally equalized as intelligible Being}, exactly
 because the dialectical number of the parts of the One is finite and
 equal to the dialectical number of the parts of the Being (\textit{Parmenides}
 144e1--3). In the Sophist, the paradigmatical initial dyad of the
 \textit{Parmenides} \textit{One, Being}, is renamed \textit{Being, not-Being}, generalizing
 intelligible Beings with anthyphairetic dyads \textit{Beautiful,
 not-Beautiful, Just, not-Just}, and so on.
 
\medskip

\noindent-- Plato’s fundamental property is its Oneness, in fact every intelligible Being is One and Many. Intelligible Oneness is a mystery for the modern Platonists, but the Oneness of an intelligible Being is now essentially evident, since anthyphairetic periodicity implies that 
\begin{center}
  \begin{minipage}{0.8\textwidth}
    the dialectical number of parts of every anthyphairetic remainder,
    whether it be $\text{One}_k$ or $\text{Being}_k$, is finite and equal to the
    dialectical number of the parts of the One;
  \end{minipage}
\end{center}
Thus 
\begin{center}
  \begin{minipage}{0.8\textwidth}
\textit{an intelligible Being is One in a self-similar sense}: 
each of its infinite in multitude parts is equalized to the One (144e1--3, 145a2--4).
\end{minipage}
\end{center}
Plato expresses the same thing by his statement that 
\begin{center}
  \begin{minipage}{0.8\textwidth}
the infinite multitude of parts into which a power is divided is \textit{collected into one}
(ἐπειδὴ ἄπειροι τὸ πλῆθος αἱ δυνάμεις ἐφαίνοντο, πειραθῆναι \textit{συλλαβεῖν εἰς ἕν,  
Theaetetus} 147d7-8: ταύτας πολλὰς οὔσας \textit{ἑνὶ εἴδει περιέλαβες, Theaetetus} 148d5--6) 
\end{minipage}
\end{center}

\medskip

\noindent(Section 8).

\medskip

\noindent--The Compresence of the almost contradictory properties, \\
Infinite \&\ Finite, One \&\ Many, In itself \&\ in the Other, Motion \&\ Rest,\\
in Plato’s intelligible Being; and the use of these properties to separate the intelligible Beings  from the sensibles (Section 9).

It should be noted that according to the Introduction in Plato’s
\textit{Parmenides}, Zeno of Elea constructs his true Beings in such way that
they satisfy the almost contradictory simultaneous compresence of
opposite properties, such as Infinite and Finite, One and Many, and in
Motion and at Rest (128e5--130a2), and then uses these properties in
order to establish that the true Beings are different and separate
from the sensibles (127d6--128a4). This suggests that Zeno’s true
Beings are quite similar precursors to Plato’s intelligible Beings, in
particular that Zeno’s philosophical thought is anthyphairetic in
nature, and thus necessarily was influenced by the Pythagorean
Mathematics of incommensurability. This is in fact studied and
established in Negrepontis, 2023~\cite{Negrepontis2023}. Thus, Plato’s
debt to Zeno and Zeno’s debt to the Pythagoreans is considerable.

\medskip

\noindent--The Third Man Argument in the Parmenides 132a1--b2 is now
easily resolved.  Vlastos’ interpretation of the Third Man Argument,
based on his \textit{Non-Identity principle}, in the is false, exactly because
the self-similar Oneness of an intelligible Being prevents the
infinite regress. Thus, Plato’s Third Man Argument turns out to be,
not an argument menacing Plato’s theory of intelligible Beings—as
erroneously conceived by Vlastos and his followers, but simply as an
innocuous exercise testing the reader’s understanding of the nature of
intelligible Beings (Section 10).

\section{The failure of modern Platonic scholars, from the 15th century on, to come to an understanding of Plato's philosophy, especially of Plato's \textit{Parmenides}}

Here is a small sample of opinions on the status of our knowledge and understanding of Plato's \emph{Parmenides} by Platonic scholars. \\

\noindent C.~C.~Meinwald, 1991~\cite{Meinwald1991}
\begin{quote}
Plato's \emph{Parmenides} today finds itself in a strange position: it is clearly an important work, but its import remains remarkably unclear. 
The difficulty of analyzing this text is due, in part, to its complicated structure. Within three frames we find the dialogue proper, itself consisting of two parts connected by a brief transitional section. The first part of the main dialogue is a series of rather brief exchanges enlivened by humor and some dramatic incident; the second part consists of almost thirty Stephanus pages of obscure and unadorned argument. These two parts are so strikingly different that there can be no question of their being coordinate episodes of the same kind. Yet when deciding what, exactly, the character of each part is and how they fit together to make up a whole, scholars still express doubt. But understanding the second part of the dialogue has been the single most intractable task in interpreting the \emph{Parmenides}, if not in Plato scholarship as a whole. We are faced with an unbroken series of arguments---many seemingly so bad as to be embarrassing---systematically arranged to produce apparently contradictory conclusions. Interpreters are so divided about what this exercise achieves that disagreement still persists over whether it has any positive results at all.
\end{quote}

\noindent B. Frances, 1996~\cite{Frances1996} 
\begin{quote}
Of course, the reason for the neglect of the last three quarters of the dialogue is obvious enough: interpreting the rest of the dialogue has looked utterly hopeless. Despite the best efforts of the ablest scholars, most of the arguments in the latter part of the dialogue [\emph{Parmenides}] have appeared to be incomprehensible.\\
\ldots 
\begin{quote}
    For all these reasons, then, the One is and becomes older and younger than itself and the others, and neither is nor becomes older or younger than itself or the others ([\emph{Plato's Parmenides}] 155d).
 \end{quote}

\textit{The long second part of the} Parmenides \textit{includes many fantastic, apparently contradiction-ridden statements similar to the one above. The problem of formulating a convincing interpretation of the dialogue that makes sense of these statements has proven so difficult that it has simply been ignored} by most commentators on the third man argument (or arguments) found in the first part of the dialogue.
This by itself may be a serious defect of these analyses of the third man argument given by these scholars. For Plato made it clear (135d) that he thought the truth regarding the problems of the first part of the dialogue---such as that of the third man---is to be found by undertaking the intellectual `exercise' demonstrated by Parmenides in the second part of the dialogue. And it is also clear that after writing the \emph{Parmenides} Plato continued to endorse a theory of Forms.
Thus, he must have thought that the problems brought up in the first part of the dialogue were to be treated utilizing the considerations from the second---treated in such a way that the problems do not rule out a theory of Forms continuous with that presented in the earlier dialogues. 
\textit{Therefore, in order to see how Plato formulated and countered these problems one must do the dirty work; one cannot rest content, as most have done, with some interpretation that fails to locate Plato's response to the third man argument firmly in the second part of the } Parmenides. \textit{One must, under pain of misinterpretation and irrelevance, explain the significance that the paradoxical arguments of the second part have for the problems of the first part. I suspect that we have all known this is true.}
In spite of that uncomfortable knowledge, reams have been written on the first quarter of the dialogue, while the last three quarters are usually completely ignored, as if to say that Plato was not really serious when he insisted that we exercise as Parmenides demonstrated in order to see our way through the difficulties such as that of the third man. This is not to say that no good work can be done while ignoring the opaque parts of the dialogue; Vlastos' catalytic 1954 paper was a watershed despite its neglect of the latter three quarters of the dialogue.
\textit{But it does imply, despite many scholars' hopes, that in the absence of a detailed understanding of the multitude of paradoxical arguments, any claim to know what Plato thought of the third man argument is simply without any foundation. Of course, the reason for the neglect of the last three quarters of the dialogue is obvious enough: interpreting the rest of the dialogue has looked utterly hopeless. Despite the best efforts of the ablest scholars, most of the arguments in the latter part of the dialogue have appeared to be incomprehensible.}
B.~Frances, 1996~\cite{Frances1996}.
%\textit{Plato's Response to the Third Man Argument in the Paradoxical Exercise of the} Parmenides, Ancient Philosophy 16 (1996), 47--64
\end{quote}

\noindent R. E. Allen, 1997~\cite{Allen1997}
\begin{quote}

The \emph{Parmenides} has been read as everything from a joke to an exercise in the detection of fallacies to a revelation of the Unknown God. In actual fact, the dialogue is aporetic,\ldots and in no sense a budget of fallacies.\\
\ldots\\
What follows is not easy reading, a fact which no amount of hard writing could change. But if the reader will take time to follow the argument, there is a reward. Beneath the surface complexity of the \emph{Parmenides} is a simplicity, a rhythmical, dynamic beauty like that of music. I count myself fortunate to have been able to live, for brief periods of time, in the presence of such great beauty
\end{quote}

\noindent S. Scolnicov, 2003~\cite{Scolnicov2003}
\begin{quote}

Of all Plato's dialogues, the \emph{Parmenides} is notoriously the most difficult to interpret. Scholars of all periods have violently disagreed about its very aims and subject matter. The interpretations have ranged from reading the dialogue as an introduction to the whole of Platonic---and more often Neoplatonic---metaphysics [For a summary of Neoplatonic interpretations, see Dodds, 1928~\cite{Dodds1928}, Wundt, 1935~\cite{Wundt1935}. The esotericist interpretation (e.g., Migliori,1990~\cite{Migliori1990}), influenced by Krämer, can be seen as a variant of this trend. In the same vein, Séguy-Duclot, 1998~\cite{Seguy-Duclot1998} interprets the dialogue as pointing beyond itself, to higher levels, up to a henological point of view above ontology.] to viewing it as a record of unsolved (and perhaps unsolvable) ``honest perplexities," [Vlastos (1965b [1954])~\cite{Vlastos1954} 145] as protreptic ``mental gymnastics," [Grote, 1875~\cite{Grote1875}, III, chap. 27; Peck, 1953--54~\cite{Peck1953}; cf. Kutschera, 1995~\cite{Kutschera1995}. See also Wilamowitz (1948), I 402; most recently Gill, ``Introduction," in Gill and Ryan, 1996~\cite{Gill1996}. Klibansky, 1943~\cite{Klibansky1943} n. 1) attributes such a view already to Alcinous (Albinus), possibly on the strength of chaps. 5 and 6 of his \emph{Didaskalikos}.] as a collection of sophistic tricks, [E.g., Owen (1986 [1970]~\cite{Owen1986}).] or even as an elaborate (though admittedly tedious) joke. [Cf., e.g., Taylor, 1934~\cite{Taylor1934} 29.]
\end{quote}

\noindent S. Peterson, 2011~\cite{Peterson2011}
\begin{quote}

The dialogue's second part, 137c--166c, is the longest passage of unrelenting argument in Plato's writings. Its arguments are his most puzzling. p.383
\end{quote}

\noindent S. Peterson, 2019~\cite{Peterson2019} 
\begin{quote}

The \emph{Parmenides's} reputation as Plato's most challenging [S. Scolnicov, Plato's ``Parmenides" [\emph{Parmenides}] (Berkeley, 2003), 1;\\
G. Priest, ``The Parmenides. A Dialethic Interpretation" [\emph{Parmenides}], \emph{Plato Journal} 12, 2012~\cite{Priest2012}, 1–63, 1; \\
C. Kahn, \emph{Plato and the Post- Socratic Dialogue} [\emph{Plato}] (Cambridge, 2013~\cite{Kahn2013}), 1], indeed enigmatic [M. L. Gill, ``Introduction" to M. L. Gill and P. Ryan, \emph{Plato. ``Parmenides"} [\emph{Parmenides}] (Indianapolis, Ind., 1996~\cite{Gill1996}),1. Kahn \emph{Plato}, 2;\\ S. Rickless, ``Plato's Parmenides," [\emph{Parmenides}] \emph{Stanford Encyclopedia of Philosophy}, revision of July 30, 2015), 1.] work continues. \cite{Peterson2019} p.232
\end{quote}

\noindent S. Rickless, 2020 \cite{Rickless2020}
\begin{quote}

The \emph{Parmenides} is, quite possibly, the most enigmatic of Plato's dialogues.
\end{quote}

\section{The basics of periodic anthyphairesis: Our study of Plato's philosophy, especially of \emph{Meno, Theaetetus, Sophist, Statesman}, leads us to expect that the dyad One and Being of the Second Hypothesis in the \emph{Parmenides} satisfies the philosophical analogue of periodic anthyphairesis}

\subsection{The One of the second hypothesis in the \emph{Parmenides} is expected to satisfy with its part the Being the philosophical analogue of periodic anthyphairesis}

In earlier works we have found that \textit{Plato's intelligible Being is the philosophical imitation of a dyad in periodic anthyphairesis}. In fact, the knowledge (episteme) of the intelligible Being Angler and of the intelligible Being Sophist in Plato's dialogue \emph{Sophist}, and of the intelligible Being Statesman in Plato's dialogue \emph{Statesman} (Negrepontis, 2018~\cite{Negrepontis2018}, Negrepontis, Farmaki, Brokou, 2024~\cite{Negrepontis2024c}), 
and of the diameter to the side of a square in Plato's dialogue \emph{Meno} (Negrepontis, 2024~\cite{Negrepontis2024a}) have been found to be a dyad satisfying (the philosophical analogue of) of Theaetetus' Logos Criterion for periodic anthyphairesis. 

Since the One of the second hypothesis in the \emph{Parmenides} is declared to be (in 155d--e) a paradigmatical intelligible Being, 
\textit{we certainl}y expect that the One of the Second Hypothesis in the \emph{Parmenides} possesses the philosophic analogue of periodic anthyphairesis.

\subsection{The concept of anthyphairesis}

Let $a, b$ be two natural numbers or magnitudes (e.g. line segments), with $a>b$. 

We measure the great $a$ with the small $b$, and we find how many small $b$'s are needed to cover the large $a$. 
There are then a natural number $k_1$, the quotient, and $c_1$, the remainder, such that
$$a=k_0b+c_1, \quad \text{with } c_1<b.$$
This is the first step of the anhyphairesis of $a$ to $b$. If $c_1$ is zero, this is the end of the process. Otherwise a new pair of great and small is formed, namely the pair $b>c_1$. Notice that the role of $b$ is reversed: in the first step $b$ is small, in the second step $b$ is great.

We repeat the process. There are a natural number $k_2$ and $c_2$, such that
$$b=k_1c_1+c_2, \quad \text{with } c_2<c_1.$$
We continue in this way:
$$c_1=k_2c_2+c_3, \quad \text{with } c_3<c_2,$$
$$\ldots$$
$$c_n=k_{n+1} c_{n+1}+c_{n+2}, \quad \text{with } c_{n+2}<c_{n+1}$$

\subsection{Finite and infinite anthyphairesis}

\begin{itemize}
    \item If this process ends at some point (namely for some index $n$, $c_n$ is zero, hence $c_{n-2}=k_n c_{n-1}$), then we say that the anthyphairesis of $a$ to $b$ is \textit{finite}. This is always the case if $a, b$ are natural numbers, and the last non-zero \textit{remainder} $c_{n-1}$ is \textit{the greatest common-measure-divisor} of $a$ and $b$ (Propositions VII.1 \& 2 of the \emph{Elements}).\\ 
%    X.3 για σύμμετρα μεγέθη.
    \item If this process does not end, then we say that the anthyphairesis of $a$ to $b$ is \textit{infinite}. This can happen only if $a$ and $b$ are magnitudes. It is then proved that $a$ and $b$ do not have any common measure, that $a$ and $b$ are \textit{incommensurable}:
\end{itemize}

\begin{quote}
Proposition X.2 of the \emph{Elements}. If $a,b$ are magninudes, $a>b$, and the anthyphairesis of $a$ to $b$ is infinite, then $a$ and $b$ are incommensurable.
\end{quote}

Note that in this case there is an infinite, strictly decreasing sequence of remainders:
$$a>b>c_1>c_2>c_3>\ldots>c_n>c_{n+1}>\ldots .$$

\subsection{Periodic anthyphairesis}

The most important case of infinite anthyphairesis is the \textit{periodic} one. An anthyphairesis is \textit{periodic} if the sequence of \textit{quotients}
$$k_0, k_1, k_2, \ldots, k_n, k_{n+1},\ldots$$
of natural numbers is a periodic sequence of numbers.

An anthyphairesis, not necessarily infinite a priori, is recognized as periodic according to the following \textit{Logos Criterion}: there is an index $n$ such that
$$\frac{a}{b}=\frac{c_n}{c_{n+1}};$$
Indeed then $\text{Anth}(a,b)=\text{Anth}(c_n,c_{n+1})$, readily implying
$$\frac{b}{c_1}=\frac{c_{n+1}}{c_{n+2}}, \quad \frac{c_1}{c_2}=\frac{c_{n+2}}{c_{n+3}}, \quad \ldots,$$
and hence that the ratios of successive parts-remainders
$$\frac{a}{b}, \frac{b}{c_1}, \frac{c_1}{c_2}, \frac{c_2}{c_3},\ldots, \frac{c_n}{c_{n+1}}, \ldots$$
is a periodic sequence of ratios.\\
The first theory of proportion for magnitudes, including the Logos Criterion, was developed by Theaetetus, and certainly this is the one used by Plato (cf. Negrepontis \&\ Protopapas, 2025~\cite{Negrepontis2025b}).

\textit{Example}. The anthyphairesis of the diameter $a$ to the side $b$ of a square (so that $a^2=2b^2$) has sequence of quotients the sequence $[1,2,2,2,\ldots]$, periodic from the second stage on.

The modern concept, essentially equivalent to the ancient anthyphairesis, is the \textit{continued fraction of a (positive) real number}. The theory was developed during the 16th, 17th (Fermat, Wallis), 18th (Lagrange, Euler), and 19th (Galois, Gauss) century.

\section{The One of the second hypothesis in the \emph{Parmenides}, an intelligible Being according to the \emph{Parmenides} 155d, whose knowledge is given in terms of Name plus Logos, forms, together with its part the Being, the dyad One and Being, satisfying the philosophical analogue of \emph{infinite anthyphairesis}, according to Plato's \emph{Parmenides} 142b1--143a3 and Proclus' \emph{Platonic Theology} 3, 89, 21--26 and 4, 94, 11--16}

\textit{Heuristics}. The One of the second hypothesis in Plato's \emph{Parmenides} (142b--155e) is a paradigmatical \textit{intelligible Being} (as stated in 155d), and one whose knowledge is described as \textit{Name plus Logos} (as stated in 155d8--e1). Since we have already shown in earlier works (Negrepontis 2012~\cite{Negrepontis2012}, 2018~\cite{Negrepontis2018}, 2024~\cite{Negrepontis2024a}) that the knowledge of the examples of intelligible Beings, the Angler and the Sophist in the \emph{Sophist}, given in terms of \textit{Name plus Logos}, both have a structure that satisfy the philosophical analogue of a geometric dyad satisfying \textit{the Logos Criterion} and thus \textit{in periodic anthyphairesis}, we expect that the same will be true for the One of the second hypothesis.

The first part of the second hypothesis in the \emph{Parmenides} (142b1--143a3) is devoted in proving that the One together with a part of it, the Being, forms a dyad that satisfies \textit{the philosophical analogue of infinite anthyphairesis}. We obtain this novel interpretation, for the first two steps of the anthyphairesis, first, by a careful linguistic analysis of the passage (142c7--e5; sections 3.3, 3.4), and, secondly, by taking into account the brief but revealing comments on two passages by Proclus in \emph{Platonic Theology} 3, 89, 21--26 (section 3.5); and for the general step of the infinite anthyphairesis, first by an analysis of the passage, in sections 3.6,3.7, and secondly by taking into account the brief but revealing comments on two passages by Proclus in \emph{Platonic Theology 4, 94, 11--16 (section 3.8)}.

We note the striking similarity that exists between Plato taking an upward step to the Whole with two geometric proofs of incommensurability (a) of the diameter to the side of a square (a modern proof with help from the Elegant Theorem in Proclus' Commentary to Plato's \emph{Republic}) and (b) of the mean extreme ratio (Proposition XIII.5 of the \emph{Elements}). We finally establish a remarkably close relation of the infinite multitude of parts of the dyad One Being in the passage 142b1--143b3 with the infinite multitude of parts in Zeno's \textit{Fragment B3}.

\subsection{``The One is" (hen estin) is the hypothesis that the One is an intelligible Being 142b1--5}

\begin{quote}\selectlanguage{polutonikogreek}
—Βούλει οὖν ἐπὶ τὴν ὑπόθεσιν πάλιν ἐξ ἀρχῆς ἐπανέλθωμεν, ἐάν τι ἡμῖν ἐπανιοῦσιν ἀλλοῖον φανῇ; \\
— Πάνυ μὲν οὖν βούλομαι. \\
— Οὐκοῦν ἓν εἰ ἔστιν, φαμέν, τὰ συμβαίνοντα περὶ αὐτοῦ, ποῖά ποτε τυγχάνει ὄντα, διομολογητέα ταῦτα· οὐχ οὕτω; \\
— Ναί. \\
\selectlanguage{english}Shall we then return to our hypothesis (hupothesin) and see if a review of our argument discloses any different (alloion) point of view? By all means. We say, then, that ``if the one is" (hen ei estin), we must come to an agreement about the consequences (ta sumbainonta), whatever they may be, do we not? Yes. 142b1--5
\end{quote}

The meaning of the hypothesis ``hen estin" is that \textit{the One is (a true intelligible Being)}. According to the Eleatic philosopher Parmenides the true, intelligible Being is a \textit{One}. But which One, a One in what sense and with which properties? In the first hypothesis of the \emph{Parmenides} a partless One, a One similar to a geometric point, was assumed, but according to the conclusion of the first hypothesis such a One definitely failed to be an intelligible Being. So, in the second hypothesis Plato introduces a second, non-partless One, radically different with the One of the first hypothesis, that will hopefully succeed this time to be seen as the true intelligible Being---and which in fact does succeed, as clearly stated at the conclusion of the second hypothesis (155d).

\subsection{The One of the second hypothesis forms a dyad of parts with the Being, in such a way that the One \textit{participates} into the Being 142b5--7, c5--7}

\begin{quote}
— Ὅρα δὴ ἐξ ἀρχῆς. \textit{ἓν εἰ} ἔστιν, ἆρα οἷόν τε αὐτὸ εἶναι μέν, \textit{οὐσίας} δὲ μὴ \textit{μετέχειν}; \\
— Οὐχ οἷόν τε. 142b5--7 \\
Now consider the first point. If the One is, can it be and not \textit{ participate} (metechein) in the Being (ousias)? \\
No, it cannot.

— Ἆρα οὖν ἄλλο ἢ ὅτι \textit{οὐσίας μετέχει τὸ ἕν}, τοῦτ' ἂν εἴη τὸ λεγόμενον, ἐπειδάν τις συλλήβδην εἴπῃ ὅτι \textit{ἓν ἔστιν}; \\
— Πάνυ γε. 142c5--7 \\
Then if we say concisely (sullebden) `the One is,' it is equivalent to saying that `\textit{the One participates in Being} (to hen metechei ousias)'? \\
Certainly.
\end{quote}

The hypothesis `the One is' is interpreted as implying that there are two initial parts in the intelligible Being, 
the One and the Being, such that `the One \textit{participates }(metechei) in the Being'\footnote{ ‘Ousia’ in the Parmenides
has always the meaning of the ‘Being part’.\par In  142b6, 7, 8 ousia is first
introduced in the second deduction, and it is stated that \par ‘hen’
participates in ousia. \par In the next occurrences \par 142d3 ousia \par 143
a4,6, b1,2,3,4,5,7, c1,2,5,7, \par 144 a7,8,  b1--e 3 \par is clearly a part in
the One Being.\par Throughout the whole passage ‘ousia’ refers to the Being
part. \par  (cf.  the following list of references to ‘ousia’ during this
passage:\par Επὶ πάντα ἄρα πολλὰ ὄντα ἡ οὐσία νενέμηται 144b1\par οὐσία γε τῶν
ὄντων του ἀποστατοῖ; 144b4\par  μέρη ἀπέραντα τῆς οὐσίας 144c1\par  μέρη αὐτῆς
144c2 \par ἔστι μὲν μέρος τῆς οὐσίας, οὐδὲν μέντοι μέρος; 144c3,\par Πρὸς
ἅπαντι ἄρα [ἑκάστῳ] τῷ τῆς οὐσίας μέρει 144c6,\par οὔτε σμικροτέρου οὔτε
μείζονος μέρους οὔτε ἄλλου οὐδενός 144c7--8 \par ἅπασι τοῖς τῆς οὐσίας μέρεσιν
144d3\par πολλὴ ἀνάγκη εἶναι τοσαῦτα ὅσαπερ μέρη 144d5\par πλεῖστα μέρη ἡ οὐσία
144d6\par \par This can be seen because the final conclusion 144e1--3 is clearly
about the Being part, and the conclusion follows from all the previous numerous
statements about ‘ousia’ in the passage 144b1--e3.}. 

\textit{Question 3.2.} What does Plato mean by the statement (in 142b5--7, c5--7) that `the One \textit{participates} into the Being'?

The basic use of participation in Plato is that ``a sensible entity \textit{participates} in the intelligible Being" (cf. \emph{Politeia} 476d, \emph{Phaedo} 100c4--8), and through this participation the sensible acquires approximately, imperfectly some similarity to the intelligible Being. We thus expect that if $a$ \textit{participates} in $b$, then $a$ becomes somewhat, approximately similar to $b$. Furthermore, we expect that the participated entity $b$ is the active one, and the participating entity $a$ the passive one. Here of course both parts One and Being are intelligible, but it is reasonable to expect that participation will have a similar meaning. Our \textit{Answer}, based on Plato's text and Proclus' revealing comments is given in Section 3.5.

\subsection{The One, by participating into the Being, generates parts, in fact two parts, a Being and a One, and becomes a whole consisting of two parts, a One and a Being 142c7--e3}

\subsubsection{The fact that `the One participates in the Being' implies that `the One has parts'}
\begin{quote}
— Πάλιν δὴ λέγωμεν,\\
\textit{ἓν εἰ ἔστιν},\\
τί συμβήσεται. σκόπει οὖν εἰ οὐκ ἀνάγκη ταύτην τὴν ὑπόθεσιν τοιοῦτον ὂν \textit{τὸ ἓν} σημαίνειν, οἷον \textit{μέρη ἔχειν}; \\
— Πῶς; 142c7--d1 \\
``Let us again say what will follow\\
`if the One is'\\
and consider whether this hypothesis must not necessarily show that the One is of such a nature as to have parts (mere echein).''\\
``How?"
\end{quote}

We have already interpreted the hypothesis ``the One is" by the statement ``the One participates into the Being"; now it is claimed that the hypothesis ``the One is" implies that ``the One has parts". Thus, the generation of parts of the One must be closely related to the participation of the One into the Being. The process by which the One generates parts is not yet clear to Socrates, cf. the question how (pos, 142d1)

\textit{Question}\\
How the \textit{participation} of the One into the Being is a cause for the generation of parts of the One? Our \textit{Answer} to be given in Sections 3.5, 3.8.

\subsubsection{A first suggestion of how the One has parts is given by going upwards and considering the way the whole One \& Being has generated parts the One and the Being 142d1--e3}

Parmenides sets out to explain how (hode), 142d1).
\begin{quote}
— Ὧδε· εἰ τὸ \textit{ἔστι τοῦ ἑνὸς ὄντος λέγεται καὶ τὸ ἓν τοῦ ὄντος ἑνός},\\
ἔστι δὲ \textit{οὐ τὸ αὐτὸ} ἥ τε οὐσία καὶ τὸ ἕν,\\
\textit{τοῦ αὐτοῦ δὲ ἐκείνου} οὗ ὑπεθέμεθα, \textit{τοῦ ἑνὸς ὄντος,}\\
ἆρα οὐκ ἀνάγκη\\
\textit{τὸ μὲν ὅλον ἓν ὂ}ν εἶναι αὐτό, \\
\textit{τούτου δὲ γίγνεσθαι (142d5) μόρια τό τε ἓν καὶ τὸ εἶναι;} \\
— Ἀνάγκη. \\
— Πότερον οὖν ἑκάτερον \textit{τῶν μορίων} τούτων \textit{μόριο}ν μόνον προσεροῦμεν, ἢ \textit{τοῦ ὅλου μόριον} τό γε μόριον προσρητέον; \\
— \textit{Τοῦ ὅλου}. \\
``How does that come about?" ``In this way. If the Being is part of the One Being and the One is part of the Being One, and the Being and the One are \textit{not the same }(ou tauton), but belong to the \textit{One Being }of our hypothesis, must not the One Being be a \textit{whole (holon)}, of which the One and the Being are \textit{parts (moria)} \textit{generated} (gignesthai, 142d5)?" ``Inevitably." ``And shall we call each of these \textit{parts (moria) }merely a \textit{part (morion)}, or must it, in so far as it is a \textit{part (morion)}, be called a \textit{part (morion)} of the whole?" ``A \textit{part (morion)} of the whole." 142b1--8
\end{quote}

Towards elucidating the manner in which the One, by participating in the Being, gets to have parts, Plato forms the \textit{whole (holon) }One \& Being (hen on), and notes that the One \& Being is a whole and has parts the One and the Being. It is important to note that the One and the Being are regarded as generated parts. But the nature of generation of new parts from existing ones is not clear, so the question takes the following form:

\subsubsection{Analogously to the whole One \& Being, the One has a part and is a whole, namely it consists of two parts 142d8--9}

\begin{quote}
— Καὶ \textit{ὅλον} ἄρα ἐστί, \textit{ὃ ἂν ἓν ᾖ}, καὶ \textit{μόριον} ἔχει. \\
— Πάνυ γε. 142d8--9 \\
Hence (ara) that which is One both is a whole (holon) and has a part (morion). Certainly.
\end{quote}

The part One possesses a part and is a whole. It appears rather clear that the part that the One possesses is the other initial part, the Being. On the other hand, to be a whole, in analogy to the upward situation (3.4.2), means that the One will consist of two parts. It is a question of how this will happen, but shortly below a more specific claim occurs, hence a more specific question is formed.

\textbf{Question 3.3.2 \& 3.} How does the upwards step (142d1--e3, d8--9) help understand the method by which the parts of the One are generated? Our Answer to be given in Sections 3.4--5; see also Section 3.9 for its mathematical origin.

\subsubsection{In fact, the One has as a part the Being; and, furthermore, the Being has a part that is a One 142d9--e3}

\begin{quote}
— Τί οὖν; τῶν μορίων ἑκάτερον τούτων τοῦ ἑνὸς ὄντος, τό τε ἓν καὶ τὸ ὄν, ἆρα \textit{ἀπολείπεσθον} ἢ τὸ ἓν τοῦ εἶναι μορίου ἢ τὸ ὂν τοῦ ἑνὸς μορίου; \\
— \textit{Οὐκ} ἂν εἴη. 142d9--e3 \\
Well then, can either of these two parts of the One \& Being (henos ontos), the One and the Being, be defective (apoleipesthon) with respect to the other, and the One not be a part (moriou) of the Being or the Being not be a part (moriou) of the One? No.
\end{quote}

The first claim is that the One contains as part the Being, something already suggested previously (cf. 3.3.1), but the second claim that\textit{ the Being has as part a One} is as yet by no means clear, in fact something of a mystery. The One that will be a part of the Being clearly cannot be the original One, since the Being is a proper part of the One.

\textbf{Question 3.3.4.} How does it happen that the part Being contains a One as part, as claimed in 142d9--e3?\\
Our \textit{Answer} to be given in Sections 3.5 and 3.6 (the crucial Proclus' comments).

\subsubsection{As a consequence of 3.3.4, each of the One and the Being possess as parts a One and a Being}

\begin{quote}
— Πάλιν ἄρα καὶ τῶν μορίων ἑκάτερον τό τε ἓν ἴσχει καὶ τὸ ὄν, 142e3--4 \\
Therefore (ara) again, each of these two parts possesses (ischei) the One and the Being
\end{quote}

\textbf{Question 3.3.5.} In what way the fact that the One has as part the Being implies (as suggested by `ara', 142e3) that the One consists of a Being part and a One part, and similarly for the Being? \\
Our Answer is given in Sections 3.4 and fully confirmed in 3.5.

\subsection{Our anthyphairetic interpretation of the crucial passage \emph{Parmenides} 142d9--3: the One contains as part the Being, therefore (ara) the One consists of the Being and the One$_1$, with One$_1$ a \textit{generated} part smaller than the Being, and the Being contains part One$_1$, therefore (ara) the Being consists of the One$_1$ and the Being$_1$, with Being$_1$ a \textit{generated} part smaller than the One$_1$ 142e3--5}

We reach this interpretation by a careful linguistic analysis of the text.

\subsubsection{The crucial passage 142e4--5}

καὶ \textit{γίγνεται τὸ ἐλάχιστον ἐκ} δυοῖν αὖ μορίοιν \textit{τὸ μόριον}' 142e4--5

This is the crucial passage, in which Plato, in a condensed manner, describes the fundamental process of generation by the dyad One and Being; for the correct interpretation of the passage, it is essential to render correctly the meaning of the two terms

\textit{`to elachiston'} and \textit{`gignetai ek'}

\subsubsection{The interpretation of the crucial terms `to elachiston' and `gignetai ek' in 142e3--5 by previous scholars}

Here is a sample of the way these have been rendered by previous scholars:
\begin{itemize}
    \item[] `and \textit{the smallest of parts} is \textit{composed} of these two parts' H.N. Fowler, 1926 \cite{Fowler1926} p.255;
    \item[] `and the part\textit{ turns out to come from at least} two parts' Meinwald, 1991~\cite{Meinwald1991} p.109 [where it is noted (in p.181) that `to elachiston' is taken adverbially, `at least', by Jowett, Cornford, Dies, Curd];
    \item[] `and \textit{the least part} also turns out \textit{to consist} of two parts' Allen, 1997~\cite{Allen1997} p.26:
    \item[] `and the part comes to be \textit{composed of at least }two parts again' Palmer, 1999~\cite{Palmer1999} p.222.
\end{itemize}
We note that `to elachiston' is rendered by Jowett, Cornford, Dies, Curd, Meinwald, and Palmer as an \textit{adverb to `gignetai ek' `at least} two parts', and by Fowler, and Allen as an \textit{adjective} to `morion', \textit{the least} part; while `gignetai ek' is rendered by Fowler, Allen, and Palmer as `composed/consist of', and by Meinwald as `to come from'.

\subsubsection{Our interpretation of the terms `to elachiston' and `gignetai ek' in 142e3--5}

We examine separately four of the terms in the 142e3--5 passage.

\textbf{3.4.3.(i) `palin', 142e3}

The term \textit{`palin', `again'}, indicates repetition. It is clear that the term confirms the analogy role of the One \& Being. Thus, the part One \textit{again (palin)} will consist of two parts, of a Being part and of a One part (and similarly for the part Being), as was the case with the One \& Being part.

\textbf{3.4.3.(ii) `ara', 142e3}

The term `therefore'(`ara', 142e3) indicates an implication: The second statement, that the One part consists of two parts, one of the type Being and another of the type One, \textit{follows, is a consequence }from the previous statement (3.3.4), that the One part possesses as part the Being part, and the Being part possesses a One part. We postpone discussing a similar claim about the Being part. Note that the One with the Being part it contains forms a dyad; thus what is needed is to obtain somehow a second part of the One complementary to the Being part.

\textbf{3.4.3.(iii) The meaning of \textit{`to elachiston' }in 142e4--5}

Meinwald, 1991~\cite{Meinwald1991} and Palmer, 1999~\cite{Palmer1999} attach `to elachiston' as an \textit{adverb} to `ek duoin morioin' rendering the phrase as `\textit{from at least two parts}'. But such a rendering is in contrast with Plato's use of the word `to elachiston' in the totality of his work; there is, besides the present one, a total of nine occurrences of `to elachiston' in all of Plato; in all of them the meaning is some `least' thing, not the adverbial one.

\begin{itemize}
    \item[{[1]}] ΣΩ. Ἐλάχιστον τοίνυν μοι δοκοῦσι τῶν ἐν τῇ πόλει δύνασθαι οἱ ῥήτορες.\\
    Socrates Then, to my thinking, the orators have the \textit{smallest (elachiston)} power of all who are in their city. \emph{Gorgias} 466b9--10
    
    \item[{[2]}] ἐὰν δὲ τύχῃ πάντων ἀσθενέστατος ὤν, πάντων \textit{ἐλάχιστον }τῷ βελτίστῳ, ὦ Καλλίκλεις;\\
    Or if he chance to be the weakest of all, ought he not to get the \textit{smallest (elachiston)} share of all though he be the best, Callicles? \emph{Gorgias} 490c5--7
    
    \item [{[3]}] [τόδε δέ μοι εἰπέ, σὺ αὐτὸς πόθεν πλεῖστον ἀργύριον ἠργάσω τῶν πόλεων εἰς ἃς ἀφικνῇ; ἢ δῆλον ὅτι ἐκ Λακεδαίμονος, οἷπερ καὶ πλειστάκις ἀφῖξαι; ΙΠ. Οὐ μὰ τὸν Δῖα, ὦ Σώκρατες.\\
    ΣΩ. Πῶς φῄς; ἀλλ' \textit{ἐλάχιστον}; \\
    ΙΠ. Οὐδὲν μὲν οὖν τὸ παράπαν πώποτε.\\
    Socrates] and the test of this is, who makes the \textit{most (pleiston)} money. Well, so much for that. But tell me this: at which of the cities that you go to did you make the most money? Or are we to take it that it was at Lacedaemon, where your visits have been most frequent?\\ 
    Hippias No, by Zeus, it was not, Socrates.\\
    Socrates What's that you say? But did you make \textit{least (elachiston)} there? \\Hippias \textit{major} [sp] 283b4--c1
    
    \item [{[4]}] ἀλλὰ ἕν γε ἀνθ' ἑνὸς οὐκ \textit{ἐλάχιστον }ἔγωγε θείην ἂν εἰς τοῦτο ἀνδρὶ νοῦν ἔχοντι, ὦ Σώκρατες, πλοῦτον χρησιμώτατον εἶναι\\
    so to depart in fear to that other world---to this result the possession of property contributes not the \textit{least (elachiston)}. \emph{Politeia} 331b5--7
    
    \item [{[5]}] Ἀπορρᾳθυμεῖν ἡμῖν δοκεῖς, ἔφη, καὶ εἶδος ὅλον οὐ τὸ \textit{ἐλάχιστον} ἐκκλέπτειν τοῦ λόγου ἵνα μὴ διέλθῃς, καὶ λήσειν οἰηθῆναι εἰπὼν αὐτὸ φαύλως, ὡς ἄρα περὶ γυναικῶν τε καὶ παίδων παντὶ δῆλον ὅτι κοινὰ τὰ φίλων ἔσται. and are trying to cheat us out of a whole division, and that not the \textit{least (to elachiston)}, of the argument to avoid the trouble of expounding it, and expect to `get away with it' by observing thus lightly that, of course, in respect to women and children it is obvious to everybody that the possessions of friends will be in common. \emph{Politeia} 449c2--5
    
    \item [{[6]}] Ἐν τοῖς γάμοις τοίνυν καὶ παιδοποιίαις ἔοικε τὸ ὀρθὸν τοῦτο γίγνεσθαι οὐκ \textit{ἐλάχιστον}. Ἐφάνησαν δὲ \textit{πλεῖστον} ἀφεστῶσαι οὐχ αἱ ἐρωτικαί τε καὶ τυραννικαὶ ἐπιθυμίαι; Πολύ γε.\\
    ``In our marriages, then, and the procreation of children, it seems there will be not \textit{the least (elachiston)} need of this kind of `right.'" \emph{Politeia} 459d4--5
    
    \item [{[7]}] \textit{Ἐλάχιστον} δὲ αἱ βασιλικαί τε καὶ κόσμιαι; Ναί.\\
    Then the tyrant's place, I think, will be fixed at the furthest remove from true and proper pleasure, and the king's at \textit{the least (elachiston)} ``\\
    ``Necessarily." \emph{Politeia} 587a13--b4
    
    \item [{[8]}] ΞΕ. Ὅταν εἴπῃ τις· ``ἄνθρωπος μανθάνει," λόγον εἶναι φῂς τοῦτον \textit{ἐλάχιστόν} τε καὶ πρῶτον;\\
   \textit{ Stranger} When one says ``a man learns," you agree that this is \textit{the least (elachiston)}  and first of sentences, do you not? \emph{Sophistes} 262c9--10
    
    \item [{[9]}] παιδὸς τυγχάνων οὐκ ἀπέσχετο \textit{τῶν μεγίστων }κακῶν. δίκη δὴ τούτῳ θάνατος, \textit{ἐλάχιστον} τῶν κακῶν, τοὺς δὲ ἄλλους\\
    For him the penalty is death, \textit{the least (elachiston)}  of evils \emph{Nomoi} 854e6--7
    
\end{itemize}

\textbf{3.4.3. (iv)} The meaning of \textit{`gignetai ek' in 142e4--5 is `is generated from'}

The general form of the sentence in 142e3--5, consisting of two parts \\
\begin{center}
{%\renewcommand{\arraystretch}{2.8}
\begin{supertabular}{|p{\cellwidth}|p{\cellwidth}|p{\cellwidth}|p{\cellwidth}|}
\hline
\cellcenter{{\itshape ek}+ }&
\cellcenter{entities A, B\\ in genitive (or dual)\\ case}
&
\cellcenter{\textit{gignetai}+} &
\cellcenter{entity C\\ in accusative\\ case}
\\[1.5ex] \hline
\end{supertabular}}
\end{center}
 is always rendered as (*)

\begin{center}
\tablefirsthead{}
\tablehead{}
\tabletail{}
\tablelasttail{}
{%\renewcommand{\arraystretch}{2}
\begin{supertabular}{|p{\cellwidth}|p{\cellwidth}|p{\cellwidth}|p{\cellwidth}|}
\hline
\cc{\textit{from}} & & & \\[.6ex] \hline
 & \cc{the pre-existing\\ entities A, B,…,} & & \\[.6ex] \hline
 & & & \cc{the entity C}\\[.6ex] \hline
 & & \cc{is \textit{generated'}} & \\[.6ex] \hline
\end{supertabular}}
\end{center}
Some kinds of generation may well be by composition, and in such cases it is
also correct to say that  (**) 

\begin{center}
\tablefirsthead{}
\tablehead{}
\tabletail{}
\tablelasttail{}
\begin{supertabular}{|p{\cellwidth}|p{\cellwidth}|p{\cellwidth}|p{\cellwidth}|}
\hline
 &
 &
 &
\arraybslash \cc{entity C}\\\hline
 &
 &
\cc{\textit{is composed}} &
\\\hline
 \cc{of} &
 &
 &
\\\hline
 &
\cc{entities A, B} &
 &
\\\hline
\end{supertabular}
\end{center}

but there are other ways of generation, and so to render (*) by (**) will then
be a mistake. Here are some examples:

\medskip

 [1]  πνεῦμα \textit{ἐξ αὐτοῦ} [αέρος] \textit{γίγνεται} ῥέοντος \textit{Cratylus} 410b3

\begin{center}
\tablefirsthead{}
\tablehead{}
\tabletail{}
\tablelasttail{}
  \begin{supertabular}{|p{\cellwidth}|p{\cellwidth}|p{\cellwidth}|p{\cellwidth}|}
\hline
 &
 &
 &
\cc{wind (pneuma)}\\ \hline
 &
 &
\cc{\itshape is generated\\ (gignetai)} &
\\\hline
\cc{\textit{from}} &
 &
 &
\\\hline
 &
\cc{the flow of air\\ (aeros rheontos)} &
 &
\\ \hline
\end{supertabular}
\end{center}
The air is preexisting, and the wind is generated from [not composed of]  air.

\medskip

 [2] \textit{γίγνεται} πάντα, οὐκ ἄλλοθεν ἢ \textit{ἐκ τῶν ἐναντίων} τὰ ἐναντία 

 \textit{Phaedo} 70e1--2; also \textit{Phaedo} 71a9, 71c6.

 \null
 
 \begin{center}
%\tablefirsthead{}
%\tablehead{}
%\tabletail{}
%\tablelasttail{}
  \begin{supertabular}{|p{\cellwidth}|p{\cellwidth}|p{\cellwidth}|p{\cellwidth}|}
\hline
 &
 &
 &
\cc{\itshape all contraries\\ (enantia)}\\\hline
 &
 &
\cc{\textit{are generated}\\ \textit{(gignetai)}} &
\\\hline
\cc{from} &
 &
 &
\\\hline
 &
\cc{their contraries only} &
 &
\\\hline
\end{supertabular}
\end{center}
The opposites are pre-existing, and all things are generated from [not composed
of] them. For example, the number 100 is generated from the large number 103
and the small number 3, but the number 100 is not composed of 103 and 3; and a
child is generated from its parents but is not composed of them.

\medskip

[3] \textit{ἐξ ὧν} μάλιστα ταῖς πόλεσιν καὶ ἰδίᾳ καὶ δημοσίᾳ κακὰ \textit{γίγνεται Republic} 373e7

\begin{center}
\tablefirsthead{}
\tablehead{}
\tabletail{}
\tablelasttail{}
  \begin{supertabular}{|p{\cellwidth}|p{\cellwidth}|p{\cellwidth}|p{\cellwidth}|}
\hline
\cc{\textit{From (ex)}} &
 &
 &
\\\hline
 &
\cc{those things} &
 &
\\\hline
 &
 &
 &
\cc{the greatest disasters\\
in the cities\\
public and private,}
\\ \hline
 &
 &
\cc{\itshape are generated\\ (gignetai)} &
\\ \hline
\end{supertabular}
\end{center}
The pre-existing things from which disasters are generated.

\medskip

[4] \textit{ἐκ} δὲ τῆς τούτων ὁμιλίας τε καὶ τρίψεως πρὸς ἄλληλα \textit{γίγνεται ἔκγονα} πλήθει
μὲν ἄπειρα, δίδυμα δέ, τὸ μὲν αἰσθητόν, τὸ δὲ αἴσθησις \textit{Theaeteteus} 156a7--b2

\begin{center}
\tablefirsthead{}
\tablehead{}
\tabletail{}
\tablelasttail{}
\begin{supertabular}{|p{\cellwidth}|p{\cellwidth}|p{\cellwidth}|p{\cellwidth}|}
\hline
\cc{\textit{from (ek)}} &
 &
 &
\\ \hline
 &
 \cc{the interaction and\\ friction of two\\ opposite powers,\\
 an active\\
 and a passive,}
&
 &
\\ \hline
 &
 &
 &
\cc{an infinite multitude\\ of
  twin offsrprings\\ (ekgona):\\
  the feeling of sense\\ (aistheton) and the\\ sensible
(aesthesis);}
\\ \hline
 \end{supertabular}
\end{center}
From the interaction of two pre-existing opposite powers, an active and a
passive, an entity is generated, for example subtracting the number 7 from the
number 5 the number 2 is generated. We certainly cannot say that the number 2
is \textit{composed} of the numbers 7 and 5, but we can say that the number 2 is
\textit{generated} from the passive number 7 and the active number 5.

\medskip

[5] τρέφεται καὶ γίγνεται ἐκ τούτου  καὶ αὔξεται τὸ τοῦ παντὸς πῦρ \textit{ὑπὸ τοῦ παρ'
ἡμῖν πυρός},  ἢ τοὐναντίον \textit{ὑπ' ἐκείνου} ταό τ' ἐμὸν καὶ ταὸ σὸν καὶ ταὸ ταῶν
ἄλλων ζῴων  ἅπαντ' ἴσχει ταῦτα \textit{Philebus} 29c5--8

\begin{center}
\tablefirsthead{}
\tablehead{}
\tabletail{}
\tablelasttail{}
\begin{supertabular}{|p{\cellwidth}|p{\cellwidth}|p{\cellwidth}|p{\cellwidth}|}
\hline
 &
 &
 &
\cc{the fire\\ of the universe}\\ \hline
 &
 &
\cc{is nourished,\\ \textit{generated (gignetai),}\\ 
 and increased}
&
\\ \hline
\cc{\itshape from (ek)} &
 &
 &
\\\hline
 &
\cc{the fire\\
within us,} &
 &
\\\hline
\end{supertabular}
\end{center}

\noindent or, on the contrary,

\begin{center}
\tablefirsthead{}
\tablehead{}
\tabletail{}
\tablelasttail{}
\begin{supertabular}{|p{\cellwidth}|p{\cellwidth}|p{\cellwidth}|p{\cellwidth}|}
\hline
 &
 &
 &
\cc{my fire, and yours,\\ and that\\ of all living
beings,}
\\ \hline
 &
 &
 \cc{is nourished,\\ \textit{generated (gignetai),}\\
   and increased}
&
\\ \hline
\cc{\itshape from} &
 &
 &
\\\hline
 &
\cc{the universal fire} &
 &
\\\hline
\end{supertabular}
\end{center}
Fire itself, intelligible fire, is pre-existing, and from it the sensible fire
is generated; 

\noindent it would not be true to state that 

\noindent sensible fire is `\textit{composed of}' intelligible fire, 

\noindent but it would be perfectly true to state that 

\noindent sensible fire is `\textit{generated from}' intelligible fire. 

\selectlanguage{polutonikogreek}

\medskip

[6]'τά τε αὖ σμικρότερα ὅταν ἐν τοῖς μείζοσιν πολλοῖς \textit{περιλαμβανόμενα} ὀλίγα\\
διαθραυόμενα κατασβεννύηται, συνίστασθαι μὲν ἐθέλοντα εἰς τὴν τοῦ κρατοῦντος
ἰδέαν\\
πέπαυται  κατασβεννύμενα\\
\textit{γίγνεταί τε ἐκ πυρὸς ἀήρ,}\\
\textit{ἐξ ἀέρος ὕδωρ Timaeus} 57b3

\selectlanguage{english}
whenever a few of the smaller corpuscles, being caught within a great number of
larger corpuscles, are broken up and quenched, then, if they consent to be re
compounded into the shape of the victorious kind, they cease to be quenched,
and 

\begin{center}
\tablefirsthead{}
\tablehead{}
\tabletail{}
\tablelasttail{}
\begin{supertabular}{|p{\cellwidth}|p{\cellwidth}|p{\cellwidth}|p{\cellwidth}|}
\hline
 &
 &
 &
\cc{air [accusative]}\\ \hline
 &
 &
\cc{\itshape is generated\\ (gignetai)} &
\\\hline
\cc{\itshape from  (ek)} &
 &
 &
\\\hline
 &
\cc{fire [genitive],} &
 &
\\ \hline
\cc{and \textit{from (ex)}} &
 &
 &
\\\hline
 &
\cc{air\\
{[genitive]}} &
 &
\\\hline
 &
 &
 &
\cc{water [accusative];}\\ \hline
\end{supertabular}
\end{center}
From pre-existing fire, air is \textit{generated}, and 

from already generated air water is \textit{generated}. 

Again air cannot be said to \textit{consist} of fire, but it certainly is \textit{generated} from
fire.

\medskip

 [7] κατὰ φύσιν γὰρ σάρκες μὲν καὶ νεῦρα ἐξ αἵματος \textit{γίγνεται}

\medskip
 
For in the order of nature 
\begin{center}
\tablefirsthead{}
\tablehead{}
\tabletail{}
\tablelasttail{}
\begin{supertabular}{|p{\cellwidth}|p{\cellwidth}|p{\cellwidth}|p{\cellwidth}|}
\hline
 &
 &
 &
 \cc{flesh and sinews\\
   {[accusative]}}\\ \hline
 &
 &
\cc{\itshape are generated\\ (gignetai)} &
\\\hline
\cc{\itshape from (ex)} &
 &
 &
\\\hline
 &
\cc{blood [genitive],} &
 &
\\ \hline
\end{supertabular}
\end{center}
\textit{Timaeus} 82c8

\noindent From pre-existing blood, flesh and sinews are \textit{generated}

\medskip

[8] Οὐκ\\
ἐκ χρημάτων ἀρετὴ \textit{γίγνεται}, \\
ἀλλ' ἐξ ἀρετῆς χρήματα καὶ τὰ ἄλλα ἀγαθὰ τοῖς ἀνθρώποις ἅπαντα καὶ ἰδίᾳ καὶ
δημοσίᾳ. \textit{Apologia}  30b2--4

\medskip

\noindent It is not the case that  
\begin{center}
\tablefirsthead{}
\tablehead{}
\tabletail{}
\tablelasttail{}
\begin{supertabular}{|p{\cellwidth}|p{\cellwidth}|p{\cellwidth}|p{\cellwidth}|}
\hline
 &
 &
 &
\cc{virtue} \\ \hline
 &
 &
\cc{\itshape is generated\\ (gignetai)} &
\\ \hline
\cc{\itshape from (ek)} &
 &
 &
\\\hline
 &
\cc{money [genitive]} &
 &
\\\hline
\end{supertabular}
\end{center}

\noindent but 

\begin{center}
\tablefirsthead{}
\tablehead{}
\tabletail{}
\tablelasttail{}
\begin{supertabular}{|p{\cellwidth}|p{\cellwidth}|p{\cellwidth}|p{\cellwidth}|}
\hline
\cc{\itshape from} &
 &
 &
\\ \hline
 &
\cc{virtue} &
 &
\\\hline
 &
 &
\cc{\itshape is generated\\ (gignetai)} &
\\ \hline
 &
 &
 &
\cc{money and all other\\ good things to man,\\ both to the
individual\\ and to the state}
\\ \hline
\end{supertabular}
\end{center}

From pre-existing virtue, money is \textit{generated},\\
and not conversely, from money, virtue is \textit{generated}.

Also, once we interpret `gignetai' in terms of generation from re-existing parts, the adverbial rendering can only mean that one part could be generated from more than two parts, say from three, a meaningless statement in the context, since our starting point is exactly two parts.\\

There is one more question: `\textit{least part}' among which family of parts? Here we find an explanation for Plato's insistence in 142d6--8 (Section 3.3): that a part is not just a part, but a part of the whole. Thus the `least part' can only refer to `the least part of the whole', the smaller part of the two parts that add up to the whole.\\

Furthermore, if Plato wanted to say that

`the part turns out to \textit{consist}, or \textit{is composed}, of two parts'

\noindent (as, for example, Fowler, Allen, Palmer maintain),\\
then he would not say that

`to morion gignetai \textit{ek duoin morioin}'

\noindent but might rather say

`to morion gignetai \textit{duo moria}.

\smallskip

\noindent[Translations by \cite{Plato} with modifications by the author]

\subsubsection{Our rendering of the \emph{Parmenides} 142d9--e5}

Putting together 1.5.4. (i), (ii), (iii), and (iv), we obtain the following rendering of the \emph{Parmenides} 142d9--e5

\begin{quote} 
Since the One has the Being as part, (and the Being has the One as part), \\
again (`palin') [as is the case with the whole One \& Being], \\
it follows \textit{therefore (`ara')} that \\
the One part possesses both a One part and a Being part \\
and the Being part also possesses both a Being part and a One part, and \\
(in every case) the \textit{smallest (`to elachiston')} part \\
\textit{is generated from (gignetai ek)} two [pre-existing] parts.
\end{quote}

\subsubsection{Our interpretation of the \emph{Parmenides} 142d9--e5 leads to the first two anthyphairetic relations of the dyad One, Being}

Our starting point is the dyad One and Being. The One contains as part the Being (Section 3.3.4); the two parts One and Being generate, according to 3.4.4, another part, say $x$, smaller than the parts that generate it, hence smaller than Being; but One must possess two parts, a Being and a One, according to 3.4.4, hence the part $x$ must coincide with One$_1$; it follows that the new part One$_1$ of the One is generated by division/subtraction. This division of the One by the Being is seen to be\textit{ anthyphairetic},

$$\text{One}=\text{Being}+\text{One}_1, \quad \text{with } \text{Being}> \text{One}_1,$$
analogous to the division of the One \& Being by the One
$$\text{One \& Being}=\text{One}+\text{Being}, \quad \text{with } \text{One} >\text{Being}.$$
Thus, we now understand the reason why Plato, in his attempt to explain how the One has parts a Being and a One has gone to the upward anthyphairetic step, thus answering our Question in Section 3.3.2\&3.

But in fact, the generation of the One$_1$ part now explains the earlier claim, in 142d9--e3, that the Being contains a One part; in fact
$$\text{the Being contains the part One}_1.$$
In turn, this explains the earlier claim, in 142e3--4, that the Being part contains a One part and a Being part, and this sets forth the mechanism of generation of the third least part from the two already existing, in 142e4--5, that there is a Being part, say Being$_1$, such that
$$\text{Being}=\text{One}_1+\text{Being}_1, \quad \text{with } \text{One}_1>\text{Being}_1,$$
again an anthyphairetic relation. We have thus obtained the first two anthyphairetic relations of the anthyphairesis of dyad of ``hetera" One and Being.

\subsection{Proclus' comments on the \emph{Parmenides} 142d9--e5 in \emph{Platonic Theology} 3, 89, 21--26 provides full confirmation of our anthyphairetic interpretation of the second hypothesis (in Section 3.5)}

We now turn to Proclus' comments in the \emph{Platonic Theology} 3, 89, 21--26 for a strong confirmation of our interpretation of the crucial \emph{Parmenides} 142d9--e5 passage

\textbf{(a) Proclus' comments in \emph{Platonic Theology} 3, 89, 21--23 for the part One, translated 1816~\cite{ProclusTheology}, the passage translated by the author}

\begin{quote}
`Τὸ μὲν τοίνυν ἓν τοῦ ὄντος <μετέχον>\\
 διαιρεῖται πάλιν [[εἰς ὂν καὶ ἕν]]\\
  ὥστε τὸ ἓν καὶ τὸ ὂν ἀπογεννᾶν ἑνάδα δευτέραν,\\
   μοίρᾳ τοῦ ὄντος συνταττομένην' \\
On the one hand, the One, by participating in (metechon) the Being, is being divided (diaireitai) again [[in one and being]], so that the One and the Being generate (`apogennan') a secondary One (henada deuteran), equal (suntattomenen) to a part (moirai) of the Being. \emph{Proclus, Platonic Theology 3, 89, 21--23}
\end{quote}

Proclus' comments for the part One constitute a strong confirmation of our ongoing anthyphairetic interpretation. In fact:
\begin{enumerate}
    \item [(i)] Two parts One and Being `generate' (`apogennan') a new part One$_1$, as follows: the One contains the Being as part, and the two as a dyad generate a new partial One$_1$, in such a way that One$_1$ is equal to a part of Being. The Proclus description fully agrees with our interpretation of `gignetai ek' and `elachiston' in Section 3.5.5. Since the part One consists of the two parts Being and One$_1$, and since One$_1$ is generated from the dyad One$>$Being, there can be no other way of generation of One$_1$ than subtraction of the Being from the One.
    \item [(ii)] `the One participating (metechon) in Being' is divided into One and Being. means that the One has parts, in the following precise sense: the One consists of parts Being and One$_1$, with One$_1$ equal to a part of Being
\end{enumerate}

At long last we obtain the genuine meaning for the participation of the One into the Being, thus answering our Question 3.2:

\begin{quote}
The One participates in the Being (in 142b5--7, c5--7) 
\end{quote}

means precisely that

\begin{quote}
The One contains the Being as part, the (passive) One is divided by the (active) Being, and there is generated as anthyphairetic remainder One$_1$ smaller than the Being.
\end{quote}

\begin{enumerate}
    \item [(iii)] The presence of `palin' in Proclus' comments is an echo of the term `palin' in \emph{Parmenides} 142e3 and thus can only be understood in relation to the Whole, the One is in analogy, as, with the whole. In fact, the analogy is made clear with a juxtaposition of the two relations:
\end{enumerate}

Whole consists of two parts One and Being, with Being a part of One, and
One consists of two parts Being and One$_1$, with One$_1$ a part of Being.

We now turn into Proclus' comments for the Division of the Being

\textbf{(b) Proclus' comments in \emph{Platonic Theology} 3, 89, 24--26 for the part Being, translated 1816~\cite{ProclusTheology}, translation of the passage by the author}

\begin{quote}
Τὸ δὲ ὂν τοῦ ἑνὸς μετέχον\\
διακρίνεται πάλιν εἰς ὂν καὶ ἕν·\\
ἀπογεννᾷ γὰρ ὂν μερικώτερον ἑνάδος μερικωτέρας ἐξηρτημένον' \\
On the other hand, the Being, by participating (metechon) in the One,\\
is again being divided (diakrinetai) in Being and One;\\
because it generates (apogennai)\\
a Being more partial (merikoteran) dependent on (exertemenon)\\
a more partial (merikoteras) One.
\end{quote}

It is notable that Proclus deals first with the division of the One, and then, only after a byproduct of this division is the generation of One$_1$, equal to a part of Being, thus allowing the formation of the dyad One$_1 <$ Being, is in position to deal with the division of the Being. Thus the division of the Being does not take place independenty of the division of One, but on the contrary depends on it.

There is one point that needs be clarified. Is `the second one' (`henada deuteran') the same with `the more partial one' (`henados merikoteras')? Proclus probably considers this obvious, since every new part must be generated, this is done in a linear way, and, at the stage of the second relation, the only existing ones is the original One and `the second one', which is indeed a more partial one, since it is equal to a part of the Being, which is a part of the One. Note also that Proclus later in the \emph{Platonic Theology} refers to the division of the One in the more partial ones (`kermatizousa de to hen eis tas merikoteras enadas', 4,80,4--5). Thus `second one' One$_1$', generated by division, cannot but be one of the `more partial ones'.

Proclus comments interprets `gignetai ek', `to elachiston morion', `ara' and `metechei' of the \emph{Parmenides} passage in agreement with our interpretation in Section 3.4. Furthermore, both the \emph{Parmenides} passage and Proclus' comments on it confirm that the \emph{Parmenides} passage 142d9--e5 describes the first two steps of the anthyphairesis of One to Being.

\subsection{The general step of the anthyphairesis of One to Being 142e5--7}

After describing the first double step of the division of the dyad One, Being, the passage continues with the description of the general step, rendered as follows
\begin{quote}
καὶ κατὰ τὸν αὐτὸν λόγον οὕτως ἀεί,\\
ὅτιπερ ἂν μόριον γένηται,\\
τούτω τὼ μορίω ἀεὶ ἴσχει·\\
τό τε γὰρ ἓν τὸ ὂν ἀεὶ ἴσχει καὶ τὸ ὂν τὸ ἕν· \\
`and thus always (`aei') in the same manner (as before) (kata ton auton logon) whatever becomes (`genetai') a part,\\
always possesses (`ischei') a One and a Being\\
because (`gar', 142e7)\\
a (generated) One always (`aei') possesses (`ischei') as part a Being and\\
a (generated) Being always (`aei') possesses (`ischei') as part a One' 142e5--7
\end{quote}

`and thus always in the same manner' (kata ton auton logon houtos aei, 142e5) refers to the immediately preceding statement `the two parts generate a smallest part' (142e4--5) and turns the preceding statement into a general rule of generation (by division-subtraction) of a third least part from a dyad of already generated unequal parts.

Note that in 142e5--7 the generalized statement `every part possesses a One and a Being' holds because (`gar', 142e7) the generalized statement `every One possesses a Being, and `every Being posseses a One' holds, precisely in the same way that the initial statement `the One possesses a Being and the Being possesses a One' holds true, hence (`ara', 142e3) both the One and the Being possess a One and a Being'. Thus the term `gar' establishes precisely the same causal relation and has exactly the same force for the general case (142e5--7) that the term `ara' (142e3) for the initial case (142d9--e5).

We are then led in modern notation to the precise statement that
\begin{itemize}
    \item for every natural number $k$ the One$_k$ possesses (`ischei') the Being$_k$ and the One$_{k+1}$, with the One$_{k+1}$ equal to a part of the Being$_k$, and
    \item the Being$_k$ possesses (`ischei') the One$_{k+1}$ and the Being$_{k+1}$, with the Being$_{k+1}$ equal to a part of the One$_{k+1}$.
\end{itemize}
It is clear that the reasoning is by mathematical recursion.

\subsection{The dyad One and Being satisfies the philosophical analogue of infinite anthyphairesis 142e7--143a3}

The result of the inductive argument is that every part is ad infinitum divided into two parts
\begin{quote}
`ὥστε ἀνάγκη δύ' ἀεὶ γιγνόμενον μηδέποτε ἓν εἶναι' (142e7--143a2),
\end{quote}

hence the dyad One and Being is in fact infinite in multitude
\begin{quote}
    `Οὐκοῦν ἄπειρον ἂν τὸ πλῆθος οὕτω τὸ ἓν ὂν εἴη' (143a2--3),
\end{quote}

the philosophic anthyphairesis is infinite. Thus, the anthyphairetic division of the `indefinite dyad $\langle\text{One, Being}\rangle$ yields
\begin{eqnarray*}
  \text{One} &=& \text{Being} +\text{One}_1, \quad \text{Being} > \text{One}_1,\\
  \text{Being} &=& \text{One}_1+\text{Being}_1, \quad \text{One}_1 > \text{Being}_1, \\
  &\vdots& \\
  \text{One}_n &=& \text{Being}_n +\text{One}_{n+1}, \quad \text{Being}_n > \text{One}_{n+1}, \\
\text{Being}_n &=& \text{One}_{n+1}+\text{Being}_{n+1}, \quad \text{One}_{n+1} > \text{Being}_{n+1},\\
             &\vdots& 
\end{eqnarray*}
And there is then an
infinite multitude (`apeiron to plethos', 143e2) of remainders-parts
of the anthyphairetic division:
$$\text{One} > \text{Being} > \text{One}_1 > \text{Being}_1 > \ldots > \text{One}_n > \text{Being}_n > \ldots .$$

\subsection{Proclus' comments on the \emph{Parmenides} 142e5--143a3-- in \emph{Platonic Theology} 4, 94, 11--16 completes his comments in \emph{Platonic Theology} 3, 89, 21--26, given in Section 3,5, and provides full confirmation of our anthyphairetic interpretation of the second hypothesis (in Sections 3.4,6,7)}

In a remarkable passage in the \emph{Platonic Theology}, Proclus describes the general step of the ge where Proclus confirms that we are in precisely this seemingly paradoxical situation. Proclus first describes the infinite multitude of parts of the One and of the Being, generated by anthyphairetic division, in the \emph{Parmenides} 142b1--143a3:
\begin{quote}
The One and Being becomes Many by the principle of the Otherness dividing each of them (καὶ τὸ ἓν πολλὰ καὶ τὸ ὂν γίγνεται, τῆς ἑτερότητος διακρινούσης ἑκάτερον). Thus every part of the Being participates in the One, and (Καὶ πᾶν μὲν τοῦ ὄντος μόριον μετέχει τοῦ ἑνός), While every part of the one participates in a part of the Being (πᾶσα δὲ ἑνὰς ἐποχεῖται μοίρᾳ τινὶ τοῦ ὄντος·) and so each of them is increased in multitude and is divided intelligibly and is fragmented \textit{ad infinitum} (πληθύεται δὲ ἑκάτερον καὶ διακρίνεται νοερῶς καὶ κατακερματίζεται καὶ ἐπ' ἄπειρον πρόεισιν.) \emph{Proclus, Platonic Theology 4, 94, 11--16;}
\end{quote}

\subsection{The mathematical origin of the upward step in One\&Being as attested by Proposition XIII.5 of the \emph{Elements} and by Proclus' Elegant theorem}

\subsubsection{Once we realize}
\begin{itemize}
    \item [(a)] that the aim of the \emph{Parmenides} 142b1--143a3 passage is to show that the dyad One and Being, formed within the context of the second hypothesis, satisfies the philosophic analogue of infinite anthyphairesis,
    \item [(b)] that the aim of the upward passage to the whole in the 142d1--9 passage, as explained in Section 3.3, is to elucidate and facilitate this philosophic analogue of anthyphairesis, and
    \item [(c)] that there is a strikingly similar Pythagorean geometrical method of passing to the whole for the purpose of obtaining full knowledge of the anthyphairesis of the diameter to the side of a square and the mean and extreme ratio,
\end{itemize}
we may come to the natural conclusion that Plato's philosophic use of passing to the whole is inspired from the corresponding mathematical use. We suggest that Plato took this upward step to the whole in 142d1--9 in close imitation of similar upward steps to the whole in geometry, in fact in Proposition XIII.5 of the \emph{Elements} and in the so called by Proclus elegant theorem.

\subsubsection{`palin' (142e3) means `as with the upward step to the whole'}

Thus 142d8--e3 indeed explains how because of the participation of the One in the Being, the Being has parts, is divided into parts, as claimed in 142c8--d1. 
But we still have to provide a satisfactory explanation for the term `palin' appearing in 142e3. 
The term, rendered as `again', indicates repetition, and there is only one conceivable previous situation of which the present one is a repetition. 
Indeed `palin' cannot but refer to the fact that the statement that the One consists of the parts Being and One$_1$, and in fact with One$_1$ a part of Being, 
is a repetition of the statement about the Whole that the One \& Being consists of the parts One and Being. 

Thus the fact that Plato went one step upward to the whole, was for the purpose of elucidation the way in which the part One is being divided into parts because of its participation in Being. 
This way has turned out to be a philosophic version of anthyphairesis, and also similar to the way in which the whole One \& Being is being divided into parts. 

Thus the upward step is towards elucidation and guidance for the meaning of participation of the One in Being, and for the manner in which the One is divided by the Being.

\subsubsection{Two similar upward anthyphairetic steps to the whole in geometry}

The Pythagorean incommensurabilities are essentially two, that of the diameter to the side of a square and that of the mean and extreme ratio. From ancient sources it is possible to conjecture that the proof of these incommensurabilities followed a similar course involving what I will call the upward anthyphairetic step to the whole.

\begin{itemize}
    \item [{[1]}] \textbf{Proposition XIII.5 of the \emph{Elements}} 
    
    If the dyad $a>b$ is in mean and extreme ratio, then the dyad $a+b>a$ is also in mean and extreme ratio.
    \end{itemize}

\textbf{Proof}. By Proposition II.11 of the \emph{Elements}, assumption $a^2=ab+b^2$. Add to both sides $a^2+ba$, and use Proposition II.4 of the \emph{Elements}, to conclude $(a+b)^2=(a+b)a+a^2$.

\textbf{Corollary}. If $a>b$ is in mean and extreme ratio, then $\text{Anth} (a,b)=[\text{period}(1)]$ and $a,b$ are incommensurable. [Here and below by $\text{Anth}(a,b)$ we denote the sequence $[k_1,k_2,k_3,\ldots]$ of successive quotients of the anthyphairesis of $a$ to $b$.]

\textbf{Proof of Corollary}. By Proposition XIII.5, $\text{Anth}(a, b)=\text{Anth}(a+b, a)$. 

Then $\text{Anth} (a, b) = [1, \text{Anth} (b, a-b)]$ 

Hence $[1, \text{Anth} (b, a-b)] =[1, 1, \text{Anth} (b, a-b)]$. 

Hence $\text{Anth} (b, a-b) = [1, \text{Anth} (b, a-b)]$. 

Hence, inductively, $\text{Anth} (b, a-b) = [1,1,1,\ldots]$. 

Hence $\text{Anth} (a, b)=[1,1,1,1,\ldots]$. 

Incommensurability of $a$ to $b$ follows from Proposition X.2 of the \emph{Elements}.

\begin{itemize}
    \item [{[2]}] \textbf{Proposition on the diameter to the side of a square} 
    
    (the Elegant Theorem in Proclus, \emph{In Platonis rem publicam} 2,27,11--28, 10: the step upward to the whole) 
    
    If $a$ is the diameter to the side $b$, then $a+2b$ is the diameter, to the side $a+b$.
\end{itemize}

\textbf{Proof}. By assumption $a^2=2b^2$.

By Proposition II.10 of the \emph{Elements}, $(a+2b)^2+a^2=2(a+b)^2+2b^2$.

It follows immediately that $(a+2b)^2=2(a+b)^2$.

\textbf{Corollary}. $\text{Anth}(a,b) =[1,2,2,2,\ldots]$, and $a, b$ are incommensurable.

\textbf{Proof} (D. Fowler, 1994 \cite{Fowler1994}). (ii) It is easy to see that $\text{Anth}(a,b)=[1, \text{Anth}(b,a-b)]$. By (i) (Elegant theorem: the step upward to the whole), $\text{Anth}(a,b)=\text{Anth}(a+2b, a+b)$. By the Elegant theorem, $\text{Anth}(a+2b, a+b)=[1,2, \text{Anth}(b,a-b)]$. Hence $[1, \text{Anth}(b,a-b)]=[1, 2, \text{Anth}(b,a-b)]$. Hence $\text{Anth}(b,a-b)=[2, \text{Anth}(b,a-b)]$. Hence, recursively, $\text{Anth}(b,a-b)=[2,2,2,\ldots]$. Hence $\text{Anth}(a,b)=[1, 2,2,2,\ldots]$. Incommensurability of $a$ to $b$ follows from Proposition X.2 of the \emph{Elements}.

\subsection{The close connection of the statement the One participates in the Being with the statement the One and the Being are other and opposite to each other 142b7--c5, 143b1--8}

\subsubsection{We first examine passage 142b7--c5}
\begin{quote}
— Οὐκοῦν καὶ ἡ οὐσία τοῦ ἑνὸς εἴη ἂν οὐ ταὐτὸν οὖσα τῷ ἑνί·\\
οὐ γὰρ ἂν ἐκείνη ἦν ἐκείνου οὐσία,\\
οὐδ' ἂν ἐκεῖνο, τὸ ἕν, ἐκείνης μετεῖχεν,\\
ἀλλ' ὅμοιον ἂν ἦν λέγειν ἕν τε εἶναι καὶ ἓν ἕν.\\
νῦν δὲ\\
οὐχ αὕτη ἐστὶν ἡ ὑπόθεσις, εἰ ἓν ἕν, τί χρὴ συμβαίνειν,\\
ἀλλ' εἰ ἓν ἔστιν· οὐχ οὕτω;\\
— Πάνυ μὲν οὖν.\\
— Οὐκοῦν ὡς ἄλλο τι σημαῖνον τὸ ἔστι τοῦ ἕν;\\
— Ἀνάγκη. \\
``Then the Being (ousia) of the One will be,\\
but will be \textbf{not same} (ou tauton) with the One;\\
for if it were identical with the One,\\
it would not be the Being of the One (ekeine ekeinou ousia),\\
nor would the One participate in it,\\
but the statement that `One is' would be equivalent to\\
the statement that `One is One' (hen hen);\\
but our hypothesis is not if `one is one' (hen hen), what will follow,\\
but if One is. Do you agree?"\\
``Certainly."\\
``[Do you agree] that the Being is something \textbf{other} (allo) than the One?"\\
``Most assuredly." 142b7--c5
\end{quote}

The two parts, the One and the Being, are \textbf{not same} (`ou tauton',142b7), 
but \textbf{different, other} (allo, 142c4) to each other, and thus form a dyad of others. 
Thus 
\begin{quote} 
    \textit{the One and the Being form a dyad of other and opposites to each other}.
\end{quote} 
With this statement the second hypothesis is differentiated right from the start from the first hypothesis. 
Indeed, the One of the first hypothesis is partless, but, according to the second interpretation of the hypothesis ``the One is", the One of the second hypothesis has part the Being. 
The passage also sets the two parts One and Being as, not simply different, but \textbf{other} to each other. 
The non-symmetrical expression ``the Being of the One" 142b7, 142b8 (there is no corresponding expression ``the One of the Being") suggests that 
\begin{quote}
    \textit{the Being is a part of the One}.
\end{quote}
This suggestion will be reinforced in the sequel.

\subsubsection{The second passage correlating participation and otherness is 143b1--8}

In the later passage 143b1--8 the nature of the two parts One and Being as ``other" to each other is emphasized
\begin{quote}
άλλο τι ἕτερον μὲν ἀνάγκη τὴν οὐσίαν αὐτοῦ εἶναι, ἕτερον δὲ αὐτό, εἴπερ μὴ οὐσία τὸ ἕν, ἀλλ' ὡς ἓν οὐσίας μετέσχεν. — Ἀνάγκη. — Οὐκοῦν εἰ ἕτερον μὲν ἡ οὐσία, ἕτερον δὲ τὸ ἕν, οὔτε τῷ ἓν τὸ ἓν τῆς οὐσίας ἕτερον οὔτε τῷ οὐσία εἶναι ἡ οὐσία τοῦ ἑνὸς ἄλλο, ἀλλὰ τῷ ἑτέρῳ τε καὶ ἄλλῳ ἕτερα ἀλλήλων. — Πάνυ μὲν οὖν. — Ὥστε οὐ ταὐτόν ἐστιν οὔτε τῷ ἑνὶ οὔτε τῇ οὐσίᾳ τὸ ἕτερον. — Πῶς γάρ; 143b1--8 \\
\ldots the Being of the One is \textbf{other} (heteron) and One is \textbf{other} (heteron), since the One is not the Being, but as One, participates (meteschen) in Being?" ``Yes, that is necessary. ``Then if the Being is \textbf{other} (heteron) and the One is \textbf{other} (heteron), the One is \textbf{other} (heteron) than the Being not because it is the One, the Being is \textbf{other} (allo) than the One not because it is the Being, but they are \textbf{other} (hetera) to each other, differ from each other by virtue of \textbf{otherness} (toi heteroi) and difference (alloi)." ``Certainly." ``Therefore, the \textbf{other} (heteron) is neither the same (ou tauton) as the One nor as the Being." ``Certainly not." 143b1--8
\end{quote}

The passage makes clear that otherness is a relation between the two parts, not something in the nature of each part separately; and that this otherness is a consequence of the participation of the One into the Being.

Note. In the present section we are intereothrness does not play some central role.

\subsection{Our interpretation of the principle of otherness (a) as having its roots in Fragment B3 of Zeno, and (b) as expressing the process of the philosophical analogue of infinite anthyphairesis}

In fact, otherness of the One and the Other is equivalent to participation of the One in the Other. There is close similarity between the infinite part of Zeno's Fragment B3 and the generation of an infinite multitude of parts in the One of the Second Hypothesis in Plato's \emph{Parmenides} 142b1--143b8 Zeno's Fragment B3 antedates the description of the dyad One, Being as a dyad of `heteron to heteron' in the \emph{Parmenides} 142b7-c7, 143b1-8.

\begin{quote}
Καὶ τί δεῖ πολλὰ λέγειν, ὅτε καὶ ἐν αὐτῷ φέρεται τῷ τοῦ Ζήνωνος συγγράμματι; πάλιν γὰρ δεικνύς, ὅτι εἰ πολλά ἐστι, τὰ αὐτὰ πεπερασμένα ἐστὶ καὶ ἄπειρα, γράφει ταῦτα κατὰ λέξιν ὁ Ζήνων· \\
``εἰ πολλά ἐστιν, ἀνάγκη τοσαῦτα εἶναι ὅσα ἐστὶ καὶ οὔτε πλείονα αὐτῶν οὔτε ἐλάττονα. εἰ δὲ τοσαῦτά ἐστιν ὅσα ἐστί, πεπερασμένα ἂν εἴη." \\
``εἰ πολλά ἐστιν, ἄπειρα τὰ ὄντα ἐστίν. ἀεὶ γὰρ ἕτερα μεταξὺ τῶν ὄντων ἐστί, καὶ πάλιν ἐκείνων ἕτερα μεταξύ. καὶ οὕτως ἄπειρα τὰ ὄντα ἐστίν." \\
And why should I say any more, for it also exists in the book of Zeno? For again, showing that if the Many are true Beings, then the same (ta auta) will be finite (peperasmena) and infinite (apeira), Zeno writes thus verbatim: [the finite property of Zeno's Fragment B3] `If the Many are true Beings, then necessarily the Many are as many (tosauta) as (hosaper) they are, and neither more (oute pleiona) of them nor fewer (oute elattona). But if they are as many (tosauta) as they (hosa) are, they would be finite. [the infinite property of Zeno's Fragment B3] If the Many are true Beings, the Many are infinite beings. For there are always \textbf{other} (hetera) \textbf{in between} (metaxu) beings, and again \textbf{other} (hetera) \textbf{in between} (metaxu) them. And in this way the beings will be infinite.' \emph{Simplicius, Commentary to Aristotle's Physics 140,27--33}
\end{quote}

We note that Zeno in his Fragment B3 shows that his true Beings have an infinite multitude of parts, by having always new hetera appearing between (metaxu) previous hetera ad infinitum. The term in between (metaxu) has been rendered by earlier scholars as between in the sense of order, similar to the image of a rational appearing between two reals. But as I have argued in detail in Negrepontis, 2023~\cite{Negrepontis2023}, the real meaning of metaxu two hetera in Zeno's Fragment B3 is that of an offspring (ekgonon) from two previous parts/hetera to each other, as used often by Plato himself and as noted by Simplicius in his work, where Zeno's Fragment B3 appears. A term fully equivalent to this meaning of ``metaxu+genitive" is ``gignetai ek+genitive", thus establishing a very close similarity of the description of the infinite in Zeno's Fragment B3 in terms of ``heteron" and ``metaxu" with Plato's description of the infinite in the Second Hypothesis in the \emph{Parmenides} in terms of ``heteron" or ``allo" (142c4, 143b1,b2,b3,b4,b5,b6,b7) and ``gignetai" (142d5,e4,e6, 143a1). It is clear that the statement in the \emph{Parmenides} 142c4 and 143b that One and Being is a dyad of parts hetera to each other is an expression equivalent to a later Platonic term that of indefinite dyad, or a Philebean Infinite, essentially meaning that it is a dyad in philosophical infinite anthyphairesis. We note that the close similarity between Zeno's Fragment B3 and the Second Hypothesis in the \emph{Parmenides} does not stop here; Plato follows most closely the finite part of Zeno's Fragment B3, as well. The reader is referred to Negrepontis, 2023~\cite{Negrepontis2023} and to a future work by the author on the close connection between Zeno's and Plato's philosophical arguments.

We conclude, answering our Question 1.3, that the Principle of Otherness that Plato enunciates in the \emph{Parmenides} Second Hypothesis \textbf{(a)} has its origin in the infinite part of Zeno's Fragment B3, and \textbf{(b)} expresses the philosophical analogue of infinite anthyphairesis of the paradigmatical intelligible dyad One and Being.

\subsection{A second proof of the anthyphairetic infinity of the dyad One, Being, dialectical, instead of the present eristic one? 143a4--b1}

It has been just proved that the dyad One, Being (`to hen on') is a dyad satisfying the philosophical analogue of infinite anthyphairesis. As a consequence, an infinite multitude of parts of the One is generated anthyphairetically. At this point Plato poses a question:

\begin{quote}
Ἴθι δὴ καὶ τῇδε ἔτι. — Πῇ; — Οὐσίας φαμὲν μετέχειν τὸ ἕν, διὸ ἔστιν; — Ναί. — Καὶ διὰ ταῦτα δὴ τὸ ἓν ὂν πολλὰ ἐφάνη. — Οὕτω. — Τί δέ; αὐτὸ τὸ ἕν, ὃ δή φαμεν οὐσίας μετέχειν, ἐὰν αὐτὸ τῇ διανοίᾳ μόνον καθ' αὑτὸ λάβωμεν ἄνευ τούτου οὗ φαμεν μετέχειν, ἆρά γε ἓν μόνον φανήσεται ἢ καὶ πολλὰ τὸ αὐτὸ τοῦτο; — Ἕν, οἶμαι ἔγωγε. — Ἴδωμεν δή· 143a4--b1 \\
``Let us make another fresh start." ``In what direction?" ``We say that the one participates (`metechein') of being, because it is?" ``Yes." ``And for that reason, the One\&Being (`to hen on') was found to be many (`polla')." ``Yes." ``Well then, will the One Itself (`auto to hen'), which we say participates (metechein') of Being, if we consider it intelligibly (tei dianoiai labomen') alone by itself (monon kath' hauto', `auto touto'), without that of which we say it participates (metechein'), be found to be only (`monon') one, or many (`polla')?" ``One, I should say." ``Just let us see;'' 143a4--b1
\end{quote}

Now the question that Plato poses, in 143a4--b1, is: consider the One, not as a part of the One Being, an indefinite dyad, but as a Platonic intelligible Being (`auto to nen', `to auto touto'). Will then this Platonic Being itself have an infinite multitude of parts? The answer cannot be given yet, since the nature of a Platonic Being has not been explained yet, but will be given later, in Section 5, in 144c2--d4, 138a3--7, 148d--149d, 155d--e, in the positive, just after a Platonic Being has been at last understood. It has been proved that the dyad $\langle\text{One, Being}\rangle$ possesses the philosophic analogue of infinite anthyphairesis. Proposition X.2 of the \emph{Elements} states that if the anthyphairesis of a pair of magnitudes $\langle a, b\rangle$ is infinite, then $a, b$ are incommensurable (namely there is no magnitude $c$ and natural numbers $m, n$ such that $a=mc, b=nc$). Thus, we may think of the Platonic statement that the philosophical anthyphairesis of the dyad One to Being is infinite to be a philosophic analogue of the incommensurability of the One and the Being (jn analogy to Proposition X.2 of the \emph{Elements}). This analogy with Mathematics is useful for understanding the nature of the question set by Plato. An analogous problem in Mathematics would be: it is a fact that the diameter to the side of a square is incommensurable. There is an arithmetical proof of this incommensurability involving odd and even number, as reported by Aristotle in \emph{Analytics Prior} 41b and, much later, as Proposition X.117, an addition to Book X of the \emph{Elements}, perhaps by Theon of Alexandria. The arithmetical proofs of incommensurability originated by Archytas and use the techniques of Book VII and VIII of the \emph{Elements}. In some sense the proof, making an appeal to the converse of X.2, shows that the anthyphairesis of the diameter to the side is infinite. But note that this arithmetical proof of incommensurability does not provide knowledge and does not inform us in any other way about the anthyphairesis, which remains wholly unknown, except for the fact of being infinite. Plato, by asking for a dialectical proof of infinity of the anthyphairesis, and by showing, as we shall see in the next Section 4, that only infinite anthyphairesis is not sufficient to construct ``philosophical"/dialectical numbers, namely numbers with equalized units, is possibly making an indirect criticism of Archytas arithmetical proofs. He will be much blunter in the \emph{Philebus} 16e--17a passage, where, still not naming Archytas, he refers to such arithmetical methods as eristic and not dialectical.

Now if we first obtain a complete knowledge of the anthyphairesis of the diameter to the side, and know that it is equal to $[1,2,2,2,\ldots]$, then we have a second and superior proof of the infinity and incommensurability. We will see in Sections 5,6, below, that this is precisely the nature of Plato's question.

\section{Eristic numbers, generated with units the terms of the infinite sequence of the anthyphairetic remainders of the dyad One \& Being of the second hypothesis of the Parmenides, unequal to each other, introduced and eventually rejected (143c1--144c2)}

In 143d5--144c2 there is an abortive attempt to generate in the One of the Second Hypothesis, eristic non-dialectic numbers, in the sense that they consist of unequal units, by employing only the infinite sequence of the anthyphairetic remainders of the anthyphairesis of the One, Being, resulting in \foreignlanguage{greek}{“pleista”}/most numerous parts of the Being, namely resulting in an infinite multitude of units for numbers 143d5--144c2; the attempt is rejected at 144d4--e3, exactly because the units of these numbers are mutually unequal. In 143c1--d5 we note the definition of the dialectical number Two, which is however premature, since this definition depends on anthyphairetic periodicity and logos, which have not been introduced yet, making the passage impossible to understand (as explained later in Section 6.7).

Since we have verified that the One of the Second Hypothesis satisfies a philosophical analogue of infinite anthyphairesis, we next expect, with greater confidence that, in accordance with our findings in the \emph{Theaetetus, Sophist, Meno} (a) that an intelligible Being has the structure of a dyad in the philosophical analogue of periodic anthyphairesis, and (b) that the knowledge of an intelligible Being can be given in terms of Name plus Logos,

\begin{quotation}
the dyad One Being will satisfy a philosophical analogue of periodic anthyphairesis.
\end{quotation}

Thus, following Section 3, the natural step would be to prove that

\begin{quotation}
the anthyphairesis of the dyad $\langle\text{One, Being}\rangle$, already proved to be infinite, is in fact periodic.
\end{quotation}

   We will presently show that this is indeed what Plato has in mind; But we find that right after the proof that the dyad One and Being satisfies the philosophical analogue of infinite anthyphairesis, that Plato \textit{does not} directly proceed to establish periodicity of the anthyphairesis, but instead shows interest to exploit the infinite sequence of the anthyphairetic remainders of the philosophical anthyphairesis of the dyad One to Being, as units for the generation of numbers.

Plato’s description of the definition of the number Two in 143c1--d5 is a special remarkable passage that will become understandable later, after anthyphairetic periodicity is established (4.1).

Plato proceeds from 143d5 to 144c2 to generate numbers in the One of the second hypothesis by employing as units the infinite sequence of anthyphaireic remainders. 
Heuristically, if we have already read the \emph{Philebus} 56d--e\\
generate numbers in the One of the second hypothesis is eventually rejected (at 144d5--7), because the units of these numbers are unequal to each other, thus resulting in \textit{eristic, not dialectical, numbers} (143c1--144c2), a requirement most lucidly described in passage. It seems that the purpose for this section in the second hypothesis of Plato’s \emph{Parmenides} is to emphasize that the \textit{dialectical numbers}, introduced in Section 6, need in an essential way the anthyphairetic periodicity of the dyad One\&Being.

This is an abortive attempt to introduce numbers in the One of the Second Hypothesis, solely by means of the infinite anthyphairesis, an attempt that is doomed to failure, because the units for these \textit{“eristic”} numbers are the anthyphairetic remainders, by their very structure unequal to each other (\emph{Parmenides} 143d5-144c2). The whole passage on the abortive attempt to introduce eristic numbers should be regarded as an argument to indicate that dialectic numbers cannot be introduced merely by the infinity of the anthyphairesis, but periodicity is needed.

A number in Greek Mathematics is always a (finite) multitude of units (cf. definition VII.2 of number in Euclid’s \emph{Elements}); we will be able to generate all number if we have at our disposal an infinite multitude of terms that can serve as units. What are to be the units for generating number in the One? A natural candidate of an infinite multitude of units is provided by the infinite sequence of the anthyphairetic remainders of the anthyphairesis of the dyad One and Being. So why not take the infinite sequence of the anthyphairetic remainders as the sequence of units for the generation of all numbers:

$$\text{One} > \text{Being} > \text{One}_1 > \text{Being}_1 > \ldots > \text{One}_n > \text{Being}_n > \ldots .$$

E.g. the two parts One and Being can be the units for the generation of the number Two, the parts One, Being, One$_1$ for the generation of the number Three, and so on.

This approach initially appears to be fine, but a deeper understanding of Plato would lead someone to realize that, on the contrary, the approach is doomed to failure. The place where Plato expresses in clearest terms the problem with this approach is in the \emph{Philebus} 56d4--e6 (examined in Section 4.5.b, below), also in the \emph{Republic} 526a--b): the numbers formed by the many, ignorant persons consist of unequal unis, such as two tables, two people, and so on, the numbers that might well be called \textit{eristic numbers}, but 

\begin{quotation}
the philosophical, dialectic numbers are necessarily composed of units that are equal to each other.
\end{quotation}

But the remainders of the anthyphairesis of One to Being is a strictly decreasing sequence, any two units are \textit{unequal} to each other, and thus the numbers introduced in 143d5-144c2 should rather be the eristic ones. We are not surprised that eventually the whole approach is rejected: the method ends with showing that the parts of the Being are\\
\noindent \textit{not pleista, infinite in number}, as claimed in 144c2,\\
\noindent \textit{but infinite only in multitude}, being able to generate only eristic, not dialectical numbers, finally rejected at 144d5--e3.

The first abortive attempt to introduce numbers consists of three stages.

\subsection{The premature definition of the dialectical number Two 143c1--d5 (explained in Section 6.7, after the introduction of the dialectical number in 148d--149d and 144d4--e1)}

As will become evident in Sections 5 and 6 below, the number Two is here defined in a dialectical manner, based on the explicit requirement that the One and the Being possess Logos. This is not explained so far, and thus this whole passage cannot be understood yet. We shall see below, in Section 5, that there is a passage, in 138a3--7, most crucial for the whole dialogue \emph{Parmenides}, that cannot be understood at the early point it is stated. We shall see later that both 138a3--7 and 143c1--144a5 rely on the plus one rule, explained in 148d--149d, and the presence of Logos, mentioned explicitly in 155d8--e1 and explained in 138a3--7.

These will definitely serve as the units in his attempt to generate, introduce numbers in the One. Thus the number Two can have units the One and the Being (143c1--d5):

\begin{enumerate}
    \item [1)] Can One and Being form the number two?
    \begin{quotation}
    — Τί οὖν; ἐὰν προελώμεθα αὐτῶν\\ 
    εἴτε βούλει τὴν οὐσίαν καὶ τὸ ἕτερον\\
    εἴτε τὴν \textit{οὐσίαν} καὶ τὸ \textit{ἓν}\\
    εἴτε τὸ ἓν καὶ τὸ ἕτερον,\\
    ἆρ' οὐκ ἐν ἑκάστῃ τῇ προαιρέσει προαιρούμεθά τινε\\
    ὣ ὀρθῶς ἔχει καλεῖσθαι \textit{ἀμφοτέρω};\\
    — Πῶς;
    \end{quotation}

    \item [2)] the One of the second satisfies the Logos Criterion for anthyphairetic periodicity (to be established later, by appealing to 138a3--7 and 148d--149d):
    \begin{quotation}
    — Ὧδε· ἔστιν οὐσίαν \textit{εἰπεῖν} (143c5);\\ 
    — Ἔστιν. \\ 
    — Καὶ αὖθις \textit{εἰπεῖν} ἕν (143c5); \\
    — Καὶ τοῦτο.\\ 
    — Ἆρ' οὖν οὐχ ἑκάτερον αὐτοῖν \textit{εἴρηται} (143c6); \\ 
    — Ναί.
    \end{quotation}

    \item [3)] The dyad One and Being forms a contact (hapsis), namely a ratio of two successive terms in the infinite sequence of the remainders of the anthyphairesis of One to Being (by appealing to 148d--149d):
    \begin{quotation}
    — Τί δ' ὅταν \textit{εἴπω} (143c7) οὐσία τε καὶ ἕν, ἆρα οὐκ \textit{ἀμφοτέρω};\\
     — Πάνυ γε.\\
     — Οὐκοῦν καὶ ἐὰν οὐσία τε καὶ ἕτερον ἢ ἕτερόν τε καὶ ἕν,\\
     καὶ οὕτω πανταχῶς ἐφ' ἑκάστου \textit{ἄμφω λέγω} (143c9);\\
      — Ναί.\\
      — Ὣ δ' ἂν \textit{ἄμφω} ὀρθῶς προσαγορεύησθον,\\
       ἆρα οἷόν τε \textit{ἄμφω} μὲν αὐτὼ εἶναι, δύο δὲ μή;\\
       — Οὐχ οἷόν τε.
    \end{quotation}

    \item [4)] Now we apply \textit{the plus one rule:} Dialectical number of the terms One and being$=$contact$+1=1+1=2$; thus that the One and the Being are two units. form the dialectical number Two. In fact the units One and being are equalized.
    \begin{quotation}
    — Ὣ δ' ἂν δύο ἦτον, ἔστι τις μηχανὴ μὴ οὐχ ἑκάτερον αὐτοῖν ἓν εἶναι;\\
     — Οὐδεμία. \\
     — Τούτων ἄρα ἐπείπερ σύνδυο ἕκαστα συμβαίνει εἶναι,\\
      καὶ ἓν ἂν εἴη ἕκαστον.\\
       — Φαίνεται.
    \end{quotation}
\end{enumerate}

We are not yet ready to discuss, or even translate, this deep and subtle passage. The complex justification, given in 143c1--d5, that One and Being form the number Two will be discussed below in Section 6.7, after the clear definition of the dialectical numbers.

\subsection{The generation of eristic numbers from the parts of the dyad One, Being 143d5--144a1}

\subsubsection{The number Three is generated with units the One and the Being, plus any other part, say One$_1$ (143d5--7).}

\begin{quotation}
— Εἰ δὲ ἓν ἕκαστον αὐτῶν ἐστι, συντεθέντος ἑνὸς \textit{ὁποιουοῦν ᾑτινιοῦν} συζυγίᾳ οὐ \textit{τρία} γίγνεται τὰ πάντα;\\
 — Ναί. \\
“But if each of them is one, by the addition of any sort of one to any pair whatsoever the total becomes three?”\\
 “Yes.”
\end{quotation}

\subsubsection{If numbers $a$ and $b$ are generated in the One, then the product $a \cdot b$, $a$ times $b$, is also generated in the One (143d7--144a1).}

\begin{quotation}
— Τρία δὲ οὐ περιττὰ καὶ δύο ἄρτια;\\
 — Πῶς δ' οὔ;\\
  — Τί δέ; δυοῖν ὄντοιν οὐκ ἀνάγκη εἶναι καὶ δίς,\\
   καὶ τριῶν ὄντων τρίς,\\
   εἴπερ ὑπάρχει τῷ τε δύο τὸ δὶς ἓν καὶ τῷ τρία τὸ τρὶς ἕν;\\
   — Ἀνάγκη.\\
   — Δυοῖν δὲ ὄντοιν καὶ δὶς οὐκ ἀνάγκη δύο δὶς εἶναι;\\
   καὶ τριῶν καὶ τρὶς οὐκ ἀνάγκη αὖ τρία τρὶς εἶναι;\\
   — Πῶς δ' οὔ;\\
   — Τί δέ; τριῶν ὄντων καὶ δὶς ὄντων καὶ δυοῖν ὄντοιν καὶ τρὶς ὄντοιν οὐκ ἀνάγκη τε τρία δὶς εἶναι καὶ δύο τρίς; \\
   — Πολλή γε.

— Ἄρτιά τε ἄρα ἀρτιάκις ἂν εἴη καὶ περιττὰ περιττάκις καὶ ἄρτια περιττάκις καὶ περιττὰ ἀρτιάκις.\\
 — Ἔστιν οὕτω.

“And three is an odd number, and two is even?”\\
“Of course.” \\
“Well, when there are two units, must there not also be twice, and when there are three, thrice, that is, if two is twice one and three is thrice one?”\\
“There must.”\\
“But if there are two and twice, must there not also be twice two?\\
And again, if there are three and thrice, must there not be thrice three?”\\
“Of course.”\\
“Well then, if there are three and twice and two and thrice, must there not also be twice three and thrice two?”\\
“Inevitably.”

“Then there would be even times even, odd times odd, odd times even, and even times odd.”\\
“Yes.”
\end{quotation}

\textit{Note}. This statement that if $a$ is a number generated in the One Being, then the number $2a$ is also generated in the One Being, which holds for the eristic numbers by 143d7--144a1, will be explicitly rejected for the dialectical in a revealing passage by Plotinus that will be discussed below in Section 6.8.

\subsection{The definition of eristic number in the second hypothesis of the \emph{Parmenides}}
Statement A: \textit{Every number is generated in the One Being} 144a2--5

\noindent Statement B: \textit{Generation of units for Eristic numbers: every part of the Being part can serve as unit for the generation of a number in the One of the second hypothesis} 144a5--9.

\begin{quotation}
— Εἰ οὖν ταῦτα οὕτως ἔχει,\\
 οἴει τινὰ ἀριθμὸν ὑπολείπεσθαι ὃν οὐκ ἀνάγκη \textit{εἶναι};\\
  — Οὐδαμῶς γε.\\
  — Εἰ ἄρα ἔστιν ἕν, \textit{ἀνάγκη καὶ ἀριθμὸν εἶναι.} \\
  — Ἀνάγκη

“Then if that is true, do you believe any number is left out,\\
 which does not necessarily be?”\\
 “By no means.” \\
 “Then if one is, number must also be.”\\
 “It must.” 142a2--5
\end{quotation}

Plato has been gradually approaching, and he now enunciates 
\begin{quotation}
\textit{statement [A]}. Every [eristic] number is generated in the One of the second hypothesis 144a4 
\end{quotation}

\textit{In the immediately succeeding passage 144a5--7 he describes is a transfer from numbers to units}, from hypothesis

\noindent \textbf{[A]} every number is generated in the One of the second hypothesis to the equivalent hypothesis that

\noindent \textbf{[B]} every part of the Being part can serve as \textit{unit} for the generation of a number in the One of the second hypothesis

\begin{quotation}
— Ἀλλὰ μὴν ἀριθμοῦ γε ὄντος\\
 πολλὰ ἂν εἴη\\
 καὶ πλῆθος ἄπειρον τῶν ὄντων·\\
 ἢ οὐκ ἄπειρος ἀριθμὸς πλήθει\\
 καὶ μετέχων οὐσίας γίγνεται;\\
 — Καὶ πάνυ γε.

“But if number is,\\
 there must be many, \\
 indeed an \textit{infinite multitude} (plethos aperon)\\
  of beings (onton);\\
   since\\
    \textit{number is infinite in multitude} and \\
    \textit{generated by participating in the Being}?”\\ “Certainly, it is.” 144a5--7
\end{quotation}

It is clear that Plato intends to generate every number using the parts generated by the anthyphairesis of One to Being. We interpret the statement
\begin{quotation}
‘the number $n$ is generated by \textit{participating in} the Being’
\end{quotation}
to mean
\begin{quotation}
‘the number $n$ \textit{is generated by the first $n$ parts} of the infinite sequence
$$\text{One, Being, One}_1, \text{Being}_1, \text{One}_2,\ldots$$
of the sequence of generated parts in the dyad One Being.’
\end{quotation}

   This kind of participation is similar to the participation of a sensible in an intelligible Being, equivalently to a true opinion about the intelligible being.

From the fact that number is infinite in multitude, Plato deduces that
\begin{quotation}
there is an \textit{infinite multitude of “beings”}.
\end{quotation}

It is nowhere explained what exactly is an element of the set of “beings”, but the terminology suggests that the “beings” are those parts of the Being, which are used as units for the generation of every number in the One \& Being, This interpretation is in agreement with the use of the term “beings” (149c2) in the enunciation of the plus one rule number of beings$=$contacts $+1$ in 148d--149d (examined in Section 3 below).

Thus, \textit{the first, eristic as will turn out to be, definition of number} in the One of the second hypothesis of the \emph{Parmenides} is the following: the number $n$ consists of the the first $n$ parts of the infinite sequence $$\text{One, Being, One}_1, \text{Being}_1, \text{One}_2,\ldots$$ of the sequence of generated parts in the dyad One Being.’

\subsection{The difficult passage \emph{Parmenides} 144a7--9}

We now concentrate on a seemingly innocuous statement:
\begin{quote}
— Οὐκοῦν [A] εἰ πᾶς ἀριθμὸς οὐσίας μετέχει 144a7--8,\\
καὶ [B] τὸ μόριον ἕκαστον τοῦ ἀριθμοῦ μετέχοι ἂν αὐτῆς 144a8--9;\\
— Ναί.
\end{quote}

\subsubsection{The standard rendering of 144a7--9 and its critique}

The standard rendering of 144a7--9 is:
\begin{quote}
``Then if [A] all number participates in the Being,
[B] every part of number will participate in the Being?"\\ ``Yes." 144a7--9
\end{quote}

Here are some variations of the standard rendering:
\begin{quotation}
\noindent [B] ``every part of number will partake of it?'' ``Yes.'' (H.N. Fowler, 1926 \cite{Fowler1926}),

\noindent [B] ``each part of number must have being also'' (Cornford, 1939 \cite{Cornford1939} p.142),

\noindent [B] ``each part of number also has a share of it? Yes'; (Allen, 1997~\cite{Allen1997} p.28),

\noindent [B] ``each part of number partakes of being as well'. (Palmer, 1999~\cite{Palmer1999}, p.229).
\end{quotation} 

Thus the standard rendering of the whole phrase is:
\begin{quote}
`If [A] every number partakes of existence, then [B] every part of number will partake of it'.
\end{quote}

Let us note however that \textit{`part of number'} makes no sense, at least in this context. For, what is exactly meant by part of number? If by part of number a smaller number is meant, then the standard rendering becomes nonsensical: if every number partakes of existence, then some numbers partake of existence. If by part of number a unit of the number is meant, this would be uniquely awkward. 

A weak point of this rendering is that it does not have sound internal logic. The natural syllogism is not
\begin{quotation}
\noindent [A] `if every number partakes of existence', 

\noindent [B] `then every part of number will partake of existence'
\end{quotation}
but rather
\begin{quotation}
\noindent `if some number partakes of existence', then every part of this number will partake of existence.
\end{quotation}

 Interestingly, Brumbaugh, 1961 \cite{Brumbaugh1961} p.99, is rendering [B], by changing the meaning of the first universal quantifier, as
\begin{quote}
`(144a7) `Now if the whole number partakes of existence, each part of number will also partake of it? Ar: Yes\ldots so that if the whole exists, so do the parts.';
\end{quote}
But Rickless', 2007~\cite{Rickless2007} p.346, analysis of the argument unwittingly brings out this very point; for, although he derives the standard interpretation of the consequent [B] (`each part of number partakes of being'), he does so, \textit{not} using the analogue of Plato's antecedent
\begin{quotation}
\noindent [A] if \textit{every (pas)} X partakes of being,
\end{quotation}
\textit{but}
\begin{quotation}
\noindent (P22) \textit{if X partakes of being},

\noindent then \textit{each (hekaston)} part of X partakes of being'.
\end{quotation}

     Thus the first universal quantifier `pas' (every) is omitted from the antecedent, while keeping the universal quantifier `hekaston' (`each') of the consequent.

\subsubsection{Our rendering of 144a7--9}
\begin{quote}
``if [A] every number participates in the Being [namely, if every finite collection of parts of Being generates a number with units the elements of this collection],\\
 then [B] every part of the Being will participate in number (144a8--9) [namely, then every part of the Being may be added as unit to any number to generate a greater number]?"\\
``Yes." 144a7--9
\end{quote}

In general, since every number can be generated in the One, and since the numbers generated in the One consist of units that are terms of the sequence $$\text{One}>\text{Being}>\text{One}_1>\text{Being}_1>\ldots>\text{One}_n>\text{Bing}_n>\ldots,$$ it follows that an infinite multitude (`plethos apeiron ton onton'), or in fact the infinite totality, of these terms, of the parts of the Being part serves as units (144a5--6), thus, statement A implies 
\begin{quotation}
\noindent statement B (`to morion hekaston tou arithmou metechoi an autes'). (144a8--9) every part of the Being may be added as unit to a number to generate a greater number.
\end{quotation}

\subsubsection{In support of our rendering}

We observe that in the \emph{Parmenides} (and in fact in all Platonic dialogues) the customary (if not quite exclusive) position of the object of the verb form `metechoi' is \textit{before} the verb, \textit{not after}.

Here is a complete list of the occurrences of ‘μετέχοι’ in all Platonic
dialogues.

\selectlanguage{polutonikogreek}
\begin{center}
%  \begin{supertabular}{|p{\cellwidth}|p{\cellwidth}|p{\cellwidth}|p{\cellwidth}|}
\begin{supertabular}{|p{\cellwidth}|p{\cellwidth}|p{\cellwidth}|p{\cellwidth}|}
  \hline
  \cc{Object to verb\\ \textit{before}} & & \cc{Object to verb\\ \textit{after}} & \cc{\itshape Passage}\\ \hline
\cc{τοῦ κινδύνου τε\\ καὶ τῆς τύχης} & \cc{μετέχοις} & & \cc{\textit{Euthydemus} 279e7--280a1}\\ \hline
\cc{εἰ ὀλίγοι αὐτῶν} & \cc{μετέχοιεν} & & \cc{\textit{Protagoras} 322d3}\\ \hline
\cc{οὐδὲν πλέον ἀναρμοστίας\\
 οὐδὲ ἁρμονίας}
&
\cc{μετέχοι ἄν} &
 &
\cc{\textit{Phaedo} 93e4--5}\\\hline
\cc{ἆρ' ἄν τι πλέον\\ κακίας ἢ ἀρετῆς}
&
\cc{μετέχοι\\ ἑτέρα ἑτέρας,}
&
 &
\cc{\textit{Phaedo} 93e7--8}\\\hline
\cc{ἀρετῆς} &
\cc{μετέχοι\\ πόλις} &
 &
\cc{\textit{Politeia} 432b4}\\\hline
\cc{μέρους} &
\cc{ἂν μετέχοι,} &
 &
\cc{\textit{Parmenides} 131c6}\\\hline
\cc{εἴτε εὐθέος σχήματος\\
εἴτε περιφεροῦς}
&
\cc{μετέχοι} &
 &
\cc{\textit{Parmenides} 137e5--6:}\\\hline
\cc{Καὶ σχήματος δή τινος,\\ ὡς ἔοικε,}
&
\cc{τοιοῦτον ὂν\\ μετέχοι ἂν\\ τὸ ἕν}
&
 &
\cc{\textit{Parmenides} 145b3--4}\\\hline
\cc{ἀριθμοῦ πλείονος} &
\cc{ἂν μετέχοι} &
 &
\cc{\textit{Parmenides} 153a4--5}\\\hline
\cc{τοῦ αὐτοῦ} &
\cc{μετέχοι τε καὶ\\ οὐ μετέχοι}
&
 &
\cc{\textit{Parmenides} 155e11}\\\hline
\cc{τοῦ ὅλου τε καὶ ἑνὸς} &
\cc{μετέχοι} &
 &
\cc{\textit{Parmenides} 157e3}\\\hline
 &
\cc{Μετέχοι δέ γε ἂν} &
\cc{τοῦ ἑνὸς} &
\cc{\textit{Parmenides} 158a3--4}\\\hline
\cc{οὐ γὰρ ἂν τοῦ ἑνὸς} &
\cc{μὴ μετέχοι} &
 &
\cc{\textit{Parmenides} 158c4}\\\hline
 &
\cc{μετέχοι ἂν τἆλλα} &
\cc{τοῦ ἑνός,} &
\cc{\textit{Parmenides} 159d1}\\\hline
\cc{μηδ' ἑνὸς} &
\cc{μετέχοι} &
 &
\cc{\textit{Parmenides} 159e7}\\\hline
\cc{ἑνὸς ἂν τοῦ ἑτέρου εἴδους}
&
\cc{μετέχοι} &
 &
\cc{\textit{Parmenides} 160a1--2}\\\hline
\cc{οὐσίας} &
\cc{μετέχοι} &
 &
\cc{\textit{Parmenides} 164a1}\\\hline
\cc{πασῶν} &
\cc{ὅσαι κινήσεις ἔσοιντο\\ μετέχοι}
&
 &
  \cc{\textit{Timaeus} 44d8}\\ \hline
\end{supertabular}
\end{center}
%%%%%%%%%%%%%%%%%%%%%%%%%%%%%%%
\selectlanguage{english}
%%%%%%%%%%%%%%%%%%%%%%%%%%%%%%% 
\noindent Thus with the statement under question:

καὶ τὸ μόριον ἕκαστον τοῦ ἀριθμοῦ μετέχοι ἂν αὐτῆς \textit{Parmenides} 144a8--9,

\noindent in the customary rendering \textit{the object} to the verb ‘metechoi’
appears \textit{after the verb}:

\begin{center}
\tablefirsthead{}
\tablehead{}
\tabletail{}
\tablelasttail{}
\begin{supertabular}{|p{\cellwidth}|p{\cellwidth}|p{\cellwidth}|p{\cellwidth}|}
\hline
 &
\cc{τὸ μόριον ἕκαστον\\ 
τοῦ ἀριθμοῦ\\
μετέχοι ἂν}
 &
\cc{αὐτῆς} &
\cc{\textit{Parmenides} 144a9}\\ \hline
\end{supertabular}
\end{center}
in which case the rendering would be

each part (to morion hekaston) of number (tou arithmou) 

participates (metechoi an) 

in the Being part (autes, [namely, ousias]),

\noindent\textit{while} we suggest that \textit{the object} to the verb ‘metechoi’ is \textit{before the verb},

\null\ \ 
\begin{center}
\tablefirsthead{}
\tablehead{}
\tabletail{}
\tablelasttail{}
\begin{supertabular}{|p{\cellwidth}|p{\cellwidth}|p{\cellwidth}|p{\cellwidth}|}
\hline
 &
\cc{τὸ μόριον ἕκαστον}
& & \\ \hline
\cc{τοῦ ἀριθμοῦ} &
 & &
\\ \hline
 &
\cc{μετέχοι ἂν\\
αὐτῆς} & & \cc{\textit{Parmenides} 144a8--9}
\\\hline
\end{supertabular}
\end{center}
in which case the rendering would be ‘tou arithmou metechoi’:

each part of the Being part (to morion hekaston autes, [namely, ousias]) 

participates (metechoi an) 

in number (tou arithmou) [as its unit].

\noindent As indicated by the table above, the (heavy, 16 to 2) odds are in favor of our interpretation.

% We suggest that the object to the verb ‘metechoi’ is 

% not after, in which case it would be 

% ‘metechoi an autes’ [namely, ‘ousias’]: 

% ‘each part of number participates in the Being part’. 

% as in the customary rendering; 

% Our interpretation:

% \begin{center}
% \tablefirsthead{}
% \tablehead{}
% \tabletail{}
% \tablelasttail{}
% \begin{supertabular}{|p{\cellwidth}|p{\cellwidth}|p{\cellwidth}|p{\cellwidth}|}
% \hline
%  &
% \centering τὸ μόριον ἕκαστον &
% \makebox[2cm]{} &\makebox[2cm]{}
% \\\hline
% \centering τοῦ ἀριθμοῦ &
% \centering μετέχοι ἂν &
%  &
% \\\hline
%  &
% \centering αὐτῆς &
%  &
% \\\hline
% \end{supertabular}
% \end{center}
% But we suggest that the object to the verb ‘metechoi’ is 

% before, and then it is ‘tou arithmou metechoi’:

% ‘each part of the Being part participates in number as a unit’.

\subsection{The third stage 144b1--c2 notes that a linguistic comparison of 144b1--c2 with the \emph{Philebus} 56d--e passage reveals a fatal criticism concludes that in the way the numbers have been generated, and the parts of the Being are ``most numerous" (pleista), rendered by us as ``infinite in number"}

\subsubsection{\emph{Parmenides} 144b1--c2}
\begin{quote}
Ἐπὶ πάντα ἄρα \textit{πολλὰ ὄντα} ἡ οὐσία νενέμηται\\
καὶ οὐδενὸς ἀποστατεῖ \textit{τῶν ὄντων,}\\
οὔτε \textit{τοῦ σμικροτάτου} οὔτε \textit{τοῦ μεγίστου;}\\
ἢ τοῦτο μὲν καὶ ἄλογον ἐρέσθαι;\\
πῶς γὰρ ἂν δὴ οὐσία γε\\
\textit{τῶν ὄντων} του\\
ἀποστατοῖ;\\ 
— Οὐδαμῶς.\\
— Κατακεκερμάτισται ἄρα\\
ὡς οἷόν τε \textit{σμικρότατα} καὶ \textit{μέγιστα} καὶ\\
πανταχῶς ὄντα,\\
καὶ μεμέρισται\\
πάντων\\
μάλιστα,\\
καὶ ἔστι μέρη ἀπέραντα τῆς οὐσίας.\\
— Ἔχει οὕτω.\\
— \textit{Πλεῖστα} ἄρα ἐστὶ\\
τὰ μέρη αὐτῆς.\\
— \textit{Πλεῖστα} μέντοι. \emph{Parmenides} 144b1--c2

``The Being (he ousia) then, is \textit{distributed (nenemetai)}\\
over all beings (onta),\\
which are many, and [the Being] is not wanting (apostatei)\\
in any of the beings (ton onton),\\
from \textit{the smallest (tou smikrotatou)}\\
to \textit{the greatest (tou megistou)}?\\
Indeed, is it not \textit{alogon (alogon)}\\
even to ask that question? \\
For how can the Being be wanting (apostatoi)\\
in any of the beings (onta)?"\\
``It cannot by any means." \\
Then the Being is \textit{fragmented (katakekermatistai)} into \\
\textit{the smallest (smikrotata)}\\
and \textit{the greatest (greatest)} beings\\
in all places (pantachos),\\
nothing else is so much divided (memeristai), [144c] and the parts of the Being\\
are \textit{infinite in multitude} (aperanta)."\\
``That is true."\\
``Hence (ara), the parts of the Being\\
are the \textit{most numerous} (pleista) of all."\\
``Yes, they are\\
the \textit{most numerous} (pleista)." 144b1--c2
\end{quote}

The first part of the passage is talking about the \textit{beings}, the last about the\textit{ parts of the Being}. They are the same entities, but \textit{being} refers to \textit{a part of the Being which serves as a unit} for the numbers that are being generated in the One of the Second Hypothesis. \textit{So a being is a part of the Being with the special, additional property that it serves as the unit for a number.}

\subsubsection{Plato's criticism of the definition/generation of eristic numbers inferred by linguistic comparison of the \emph{Parmenides} 144b1--c2 with the \emph{Philebus} 56d4--e6}

We might sense that the transfer from numbers to units has resulted in an imperceptible movement from wholehearted acceptance of statement A to an essential \textit{neutrality} regarding acceptance of its equivalent statement B. In order to realize that 144b1--c2 does in fact constitute \textit{criticism} and to understand the reason for the \textit{criticism}, that turns into \textit{explicit rejection} a little later, in 144d5--7, of the statement that [B] every part of the Being part can serve as unit for the numbers generated in the One, it will be helpful to have in mind Plato's fundamental distinction between the two kinds of numbers, the numbers used by the many, the inferior numbers so to speak, and the numbers of the philosophers, the superior dialectic numbers. 

The distinction is described forcefully in the \emph{Philebus} 56d--e:
\begin{quote}
ΠΡΩ. Ἀριθμητικὴν φαίνῃ μοι λέγειν καὶ ὅσας μετὰ ταύτης τέχνας ἐφθέγξω νυνδή.\\
ΣΩ. Πάνυ μὲν οὖν. ἀλλ', ὦ Πρώταρχε, ἆρ' οὐ διττὰς αὖ καὶ ταύτας λεκτέον; ἢ πῶς;\\
ΠΡΩ. Ποίας δὴ λέγεις;\\
ΣΩ. Ἀριθμητικὴν πρῶτον ἆρ' οὐκ

\begin{center}
\tablefirsthead{}
\tablehead{}
\tabletail{}
\tablelasttail{}
\begin{supertabular}{|p{\dcellwidth}|p{\dcellwidth}|}
\hline
\dcc{ἄλλην μέν τινα τὴν τῶν πολλῶν φατέον,} &
\\ \hline
 &
\dcc{ἄλλην δ' αὖ τὴν τῶν φιλοσοφούντων;}\\ \hline
\end{supertabular}
\end{center}

 ΠΡΩ. Πῇ ποτε διορισάμενος οὖν 

\begin{center}
\tablefirsthead{}
\tablehead{}
\tabletail{}
\tablelasttail{}
\begin{supertabular}{|p{\dcellwidth}|p{\dcellwidth}|}
\hline
\dcc{ἄλλην,} &
\\\hline
 &
\dcc{τὴν δὲ ἄλλην}\\ \hline
\end{supertabular}
\end{center}

 θείη τις ἂν ἀριθμητικήν; 

 ΣΩ. Οὐ σμικρὸς ὅρος, ὦ Πρώταρχε. 

 \selectlanguage{polutonikogreek}
\begin{center}
\tablefirsthead{}
\tablehead{}
\tabletail{}
\tablelasttail{}
\begin{supertabular}{|p{\dcellwidth}|p{\dcellwidth}|}
\hline
\dcc{οἱ μὲν γάρ που μονάδας ἀνίσους\\ καταριθμοῦνται τῶν περὶ ἀριθμόν,\\
οἷον στρατόπεδα δύο καὶ βοῦς δύο\\
καὶ δύο τὰ σμικρότατα\\ ἢ καὶ τὰ πάντων μέγιστα·} &
\\ \hline
 &
\dcc{οἱ δ' οὐκ ἄν ποτε αὐτοῖς συνακολουθήσειαν,\\
εἰ μὴ μονάδα μονάδος ἑκάστης τῶν μυρίων\\ μηδεμίαν ἄλλην
ἄλλης διαφέρουσάν\\ τις θήσει.}\\ \hline
\end{supertabular}
\end{center}

 ΠΡΩ. Καὶ μάλα εὖ λέγεις 

οὐ σμικρὰν διαφορὰν τῶν περὶ ἀριθμὸν τευταζόντων,

ὥστε λόγον ἔχειν δύ' αὐτὰς εἶναι. 

\selectlanguage{english}

‘\textit{Socrates} Are there not two kinds of arithmetic, 

\begin{center}
\tablefirsthead{}
\tablehead{}
\tabletail{}
\tablelasttail{}
\begin{supertabular}{|p{\dcellwidth}|p{\dcellwidth}|}
\hline
\dcc{that of the many people} &
\\ \hline
 &
\dcc{and that of philosophers?}\\ \hline
\end{supertabular}
\end{center}
\begin{center}
Protarchus\\ How can be distinguished
\end{center}
\begin{center}
\tablefirsthead{}
\tablehead{}
\tabletail{}
\tablelasttail{}
\begin{supertabular}{|p{\dcellwidth}|p{\dcellwidth}|}
\hline
\dcc{one kind of arithmetic} &
\\ \hline
 &
\dcc{from the other?}\\ \hline
\end{supertabular}
\end{center}

\textit{Socrates} The distinction is no small one, Protarchus.

\begin{center}
\tablefirsthead{}
\tablehead{}
\tabletail{}
\tablelasttail{}
\begin{supertabular}{|p{\dcellwidth}|p{\dcellwidth}|}
\hline
\dcc{For some arithmeticians reckon unequal units,\\
 for instance,
 two armies and\\
 two oxen and
 two the smallest (ta smikrotata) 
 or\\ the greatest of all 
(ta panton megista);} &
\\ \hline
 &
\dcc{whereas others refuse to agree with them\\ unless
 each of countless units\\ 
 is declared to differ not at all\\ 
 from each and every other unit.}\\ \hline
\end{supertabular}
\end{center}
\begin{center}
 \textit{Protarchus} You are certainly quite right in saying that\\
 here is a great difference between the devotees of arithmetic,\\
so it is reasonable to assume that it is of two kinds.’ \textit{Philebus} 56d4--e6
\end{center}
\end{quote}

%A `superior', \textit{dialectic} number, a number used by the philosophers, consists of equal units, an `inferior', \textit{eristic} number, a number used by the many, consists of unequal units.

\subsubsection{The description of the inferior, eristic, numbers in the \emph{Philebus} 56d9--e1}

\begin{quote}
`οἱ μὲν γάρ που \textit{μονάδας ἀνίσους} καταριθμοῦνται τῶν περὶ ἀριθμόν,\\
οἷον στρατόπεδα δύο καὶ βοῦς δύο\\
καὶ \textit{δύο τὰ σμικρότατα} ἢ καὶ \textit{τὰ πάντων} μέγιστα' \\
For some arithmeticians reckon \textit{unequal units}, for instance,two armies and two oxen and two the \textit{smallest (ta smikrotata) or the greatest of all (ta panton megista);} \emph{Philebus} 56d9--e1
\end{quote}

should be compared with the criticism of Statement B leveled in the \emph{Parmenides} 144b--c:

\begin{quote}
οὐδενὸς ἀποστατεῖ τῶν ὄντων, \textit{οὔτε τοῦ σμικροτάτου οὔτε τοῦ μεγίστου (144b2--3); }\\
κατακεκερμάτισται ἄρα ὡς οἷόν τε \textit{σμικρότατα καὶ μέγιστα} καὶ πανταχῶς ὄντα (144b4--6); \\
οὐκ ἀπολειπόμενον \textit{οὔτε σμικροτέρου οὔτε μείζονος μέρους} οὔτε ἄλλου οὐδενός (144c6--7).
\end{quote}

It is understood that the problem with Statement B is that the generation of all the numbers in the One is possible only by employing \textit{unequal units}, and hence, the numbers generated in the one are the numbers of the many, the inferior eristic numbers and \textit{not the numbers of the philosophers.}

\subsection{Plato notes that Statement B implies Statement C (144c1--2) and finally, in view of the most crucial statement 144d4--5, rejects statement C (pleista) and the definition of numbers as eristic (144d5--7)}

\begin{quote}
\textit{Statement C}\\
— \textit{Πλεῖστα} ἄρα ἐστὶ τὰ μέρη αὐτῆς.\\
— \textit{Πλεῖστα} μέντοι. 144c1--2 \\
The parts of the Being part are \textit{most numerous} (`pleista')\\
``Yes, they are most numerous."
\end{quote}

What is precisely meant by the Statement C,
\begin{quotation}
that the parts of the Being part are \textit{`pleista'?}
\end{quotation}

It \textit{cannot just mean} that 
\begin{quotation}
the parts of the Being part are just \textit{`infinite in multitude',}
\end{quotation}
since:
\begin{itemize}
    \item the statement that \textit{the parts of Being are infinitely in multitude} has been proved in Section 3, and certainly holds true, and
    \item Statement C that \textit{the parts of Being are `pleista'} is \textit{explicitly rejected} in 144d5--7 below, as not true:
\end{itemize}

\begin{quote}
— Καὶ μὴν τό γε μεριστὸν πολλὴ ἀνάγκη εἶναι \textit{τοσαῦτα ὅσαπερ} μέρη.\\
— Ἀνάγκη. 144d4--5 \\
— Οὐκ ἄρα ἀληθῆ ἄρτι ἐλέγομεν λέγοντες ὡς \textit{πλεῖστα} μέρη ἡ οὐσία νενεμημένη εἴη. 144d5--7
\end{quote}

Passage 144d4--7 is rendered and discussed in Section 6.2, below. In order to appreciate the meaning in statement C of the word \textit{`pleista'}, there should be a careful distinction between 

\textit{infinity in multitude} 

and \textit{infinity in number.}

\textit{Infinity in multitude} 

is a weak statement, 

in modern terms, 

essentially a set-theoretic one. 

For the sequence of parts of the One 

\textit{infinity in multitude} 

has been proved in Section 3.

But, given the Platonic strictures 

about proper, philosophic numbers, 

\textit{infinity in number} 

is a much stronger statement. 

It implies that, in some unspecified as yet way, 

these parts are \textit{equalized.} 

The statement C, involving \textit{‘pleista’},

 is precisely the \textit{infinity in number}.

We shall return to statement C in Section 6, below, where the dialectical numbers are introduced.

In conclusion the rejection of Statement C, and hence of Statements B and A, is due to the fact that the infinite multitude of parts of the One, proved in Section 3, in no way implies the generation of all number in the One of the second hypothesis, precisely because these parts are \textit{mutually unequal}, and some of these parts are the smallest, while some others are the greatest. At this point the situation might appear to be hopeless: all parts in the One are mutually unequal, hence, in fact we have a problem generating in the One even the number Two!
Thus, the \textit{eristic definition}/generation of all numbers with units the elements of the infinite sequence of remainder parts of the anthyphairesis of the dyad One, Being:

\begin{center}
\tablefirsthead{}
\tablehead{}
\tabletail{}
\tablelasttail{}
\begin{supertabular}{|p{\cellwidth}|p{\cellwidth}|p{\cellwidth}|p{\cellwidth}|}
\hline
\cc{is first introduced\\  
  with enthusiastic\\
  acceptance\\ (143d5--144a5,\\
 statement A\\ section 4.1),} &
 &
 &
\\ \hline
 &
 \cc{then is accepted\\  somewhat neutrally\\
  (144a5--9,\\ statement B\\
section 4.2--3),} &
 &
\\\hline
 &
 &
\cc{next is criticized\\
for having\\
unequal units\\
(144b1--c2,\\
Statement C\\
“pleista”\\
 section 4.4--5),} &
\\\hline
 &
 &
 &
\cc{and finally ends\\ with being\\ \textit{explicitly rejected}\\
(144d5--7\\
“pleista”\\
 section 4.6).}\\ \hline
\end{supertabular}
\end{center}

\section{Anthyphairetic periodicity of the dyad One Being in the \emph{Parmenides} Second Hypothesis}

%\begin{enumerate}
% \item [1)]
A contact/hapsis, introduced in 148d5--149d7, is identified. by the plus one rule, a rule of musical origin. with a ratio of successive anthyphairetic remainders of the dyad One Being (Section 5.4-5).

% \item [2)]
by the most crucial statement in 138a6--7 the hapseis/contacts/ratios of successive remainders form a circle, the philosophic analogue of the Logos Criterion and of periodic anthyphairesis for the dyad One \&\ Being (Section 5.6).

% \item [3)]
the presence of the One in the Being 144c2--d4, achieved by the contacts/hapseis 138a3--7 forming a circle, is a consequence of the Logos Criterion and anthyphairetic periodicity (Section 5.1-3, 5.6).

% \item [4)]
the structure of the One, the paradigmatical intelligible Being in the \emph{Parmenides}, is, described as True Opinion plus Logos, equivalently Name plus Logos in 155d3--e1, exactly as in the \emph{Theaetetus, Sophist, Meno}, and is the philosophical analogue of periodic anthyphairesis (Section 5.7).
%\end{enumerate}

The fundamental conclusion \textit{The dyad One\&Being in the second hypothesis of the Parmenides in fact satisfies the philosophical analogue of periodic anthyphairesis} (144c2--d4, 138a3--7, 148d--149d) is established with the following non-obvious steps:

%%%%
5.1.\\
144c2-d4 The One is [present] in the parts of the Being\\
145b6-e6 The One is [present] in the Other (and also the One is in itself).\\
It is stated that a part of the One will be present in a place /part of the Being, but otherwise it is not clear how the One will be in the Being.

\medskip

5.2. First Connection \\
138a3-7 explains the statement “the One is in the Other” (138a3) \\
for any entity that may be called One, in particular \\
not only for the One of the First Hypothesis \\
(which does not satisfy the statement “the One is in the Other”), \\
but also for the One of the Second Hypothesis \\
(which does satisfy the statement “the One is in the Other”).

\medskip

5.3. The statement “the One is in the Other” in terms of circularity of contacts\\
138a3-5 The statement “the One is in the Other” \\
is equivalent to the statements\\
“The One is in the Other cyclically”; and\\
“the One possesses contacts by means of many of its parts with many places/parts of the Being.”\\
Plato, interpreting the crucial statement\\
“the One is in the Other/Being”,\\
introduces for the first time a circular/periodic process of containment/presence;\\
and contacts (hapseis) by many parts of the One with many places/parts in the Being.\\
But the meaning of both is unclear.

\medskip

5.4-5. Second Connection to 148d5-149d7 where “contacts” and “dialectical number” are defined: \\
The sequence of terms in 148d5-149d7 coincides with the sequence of the anthyphairetic remainders\\

\medskip

5.4. Definition of contact (hapsis): \\
contact is a ratio between consecutive terms/parts, hence a ratio between consecutive anthyphairetic remainders of the philosophical anthyphairesis of the dyad One, Being.

\medskip

5.5. Definition of dialectical number:\\
dialectical number of terms=contacts+1 (the plus one rule).\\
This replaces the eristic definition of number.

\medskip

5.6. Back to 138a3--7 for the establishment of the periodicity of the anthyphairesis of the One to the Being

With the clarification of “contact” obtained in 5.4,\\
Plato’ connection in 138a5-7 of contacts 138a5 with circularity138a3-4 becomes clear:\\
the One possesses contacts, with many places/parts of the Being, forming a circle.

Since the contacts/ratios are in the sequence of anthyphairetic remainders,\\
as it has been proved in Section 3, it is clear that the statement\\
 “\textit{the contacts form a circle}”\\
states the Logos Criterion for anthyphairetic periodicity:\\
there is $k$, such that $\text{One}/\text{Being}=\text{One}_k/\text{Being}_k$.\\
Thus the dyad One Being satisfies the philosophical analogue of\\
\textit{the Logos Criterion and anthyphairetic periodicity.}

\medskip

Now it is clear how the One is in the Other:\\
$\text{One}_k$ is a is in Being.\\
(and also how the One is in itself: $\text{One}_k$ is in One).\\
Thus the knowledge of the One of the second hypothesis is obtained as Name plus Logos,\\
and the One forms with its part the Being \\
the philosophical analogue of a dyad in philosophic periodic anthyphairesis.\\
Note that, as we have seen in Section 2, the anthyphairesis of the One of the second hypothesis is
\textit{not abridged}, like that of “the Angler” and “the Sophist” in the \textit{Sophist},
\textit{but is complete}, providing, because of this, additional confirmation of our interpretation.

\medskip

5.7.  \textit{The One of the Second Hypothesis in the Parmenides is the paradigmatical intelligible Being, whose Knowledge is given by Name plus Logos}\\
The One of the Second Hypothesis \\
--forms with its part Being the philosophical analogue of a dyad in infinite anthyphairesis \\
(as shown in Section 3);  \\
--in such way that the One is in the Other/Being, and thus \\
--satisfies with its part Being the Logos Criterion for anthyphairesic periodicity\\
(as explained in 5.6 above by means of 142b1--144d4, 145b5--e6, 138a3--7, 148d--149d),\\
-- possesses (155d3--e1) Knowledge (episteme, 155d6), given by Name plus Logos (155d8--e1),  \\
-- is the paradigmatical intelligible Being (155d3--6). 
%%%%

The details follow.

\subsection{The presence of the One in the Being/Other 144c2--d4, 145b6--e6}

\subsubsection{Plato realizes now at 144c2 the difficulty and in fact the impossibility of his earlier aim, which was the attempt to define numbers in the One of the second hypothesis, by exploiting only the infinite anthyphairesis of the dyad One and Being, and no any additional structure. But we do not forget that we still expect that the dyad One and Being will somehow turn out to have not just omit infinite but infact periodic anthyphairesis. Exactlty at the position 144c2 Plato is making a new start, introducing the property: ``the One is present in the Being".}

\hspace*{2.25em}\begin{minipage}{0.8\textwidth}
\begin{enumerate}[itemsep=-4ex]
\item[—] Τί οὖν; ἔστι τι αὐτῶν ὃ ἔστι μὲν μέρος τῆς οὐσίας, \textit{οὐδὲν μέντοι μέρος;}\\
\item[—] Καὶ πῶς ἄν [τοι] τοῦτο γένοιτο;\\
\item[—] Ἀλλ' εἴπερ γε οἶμαι ἔστιν, ἀνάγκη αὐτὸ ἀεί, ἕωσπερ ἂν ᾖ, ἕν γέ τι εἶναι, μηδὲν δὲ ἀδύνατον.\\
\item[—] Ἀνάγκη.\\
\item[—] Πρὸς ἅπαντι ἄρα [ἑκάστῳ] τῷ τῆς οὐσίας μέρει \textit{πρόσεστιν} τὸ ἕν,\\
οὐκ ἀπολειπόμενον\\
οὔτε σμικροτέρου οὔτε μείζονος μέρους\\
οὔτε ἄλλου οὐδενός.\\
\item[—] Οὕτω.\\
\item[—] Ἆρα οὖν ἓν ὂν πολλαχοῦ \textit{ἅμα ὅλον} ἐστί; τοῦτο ἄθρει.\\
\item[—] Ἀλλ' ἀθρῶ καὶ ὁρῶ ὅτι \textit{ἀδύνατον.}\\
\item[—] \textit{Μεμερισμένον} ἄρα, εἴπερ \textit{μὴ ὅλον}· ἄλλως γάρ που οὐδαμῶς \textit{ἅμα} ἅπασι τοῖς τῆς οὐσίας μέρεσιν \textit{παρέσται ἢ μεμερισμένον}. \\
\item[—] Ναί. 144c2--d4 \\
``Well, is there any one of them which is a \textit{part} of the Being (meros tes ousias), but is \textit{not} a One part (ouden meros)?"\\
``How could that be?"\\
``But \textit{if it is part of the Being (estin)}, it must, I imagine,\\
as long as \textit{it is part of the Being (heosper an ei)},\\
be a One part; it cannot be not a One part.\\
``That is inevitable."\\
``Then \textit{Oneness is an attribute (prosestin to hen) of every part of the Being} and is not wanting to either \textit{a smaller (smikroterou) or greater (meizonos)} or any other part."\\
``True."\\
``Can the initial One part be \textit{in many places (pollachou) at once (hama)} and still be a whole (holon)? Consider that question."\\
``I am considering and I see that it is impossible."\\
``Then the initial One part \textit{is divided into parts (menerismenon)}, if it is not a whole; for it cannot be \textit{attached (parestai) to all the parts of the Being at once} (hama) unless it is \textit{divided (memerismenon)}."\\
``I agree." 144c2--d4
\end{enumerate}
\end{minipage}

The \textit{presence} (parestai', 144d4) of the One in the Being, is examined; it is realized that the presence of the One must be in the many parts of the Being (\textit{`pollachou'}, 144d1), and the question is posed whether the presence of the One in the parts of the Being \textit{simultaneously} (hama holon', 144d1) will be either as \textit{a whole}, or as \textit{a part}.
The first leg of the question is found to be impossible (`adunaton', 144d2), and hence rejected, and only the second possibility is left, resulting in a divided (\textit{`memerismenon'}, 144d2,4) One.

\subsubsection{The \textit{presence} of the One in the Being is in need of clarification. Let us first note that two quite similar properties of the One of the second hypothesis appear in 145b6--e6:}

\hspace*{1em}\begin{minipage}{0.8\textwidth}
\begin{quote}
τὸ ἓν ἀνάγκη αὐτό τε ἐν ἑαυτῷ εἶναι καὶ ἐν ἑτέρῳ \\
the One is in itself and the One is in the Other.
\end{quote}
\end{minipage}

\bigskip

Since the Being is Other to the One, it follows that the statements\\
``the One is in the Being" and ``the One is in the Other" are equivalent.\\
But we are still in the dark as to the meaning of these statements, since the One is a part greater than the Being/Other.

\subsection{The passages on the presence of the One in the Being/Other 144c2--d4, 145b6--e6 are similar to the passage on the presence of the One in the Other 138a3}

We are told in 144c2--d4 and 145b6--e6 that the statement ``the One is in the Being/Other" holds true for the One of the second hypothesis, but we find no reasonable explanation for this claim. On the other hand, we are told earlier in 138a3--7 that the same statement does not hold true for the One of the first hypothesis, and there is an explanation provided for this failure to hold. In trying to show why the partless \textit{One of the first hypothesis} does \textit{not} enjoy the property of being \textit{in the Other}, Plato feels obliged to explain what is really meant by the statement \textit{`the One is in the Other'}. This is done precisely in 138a3--7: it turns out that the presence of the One in the Other/Being occurs in passage 138a3--7.

\subsubsection{The similarities between ``the One is in the Being" 144c2--d4, ``the One is in the Other" 145b6--e6 and ``the One is in the Other" 138a3--7}

In order to establish the relevance of 138a3--7 to 144c2--d4 we make the following remarks.

\medskip

\noindent\textit{First, }

\smallskip

\noindent there is no question that the statement 

‘the One is \textit{in} the Other’ in 138a3--7

\noindent is equivalent to

‘the One is \textit{present} in the Other’,

\noindent and similar to the statement

“the One is \textit{present} in the Being” in 144c2--d4.

Indeed, the statement ‘the One in the Other’ (138a3--7, 145c7--e3) 

is expressed by \textit{‘eneie’} in 138a4, \textit{‘eneinai}’ in 139a7, and \textit{‘enestai’} in 145d4. 

‘\textit{paresti}’ in 144c2--d4 is equivalent to. ‘\textit{enesti}’

Thus, both passages 138a3--7 and 144c2--d4 deal with the presence of the One in an
intelligible entity,

the Other and the Being, respectively.

\noindent \textit{Second,}

\smallskip

the term \textit{Other} of the One, 
employed in 138a3--7, 
applies to 
the \textit{Being part} (143b1--3), 
and in fact,

only to the Being part, 
about which 144c2--d4 is concerned.

\begin{quote}
Cf. `ἄλλο τι ἕτερον μὲν ἀνάγκη τὴν οὐσίαν αὐτοῦ εἶναι, ἕτερον δὲ αὐτό [[\text{τὸ ἕν}]], εἴπερ μὴ οὐσία τὸ ἕν, ἀλλ' ὡς ἓν οὐσίας μετέσχεν. 143b1--3
\end{quote}

The two parts One and Being form a dyad of `hetera' and `alla' to each other.

Thus, passage 138a3--7, that describes how the One should be in the Other, in effect describes how the One should be in the Being. 

The Proposition that 

`the One is in the Other', 

for the One of the second hypothesis, appears in 145c7--e3. 

The passage is wholly concerned with rejecting all the naïve types of presence and containment (the presence of the One in itself, in one, in several and in all its parts),

since in fact

the One is certainly \textit{not contained} as a part of the Being part,

\textit{but} on the contrary the One part \textit{contain}s as a part the Being part,

ending with the statement that the One, in view of these naïve rejections, must be in the Other (`en allo', `en heteroi,145e2--3),

thus, in effect, referring the explanation to the earlier passage 138e3--7.

\medskip

\textit{Third,}

the passage 138a3--7 is \textit{the only} other passage in the \emph{Parmenides},

besides the two passages 142c2--d4 [and 130e4--131e7],

in which the term \textit{`pollachou}' appears.

These are the five occurrences of `pollach-' in the \emph{Parmenides}.

The presence of the intelligible One in the intelligible Being 144c2--e4: pollachou 144d1,

the presence of the intelligible One in the intelligible Other 138a3--7: pollachou 138a5, pollachei kukloi 138a6, and

the presence of the intelligible Being in the sensibles 130e4--131e7:131b4,131b7.

The meaning of \textit{`pollachou'} in 138e3--7

is seen to be the same as

in 144c2--d4,

namely `in many parts of the Being/Other'.

\medskip

\textit{Fourth,} 

the meaning of \textit{`pollois'} in 138a5 refers to
the many parts of the One, by means of which the One is present
in the many parts/places of the Other/Being.
Thus, the One in 138a3--7is present in the Other
by means of many of its parts,
and thus as a \textit{`memerismenon'} One,
exactly as with the One in 144c2--d4.

The careful comments by Proclus, \emph{Commentary to Plato's Parmenides} (see Appendix 5.2.3, below) confirm this interpretation of `pollois'.

Thus, the crucial passage 138a3--7 deals with the situation in which
the One,
consisting of `polla' (and thus `memerismenon'),
is present in many places-parts (`pollachou') of the Being.
Passage 138a3--7 then deals precisely with the situation of
the divided One present in
the many places-parts of the Being,
namely precisely with the situation of the passage 144c2--d4, the one we are trying to decipher.

To the best of my knowledge the crucial close connection between these two passages, a connection that, as we shall see, will lead to a complete understanding of the second hypothesis, has not been noticed by previous scholars.

\medskip

\textbf{Appendix to 5.2.1.} Proclus, \emph{Commentary in Plato's Parmenides} 1141, 5--7 \& 20--24.

The comments by Proclus, \emph{Commentary in Plato's Parmenides} 1141, 5--7 \& 20--24 confirm that the term `pollois' in the passage \emph{Parmenides} 138a3--7 refers exactly to the the presence of the One, by means of many of its parts (`pollois'), in the many parts of the Other, precisely as the divided (`memerismenon') One in the \emph{Parmenides} 144c2--d4 passage.

\begin{quote}
τὸ \textit{μὲν}\\
\textit{ἐν τόπῳ}\\
δεῖ εἶναι \textit{πολλὰ} ὃν\\
καὶ \textit{πολλοῖς} ἑαυτοῦ [[του ενός]]\\
\textit{τοῦ περιέχοντο}ς [[του τόπου, του άλλου]]\\
\textit{ἅπτεσθαι}, \\
{}[the One] being in (en) a place (topoi)\\
must involve being many (polla) [parts]\\
and touching (haptesthai) its container (tou periechontos)\\
by means of many (`pollois') [parts]\\
\emph{Proclus Commentary in Plato's Parmenides 1141, 5--7}\\ translated by Morrow \& Dillon, \cite{Proclus1987} p.497
\end{quote}

\begin{quote}
`[σημεῖον]\\
\textit{ἐν ἄλλῳ} [grammei]\\
ἐστὶ \ldots\\
\textit{πανταχόθεν}\\
\textit{περιέχεται}\\
ὑπὸ τῆς γραμμῆς\\
καὶ \textit{πολλοῖς} ἑαυτοῦ μέρεσιν\\
\textit{ἅπτεται} \\
τῆς γραμμῆς. \\
`[the point] is\\
in the other [line] \ldots \\
contained\\
on all sides\\
by the line,\\
nor does it contact the line in many of its points [by means of many of its parts] \emph{Proclus' Commentary in Plato's Parmenides 1141, 20--24} [translated by Morrow \& Dillon, \cite{Proclus1987} p.498]
\end{quote}

\subsubsection{The passage 145c7--e3 also leads to passage 138a3--7}

But the Proposition `the One is in the Other', for the One of the second hypothesis, in 145c7--e3 also leads to 138a3--7. The passage is wholly concerned with rejecting all the naïve types of presence and containment (the presence of the One in itself, in one, in several and in all of its parts), since in fact the One is certainly \textit{not contained} as a part of the Being part, but \textit{on the contrary} the One part \textit{contain}s as a part the Being part, ending with the statement that the One, in view of these naïve rejections, must be in the Other (`en allo', `en heteroi,145e2--3), thus, in effect, referring the explanation to the earlier passage 138e3--7. Thus, the passage 145c7--e3 offers very little and the essential explanation of the presence of the One in the Other-Being is contained in the passage 138a3--7.

\subsection{The presence of the One in the Being in 144c2--d4 is expressed equivalently by the circularity of contacts of parts of the One with parts of the Other 138a3--7}

\begin{enumerate}
    \item [{[1]}] — Ἐν ἄλλῳ μὲν ὂν\\
     If [the One] were in the Other 138a3
    \item [{[2]}] \textit{κύκλῳ} που ἂν περιέχοιτο ὑπ' ἐκείνου ἐν ᾧ ἐνείη
    [the One] would be \textit{cyclically} contained\\
    by that [the Other] in which [the One] would be, 138a4
    \item [{[3]}] καὶ πολλαχοῦ ἂν αὐτοῦ ἅπτοιτο\footnote{ ἅπτοιτο, ἅπτεσθαι +genitive. Hence
ἅπτοιτο αὐτοῦ [tou allou].} πολλοῖς· \\
and [the One] would contact [the Other] in many places (πολλαχοῦ) by many of its [of the One's] parts (πολλοῖς) 138a5
    \item [{[4]}] [τοῦ δὲ ἑνός τε καὶ ἀμεροῦς καὶ κύκλου μὴ μετέχοντος ἀδύνατον]\\
    πολλαχῇ \textit{κύκλῳ} ἅπτεσθαι.\\
    {}[— Ἀδύνατον]. 138a6--7\\
    {}[the One would] contact cyclically [the Other] in many places. (πολλαχῇ)
\end{enumerate}

The 138a3--7 passage is the only passage in the \emph{Parmenides}, where the meaning of the presence of the One in the Other is explained, and clearly this must be the meaning that
will explain
the presence of the One in the Being in 144c2--d4.
According to 138a3--7,
`the One is (present in) the Other'

\smallskip

[1] Ἐν ἄλλῳ μὲν ὂν 

[the One]

is in 
the Other (alloi) 138a3
means that  

\smallskip

 [2] κύκλῳ που ἂν περιέχοιτο ὑπ' ἐκείνου ἐν ᾧ ἐνείη,  

‘the One 
is circularly comprehended 
(‘kukloi… periechoito’, 138a4)
 by the Other.’\footnote{ Circularity will be discussed in Section 5.8, after
contacts-‘hapseis’ have been clarified.}
 [the One] 
would be circularly contained
 by that [the Other] 
in which [the One] would be 138a4
equivalently that  
‘the One, 
by means of many 
of its parts (‘pollois’, 138a5), 
is contacting (‘haptoito’, 138a5) 
circularly
the Other, 
in many of its places
(‘pollachou’,138a5)’

\smallskip

[3] καὶ πολλαχοῦ ἂν αὐτοῦ ἅπτοιτο πολλοῖς 

and [the One] 
would contact 
[the Other]
in many places (pollachou) 
by many parts (pollois) ·138a5

\smallskip

[4] πολλαχῇ κύκλῳ ἅπτεσθαι. 

[the One] 
would contact (haptesthai) 
cyclically (kukloi)
[the Other] 
in many places (pollachei) 138a6--7

Two new crucial terms enter the 138a3--7 description of the presence of the One in the Other, and we are led by Plato to the consideration of: 

contact (‘haptoito’) and circularity (‘kukloi’). 

\subsection{The meaning of `hapsis'
in the \emph{Parmenides} (`haptoito' in 138a3--7) is the philosophic analogue of the `ratio between two successive parts' in the sequence of remainders of the anthyphairesis of One to Being (`hapsis') in 148d--149d.}

In Section 5.4.1 the meaning of contacts/'hapseis', and the plus one rule between terms and contacts the statement 148d--149d, and in Section 5.4.2 the musical origin of the plus one rule, and that the One is \textit{contact (`haptoito')} the Other-Being part in many places, namely that the presence of the One in the Being is achieved by means of \textit{contacts (`hapseis').}

\subsubsection{According to passage 148d5--149d7}
\begin{itemize}
    \item a contact (`hapsis') is a relation between two consecutive terms next in the ordered sequence of the anthyphairetic remainders of the One to Being (148e4--7, 149a4--6), and
    \item the contacts satisfy the ``plus one rule": number of terms $=$ hapseis$+1$ (149a7--c3)
\end{itemize}
The detailed explanation is indeed given in the \emph{Parmenides} 148d5--149d7 passage on contacts (hapseis).
\begin{enumerate}
    \item [{[1]}] a contact (`hapsis') between two terms occurs if they are next to each other (`ephexes'), consecutive in their ordering 148e4--7 In [1] contact between two consecutive terms-parts is described,
    \begin{quote}
    — Τί δὲ τῇδε;\\
    ἆρ' οὐ πᾶν τὸ μέλλον ἅψεσθαί τινος\\
    ἐφεξῆς δεῖ κεῖσθαι ἐκείνῳ οὗ μέλλει ἅπτεσθαι\\
    ταύτην τὴν ἕδραν κατέχον ἣ ἂν μετ' ἐκείνην ᾖ [ἕδρα] ᾗ ἂν κέηται, \textit{ἅπτεται;} \\
    — Ἀνάγκη. \\
    `But how about this?\\
    Must not everything which is to \textit{contact (hapsesthai)} anything be \textit{next (ephexes)} to that which it is to \textit{contact (haptesthai}), and occupy that position which, being next to that of the other, \textit{contacts (haptetai)} it?" ``It must"' 148e4--7
    \end{quote}

    \item [{[2]}] generation of number two in terms of an `hapsis' 149a4--7
    
    In [2] this contact makes is related to the generation of the number two,
    
    \begin{quote}
    — Ὅτι, φαμέν,\\
    τὸ μέλλον ἅψεσθαι\\
    χωρὶς ὂν\\
    ἐφεξῆς δεῖ ἐκείνῳ εἶναι οὗ μέλλει ἅψεσθαι,\\
    τρίτον δὲ αὐτῶν ἐν μέσῳ μηδὲν εἶναι.\\
    — Ἀληθῆ.\\
    — Δύο ἄρα δεῖ τὸ ὀλίγιστον εἶναι, εἰ μέλλει ἅψις εἶναι.\\
    — Δεῖ. \\
    ``Because, as we agreed,\\
    that which is to contact (hapsesthai) anything must be\\
    outside of that which it is to contact (hapsesthai),\\
    and next (ephexes) it,\\
    and there must be no third between them."\\
    ``True."\\
    ``Then there must be two, at least (to oligiston), if there is to be contact (hapsis)"\\
    ``There must." 149a4--7
    \end{quote}

    \item [{[3]}] generation of number three in terms of two hapseis 149a8--b1
    
    In [3] two contacts between three consecutive terms-parts makes possible the generation of the number three,
    
    \begin{quote}
    — Ἐὰν δὲ τοῖν δυοῖν ὅροιν τρίτον προσγένηται ἑξῆς,\\
    αὐτὰ μὲν τρία ἔσται, αἱ δὲ ἅψεις δύο.\\
    — Ναί. \\
    `But if a third is added (prosgenetai) to the two terms (horoin)\\
   \textit{ in succession (hexes)},\\
   they will themselves be three, but the contacts (hapseis) two'. \\
   `Yes." 149a8--b1
    \end{quote}

    \item [{[4]}] the general rule: number of terms $=$ hapseis $+1$ 149b1--c3 
    
    In [4] the general plus one rule between contacts and terms-parts-existents is enunciated: the number of consecutive terms-existents is equal to the contacts between consecutive terms plus one.
   
    \begin{quote}
    — Καὶ οὕτω δὴ ἀεὶ ἑνὸς προσγιγνομένου μία καὶ ἅψις προσγίγνεται, καὶ συμβαίνει τὰς ἅψεις τοῦ πλήθους τῶν ἀριθμῶν μιᾷ ἐλάττους εἶναι. ᾧ γὰρ τὰ πρῶτα δύο ἐπλεονέκτησεν τῶν ἅψεων εἰς τὸ πλείω εἶναι τὸν ἀριθμὸν ἢ τὰς ἅψεις,\\
    τῷ ἴσῳ τούτῳ καὶ ὁ ἔπειτα ἀριθμὸς πᾶς πασῶν τῶν ἅψεων πλεονεκτεῖ·\\
    ἤδη γὰρ τὸ λοιπὸν ἅμα ἕν τε τῷ ἀριθμῷ προσγίγνεται καὶ μία ἅψις ταῖς ἅψεσιν.\\
    — ᾿Ορθῶς.\\
    — Ὅσα ἄρα ἐστὶν τὰ ὄντα τὸν ἀριθμόν,\\
    ἀεὶ μιᾷ αἱ ἅψεις ἐλάττους εἰσὶν αὐτῶν.\\
    — Ἀληθῆ. \\
    ``And thus\\
    whenever one is added (prosgignomenou),\\
    one contact (hapsis) also is added (prosgignetai), and\\
    the contacts is always one less than the multitude of number [of terms];\\
    for every succeeding number of terms exceeds the contacts\\
    just as much\\
    as the first two terms exceeded the contacts.\\
    For after the first, every additional term adds (prosgignetai) \\
    one to the number one to the contacts."\\
    
    ``Right."\\
    ``Then, whatever (`hosa') the number of beings (`onta'), the contacts are always one less." \\
    ``True." 149b1--c3
    \end{quote}
    
    Thus, Plato's interest in contacts is concentrated almost exclusively on the enunciation of \textit{the `plus one'} rule between number of terms and contacts.

    \item [{[5]}] The terms in the 148d--149d passage are the parts of the One
    
    \begin{itemize}
        \item The whole discussion in the passage is about number of terms; but we have already seen (in Section 4.1) that the number two is generated by parts One and Being as units, and the number three consisting of the parts One, Being, and another part of Being as units (Section 4.2). It is then clear that the terms (horoi'), which serve as the units of numbers in the present passage 148d--149d, are precisely the parts of the division of the One Being, which serve as the units of numbers throughout 143d--144c.
        \item It is remarkable that Plato, certainly on purpose, refers only to the number of terms, and never to the number of contacts. He phrases the plus one rule,
        \begin{quotation}
        NOT as number of terms$=$number of contacts $+1$,
   
        BUT always as number of terms$=$contacts$+1$.
        \end{quotation}
         This indicates that the units for the dialectical numbers are never the contacts, but solely the parts of the One.
        \item Furthermore, these terms are called beings (onta') in 149c2, the same name used throughout 143d1--144c8 for the units of the numbers in the One Being. These observations solidify the interpretation of the terms in the present passage as parts of the One Being.
    \end{itemize}
\end{enumerate}

\begin{quotation}
The order of the parts of the One/terms in 148d-149d is the order of their generation, namely the greater parts come first, the smaller after.
\end{quotation}

The main point of the passage 148d--149d on contacts-hapseis is the role that the contacts play for the generation of numbers. Now contacts make sense only on a set of terms, on which there is some ordering $\ll$, defined in such a way that every term if it has a successor, then it has \textit{an immediate successor} (ephexes', 148e5,8, 149a5), namely

for every $x$,
if there is a $z$ such that if $x\ll z$,
then there is a unique $x'$,
such that $x\ll x'$, and for every $t$, if $x'\ll t$, then either $x'=t$ or $x'\ll t$.

Since,

(a) as we have seen, the terms on which the plus one rule applies are the parts of the One Being,

(b) an ordering is needed on these terms to make sense of ouchings, and

(c) the only conceivable ordering on the parts of the One Being is the order of their generation

\noindent it appears that we have little choice but to conclude that the contacts are intended to apply to the set of all parts in the One Being
$$\text{One} > \text{Being} > \text{One}_1 > \text{Being}_1 > \ldots > \text{One}_k > \text{Being}_k > \text{One}_{k+1} > \text{Being}_{k+1} > \ldots$$
with the order $\ll$ defined as the order of their generation (and decreasing size), namely the order defined by setting $x\ll y$ if $y$ is a part of $x$ ($y<x$). So, denoting with $x*y$ the contact of two terms $x,y$, the sequence of the parts of the One on which the contact relation applies takes the form:
$$\text{One}*\text{Being}*\text{One}_1*\text{Being}_1*\ldots*\text{One}_k*\text{Being}_k*\text{One}_{k+1}*\text{Being}_{k+1}*\ldots$$
But the `plus one' rule on `contacts' still appears to be of little use, if we do not obtain a clear understanding of the concept of `contact'. So the next task is precisely this clarification.

\subsubsection{The plus one rule for contacts (hapseis) in 148d5--149d7: number of consecutive terms $=$ contacts $+1$ appears to have its origin in the plus one rule for musical intervals: number of consecutive terms $=$ musical intervals $+1$}

So exactly what is an ``hapsis"? It is crucial for an understanding of the One of the second hypothesis that the plus one rule appears before Plato solely in one place: music. Once we realize that the only use for Plato of the concept of `contact' is the enunciation of the `plus one rule' in 149b1--c3, it becomes possible to reveal the nature of `contacts-hapseis. Because the `plus one rule' is used elsewhere in Plato and also it appears in several passages on musical theory. Here is a list of some such prominent occurrences.

\begin{enumerate}
    \item \textit{Republic} 546b5--6  

\hspace*{1.5em}‘τρεῖς ἀποστάσεις, τέτταρας δὲ ὅρους λαβοῦσαι’ 

\hspace*{1.5em}when they have attained three intervals (apostasies) and four terms (horous)
\textit{Republic} 546b5--6

In the description of the Geometric Number in the \textit{Republic} 546b--c, Plato speaks
of four terms and three intervals:

4 terms (horoi)=3 intervals (apostasies)+1

\item Proclus, \textit{Commentary on Plato’s Republic} 2,36,21--25; 2,37,12--14

αὗται δ' οὖν αἱ αὐξήσεις μέχρι τεττάρων ὅρων προελθοῦσαι 

τρεῖς ἐχόντων ἀποστάσεις ἀλλήλων 

(πάντων γὰρ τεττάρων ὅρων συνεχῶν τρεῖς εἰσιν ἀποστάσεις)’;

These increases proceeded until there were four terms that have three intervals
between them – for in every case of four terms in succession there are three
intervals – that make all things ‘rational and expressible’ [with one another]
(546b7–c1): Proclus, \textit{Commentary on Plato's Republic}, 2, 36,2 1--25

[Proclus, \textit{Commentary on Plato’s Republic}, 2022~\cite{Proclus2022} p.252]
%volume II, Essays 7–15, 2022. 
%Translation by D. Baltzly, J. F. Finamore, G. Miles, Cambridge, p.252]

\hspace*{1.5em}‘τούτων δὴ τῶν τεττάρων ὄντων ἐφεξῆς ἐν τῷ ἐπιτρίτῳ λόγῳ ὅρων, 

\hspace*{1.5em}κζʹ λϚʹ μηʹ ξδʹ [27 36 48 64], καὶ τρεῖς ἀποστάσεις ἐχόντων’

\hspace*{1.5em}Now among these four terms (tessaron horon) in continuous (ephexes) [proportion]
in the epitritos ratio (en toi epitritoi logoi) having three intervals (treis
apostasies) – 27, 36, 48 and 64– 2,37,12--14

[translation, 2022~\cite{Proclus2022}, p.253]

\item \textit{Republic} 616d6--617d5 

Harmony of the Sirens 

8 terms (horoi)=7 intervals (diastemata)+1

\item Proclus, \textit{Commentary on Plato’s Timaeus}, 1820~\cite{Proclus1820} 2,237,3--15

\hspace*{1.5em}`τὴν διὰ πασῶν, 

\hspace*{1.5em}ἐν ὅροις μὲν ὀκτὼ θεωρουμένην, ἑπτὰ δὲ διαστήμασιν'

\hspace*{1.5em}the octave, considered in eight terms (en horois okto), and seven intervals
(hepta diastematon) 2,237, 5--6:

\item \textit{Timaeus} 34b4--35b6, 

Division of the soul by means of musical intervals.

\item Proclus,  \textit{Commentary on Plato’s Timaeus}, 1820~\cite{Proclus1820} 2,187,13--15; 2,188,6--9

Proclus has worked out Plato’s system in exhaustive detail, and he has found
that 

there are 34 terms (horoi) and 33 musical intervals-ratios.

\hspace*{1.5em}‘καὶ πᾶν τοῦτο τὸ διάγραμμα λείμματα μὲν ἔχει ἐννέα, ἐπογδόους δὲ
εἰκοσιτέσσαρας· τὰ γὰρ διαστήματα ἑνὶ λείπεται τῶν ὅρων’ 

\hspace*{1.5em}and this whole diagram consists of nine lemmata (dieseis) , and twentyfour tones
(epogdoois); because the intervals is less by one of the terms 

2,187,13--15:

\hspace*{1.5em}‘ταῦτα δ' ἦν τὰ ῥηθέντα, καὶ ἐκ τούτων εὐτάκτως ληφθέντων ἀνεφάνησαν ὅροι
τριακοντατέσσαρες μόνοι πᾶν τὸ διάγραμμα περιέχοντες.’ 

\hspace*{1.5em}From all that has been said and taken into account above, there appeared thirty
four (triakontatessares) terms (horoi) making up the whole diagram(diagram)
2,188,6--9

Plato’s division with 34 terms and 33 intervals, according to Proclus is the
following:

\begin{center}
\tablefirsthead{}
\tablehead{}
\tabletail{}
\tablelasttail{}
\begin{supertabular}{|c|c|c|c|c|c|c|c|c|c|c|c|c|c|c|}
\hline
\centering 384 &
\centering t &
\centering 432 &
\centering t &
\centering 486 &
\centering d &
\centering 512 &
\centering t &
\centering 576 &
\centering t &
\centering 648 &
\centering t &
\centering 729 &
\centering d &
\centering\arraybslash 768\\\hline
\centering [768] &
\centering t &
\centering 864 &
\centering t &
\centering 972 &
\centering d &
\centering 1024 &
\centering t &
\centering 1152 &
\centering t &
\centering 1296 &
\centering t &
\centering 1458 &
\centering d &
\centering\arraybslash 1536\\\hline
\centering [1536] &
\centering t &
\centering 1728 &
\centering t &
\centering 1944 &
\centering d &
{\centering [18]\par}

\centering 2048 &
\centering t &
{\centering [19]\par}

\centering 2304 &
\centering t &
\centering 2592 &
\centering t &
\centering 2916 &
\centering d &
{\centering 3072\par}

\\\hline
\centering [3072] &
\centering t &
\centering 3456 &
\centering t &
\centering 3488 &
\centering t &
\centering 4374 &
\centering d &
\centering 4608 &
 &
 &
 &
 &
 &
\\\hline
\centering [4608] &
\centering t &
\centering 5184 &
\centering t &
\centering 5832 &
\centering d &
{\centering [29]\par}

\centering 6144 &
\centering t &
{\centering [30]\par}

\centering 6912 &
\centering t &
\centering 7776 &
\centering t &
\centering 8748 &
\centering d &
\centering\arraybslash 9216\\\hline
\centering [9216] &
\centering t &
\centering 10368 &
 &
 &
 &
 &
 &
 &
 &
 &
 &
 &
 &
\\\hline
\end{supertabular}
\end{center}
t=tone, d=diesis

34 terms (horoi) =33 intervals (diastemata)+1

\item Timaeus Locrus, \textit{On Nature} (Peri phusios) 1972~\cite{TimaeusLocrus1972} 209,6--7

δεῖ δ' εἶμέν πως πάντας σὺν τοῖς συμπληρώμασι καὶ τοῖς ἐπογδόοις 

 ὅρους ἓξ καὶ τριάκοντα’ On Nature 

Mit den Auffüllungen und den Ganztönen (epogdoois) müssen es im Ganzen 36
Glieder (horous hex kai triakonta) sein 

Timaeus Locrus, \textit{Peri phusios} 209,6--7

[[translation Timaeus Locrus, \textit{De Natura Mundi et Animae}, translated by Walter
Marg, E.J. Brill, 1972~\cite{TimaeusLocrus1972}]]

Timaeus Locrus introduces an interesting modification of Plato’s Timaeus
division with the introduction of the apotome interval in two places:

the first, introducing between the 18th term 2048 and the 19th term 2304,
forming the interval of a tone, the term 2187, thus splitting the interval of
the tone in two intervals, an apotome (2187 to 2048) and a diesis (2304 to
2187) (cf. Proclus, \textit{Commentary to Plato’s Timaeus}, 1820~\cite{Proclus1820} 2,189,16--17), and, the
second, introducing between the 29th term 6144 and the 30th term 6912, forming
the interval of a tone, the term 6561, thus splitting the interval of the tone
in two intervals, an apotome (6561 to 6144) and a diesis (6912 to 6561) (cf.
Proclus, \textit{Commentary to Plato’s Timaeus} \cite{Proclus1820} 2,189,25--190,1), 

resulting in a division with 36 terms (horous) and 35 intervals.

Tmaeus Locrus’s division with 36 terms and 35 intervals, according to Proclus is
the following:

\begin{center}
\tablefirsthead{}
\tablehead{}
\tabletail{}
\tablelasttail{}
\begin{supertabular}{|c|c|c|c|c|c|c|c|c|c|c|c|c|c|c|c|c|}
\hline
\centering 384 &
\centering t &
\centering 432 &
\centering t &
\centering 486 &
\centering d &
\centering 512 &
\multicolumn{3}{|c|}{\centering t} &
\centering 576 &
\centering t &
\centering 648 &
\centering t &
\centering 729 &
\centering d &
\centering\arraybslash 768\\\hline
\centering [768] &
\centering t &
\centering 864 &
\centering t &
\centering 972 &
\centering d &
\centering 1024 &
\multicolumn{3}{|c|}{\centering t} &
\centering 1152 &
\centering t &
\centering 1296 &
\centering t &
\centering 1458 &
\centering d &
\centering\arraybslash 1536\\\hline
\centering [1536] &
\centering t &
\centering 1728 &
\centering t &
\centering 1944 &
\centering d &
2048 &
\multicolumn{1}{|c|}{\centering a} &
\multicolumn{1}{|c|}{2187} &
\centering d &
\centering 2304 &
\centering t &
\centering 2592 &
\centering t &
\centering 2916 &
\centering d &
\centering\arraybslash 3072\\\hline
\centering [3072] &
\centering t &
\centering 3456 &
\centering t &
\centering 3488 &
\centering t &
\centering 4374 &
\multicolumn{3}{|c|}{\centering d} &
\centering 4608 &
 &
 &
 &
 &
 &
\\\hline
\centering [4608] &
\centering t &
\centering 5184 &
\centering t &
\centering 5832 &
\centering d &
6144 &
\multicolumn{1}{|c|}{\centering a} &
\multicolumn{1}{|c|}{6561} &
\centering d &
6912 &
\centering t &
\centering 7776 &
\centering t &
\centering 8748 &
\centering d &
\centering\arraybslash 9216\\\hline
\centering [9216] &
\centering t &
\centering 10368 &
 &
 &
 &
 &
\multicolumn{3}{|c|}{} &
 &
 &
 &
 &
 &
 &
\\\hline
\end{supertabular}
\end{center}
t=tone, d=diesis, a=apotome

36 terms (horoi)=35 intervalks (diastemata)+1

Proclus in his \textit{Commentary to Plato’s Timaeus} describes the transition from the
Platonic division to the Timaeus of Locris division in 

\hspace*{1.5em}ἀνάγκη δήπου πρὸς τὸν δισχίλια τεσσαράκοντα καὶ ὀκτὼ τὴν ἀποτομὴν ποιεῖν 

\hspace*{1.5em}τὸν δισχίλια ἑκατὸν ὀγδοήκοντα ἑπτά 2,189,16--17, and 

\hspace*{1.5em}it is necessary that to the number 2048 an apotome is formed by the number 2187

\hspace*{1.5em}ἔχομεν ἀναγκαίως τὴν ἀποτομὴν ἐν τῷ λόγῳ τοῦ ἑξακισχίλια πεντακόσια ἑξήκοντα καὶ
ἑνὸς πρὸς τὸν ἑξακισχίλια ἑκατὸν τεσσαράκοντα τέσσαρα’  2,189,25--190,1. 

 \hspace*{1.5em}It is necessary that to the number 6144 an apotome is formed by the number 6561

\hspace*{1.5em} [translation by the author 1972 \cite{TimaeusLocrus1972}]

\item \textit{Timaeus} 36d2--4

\hspace*{1.5em}τὴν δ' ἐντὸς σχίσας ἑξαχῇ ἑπτὰ κύκλους ἀνίσους 

\hspace*{1.5em}κατὰ τὴν τοῦ διπλασίου καὶ τριπλασίου διάστασιν ἑκάστην, οὐσῶν ἑκατέρων τριῶν, 

\hspace*{1.5em}He split the inner Revolution in six places (hexachei) into seven (hepta)
unequal circles, according to each of the intervals of the double and triple
intervals,~Timaeus 36,2--4:

7 circles=6 intervals+1

\item Iamblichus, \textit{Theologoumena Arithmeticae}, translation, 1988~\cite{Iamblichus1988} 48,10--14; 50,5--6

\hspace*{1.5em}‘ἑπτὰ γὰρ κινημάτων ἀστερικῶν ὑπαρχόντων 

\hspace*{1.5em}παρὲξ τοῦ τῶν ἀπλανῶν ὀγδόου μέν, οὐχ ἁπλοῦ δέ, 

\hspace*{1.5em}καὶ φθόγγους ἀποτελούντων ἰσαρίθμους διὰ τῆς ῥοιζήσεως, 

\hspace*{1.5em}ἀνάγκη τὰ διαστήματα αὐτῶν καὶ οἷον μεσότητας ἓξ ὑπάρχειν.’; 

\hspace*{1.5em}for since there are seven celestial movements

\hspace*{1.5em}(apart from the movement of the fixed stars, which is eighth, but complex),

\hspace*{1.5em}and since by their hurtling they produce the same number of notes (phthongous),

\hspace*{1.5em}then their intervals (diastemata) and, as it were, means (mesotetas) are
necessarily six (hex). 48,10--14

[[Iamblichus, \textit{The Theology of Arithmetic}, translated by Robin
Waterfield, Phanes Press, 1988~\cite{Iamblichus1988} p.80]]

\hspace*{1.5em}ὅτι ἑπτὰ τῶν σφαιρῶν οὐσῶν κατὰ τὴν ἑξάδα τὰ διαστήματά ἐστι· 

\hspace*{1.5em}μονάδι γὰρ ἀεὶ ἐλάττονα. 

\hspace*{1.5em}Since there are seven (hepta) celestial spheres, the intervals (diastemata) 
fall under the hexad (hexada) : for they are always less by a monad (monadi gar
aei elattona). 50,5--6

[[translation, cite{Iamblichus1988} p.82]]

 7 notes (phthongoi)=6 intervals (diastemata)+1

\item Euclid(?), \textit{Sectio Canonis}, Proposition 9

\hspace*{1.5em}Τὰ ἓξ ἐπόγδοα διαστήματα μείζονά ἐστι διαστήματος ἑνὸς διπλασίου. 

\hspace*{1.5em} ἔστω γὰρ εἷς ἀριθμὸς ὁ Α. 

\hspace*{1.5em}καὶ τοῦ μὲν Α ἐπόγδοος ἔστω ὁ Β, τοῦ δὲ Β ἐπόγδοος ὁ Γ, 

\hspace*{1.5em}τοῦ δὲ Γ ἐπόγδοος ὁ Δ, τοῦ δὲ Δ ἐπόγδοος ὁ Ε, τοῦ Ε ἐπόγδοος ὁ Ζ, τοῦ Ζ ἐπόγδοος
\hspace*{1.5em}ὁ Η· λέγω, ὅτι ὁ Η τοῦ Α μείζων ἐστὶν ἢ διπλάσιος. 

\hspace*{1.5em}ἐπεὶ ἐμάθομεν εὑρεῖν ἑπτὰ ἀριθμοὺς ἐπογδόους ἀλλήλων, εὑρήσθωσαν οἱ Α, Β, Γ, Δ,
Ε, Ζ, Η, καὶ γίνεται ὁ μὲν Α κϚ μύρια βρμδ, [262 144] ὁ δὲ Β κθ μύρια δϠιβ,
[294 912] 

\hspace*{1.5em}ὁ δὲ Γ λγ μύρια αψοϚ, [331 776] ὁ δὲ Δ λζ μύρια γσμη, [373 248] ὁ δὲ Ε μα μύρια
θϠδ, [419 904] 

\hspace*{1.5em}ὁ δὲ Ζ μζ μύρια βτϞβ, [472 392] ὁ δὲ Η νγ μύρια αυμα, [531 441] 

\hspace*{1.5em}καί ἐστιν ὁ Η [531 441] τοῦ Α μείζων ἢ διπλάσιος. [524 288].

\hspace*{1.5em}Six epogdoic intervals are greater than one duple interval 

\hspace*{1.5em}Let A be one number. 

\hspace*{1.5em}Let B be the epogdoic of A, let C be the epogdoic of B, let D be the epogdoic of
C, let E be the epogdoic of D, let F be the epogdoic of E, and let G be the
epogdoic of F. 

\hspace*{1.5em}I say that G is more than double A. 

\hspace*{1.5em}Since we have learned how to find seven numbers that are epogdoics of one
another, let the numbers A, B, C, D, E, F, G have been found. 

\begin{center}
\tablefirsthead{}
\tablehead{}
\tabletail{}
\tablelasttail{}
\begin{supertabular}{|c|c|c|c|c|c|c|}
\hline
\centering A &
\centering B &
\centering C &
\centering D &
\centering E &
\centering F &
\centering\arraybslash G\\\hline
\centering 262,144 &
\centering 294,912 &
\centering 331,776 &
\centering 373,248 &
\centering 419,904 &
\centering 472,392 &
\centering\arraybslash 531,441\\\hline
\end{supertabular}
\end{center}

\hspace*{1.5em}and G is more than double 524 288 of A.

[translation by Andrew Barker, 1989. Greek Musical Writings, vol.2, Cambridge
University Press, 1989 \cite{Negrepontis2025b} p.199]

7 numbers (arithmoi) = 6 intervals (diastemata)+1. 

 (cf. also Porphyrius, \textit{Harmonics} 101,25--102,7) 

[contains the proof that six epogdoa composed is an interval greater than dia
pason. This is done by finding seven numbers A, B, C, D, E, F, G, such that

A*B=B*C=C*D=D*E=E*F=F*G=9:8, and noting that A*G>2:1]

\item Nicomachus, \textit{Harmonicum enchiridion}, translation 1994 {\huge???} 12,1,32--40 

[[4 arithmoi=3 diastemata+1, etc]]

\hspace*{1.5em}‘τρία διαστήματα ἐν τέσσαρσιν ἀριθμοῖς ὅ ἐστι φθόγγοις.’  

\hspace*{1.5em}three intervals (diastemata) in four (tessarsin) numbers (arithmois), 

namely notes (phthongois) Nicomachus, \cite{Barker1989} p.267 % 12, 1,37--38

\item Aristides Quintilianus, \textit{de Musica}, translation 1989 \cite{Barker1989}, Chapter 12, 3,1, 32--82 

[various instances of the ‘plus one rule’]

In at least two of these occurrences the general rule is enunciated explicitly \cite{Barker1989} p.392--535:

Iamblichus, \textit{Theologoumena Arithmeticae}, translation, 1988~\cite{Iamblichus1988}: 

\hspace*{1.5em}μονάδι γὰρ ἀεὶ ἐλάττονα

\hspace*{1.5em}[the intervals] always less (alattona) [than the terms] by a unit (monadi)
50,5--6

Proclus, \textit{Commentary to Plato’s Timaeus} 2,187,13--15:

τὰ γὰρ διαστήματα ἑνὶ λείπεται τῶν ὅρων

The intervals (diastemata) are less (leipetai) by one (heni) than the terms
(horon)

In all the cases mentioned, the ‘plus one rule’ in the \textit{Parmenides} 148d--149d

number of consecutive terms=contacts of consecutive terms +1

applies to musical intervals (diastemata’, ‘apostasies’) in place of
contacts(‘hapseis’):

number of consecutive terms=musical intervals of consecutive terms +1.
\end{enumerate}

\subsection{A contact (hapsis) in the One of the Second Hypothesis in the \emph{Parmenides} is $\alpha$ ratio of two consecutive terms in the sequence of the anthyphairetic remainders of the dyad One and Being.}
But a musical interval between two terms is just the ratio between the terms involved, and thus the musical plus one described in Section 5.3.2m is the rule
\begin{quotation}
number of consecutive terms=ratios of consecutive terms +1;
\end{quotation}
and thus, comparing Sections 5.3.1 and 5.3.2, we conclude that a contact/'hapsis' between two consecutive parts in the sequence of the parts of the One Being can only be the ratio of these two parts. This provides us with an essential clarification of the terms `hapseis'. The contacts* (hapseis) have been identified with the `logoi'$=$ratios of successive (in their order of generation) of parts of the One Being:
\begin{align*}
  \text{One} &* \text{Being},\\
  \text{Being}&* \text{One}_1,\\ \text{One}_1 &* \text{Being}_1,\\ &\vdots\\ \text{One}_k &* \text{Being}_k,\\  \text{Being}_k &* \text{One}_{k+1},\\ \text{One}_{k+1} &* \text{Being}_{k+1},\\  &\vdots
\end{align*}
We then should have then no hesitation to conclude that the rule enunciated in detail in 149b1--c3 between the number of consecutive terms and the connective hapseis is really the rule
$$\text{number of consecutive terms}=\text{ratios of consecutive terms} +1.$$

\textit{Note}. Since the absolutely partless One of the first hypothesis fails to satisfy one-by-one all the hypotheses that the One of the second hypothesis in the \emph{Parmenides} satisfies, and since there is no other mention of `haptesthai', `hapseis' in the first hypothesis beyond those in 138a3--7, it is clear that both passages, 138a3--7 and 148d5--149d7, deal with the same sort of contact. Hence, information contained in the latter about contacts is relevant to the former. It is clear that the terms `haptoito' and `haptesthai' in passage 138a in the first hypothesis correspond to the `hapseis' passage 148d5--149d7 in the second hypothesis. Therefore, the understanding of the nature of `hapseis that will result from a study of the passage 148d--149d may be gainfully employed for reaching an understanding of the role of `hapseis' in the 138a passage.

\subsection{The circularity (138a5--7) of contacts/hapseis (148d--149d), interpreted as ratios of the terms of the sequence of successive remainders, leads to (a) the Logos Criterion ($\text{One}/\text{Being}=\text{One}_k/\text{Being}_k$ for some $k$) for the anthyphairetic periodicity of the dyad One \& Being of the Second Hypothesis in the \emph{Parmenides}, and (b) the interpretation of the statements ``the One is (present) in the other/Being" (144c2--d4, 138e3--7, 145b6--e6), and ``the One is in itself" (145b6--e6)}

We now have all the pieces necessary for a satisfactory interpretation of the crucial passage 138a3--7.

\begin{enumerate}
    \item [1)] Ἐν \textit{ἄλλῳ} μὲν ὂν\\
     The One is in the Other/the Being
\end{enumerate}
The One part is greater than the Being part. So if the One will be somehow end being contained in the Being, this must take place in some non-obvious way.

\begin{enumerate}
    \item [2)] πολλαχοῦ ἂν αὐτοῦ ἅπτοιτο πολλοῖς·\\
    there are contacts/hapseis (haptoito)\\
    of many parts of the Being (pollachou)\\
    by many parts (pollois) of the One (autou)
\end{enumerate}
Hapseis presuppose a sequence of parts ordered as the natural numbers, so that there will be parts immediately next to each other. As explained in Section 5.3, the sequence referred to here is the sequence
$$\text{One} > \text{Being} > \text{One}_1 > \text{Being}_1 > \ldots > \text{One}_k > \text{Being}_k > \text{One}_{k+1} > \text{Being}_{k+1} > \ldots;$$
and, as explained in Section 5.2,

\begin{enumerate}
    \item [3)] πολλαχῇ κύκλῳ ἅπτεσθαι.\\
    these hapseis form a circle
\end{enumerate}
Here in [3] for the first time, circularity of contacts appears! Thus in the sequence with terms the parts of the One Being
$$\text{One} * \text{Being} * \text{One}_1 * \text{Being}_1 * \ldots * \text{One}_k * \text{Being}_k * \text{One}_{k+1} * \text{Being}_{k+1} * \ldots$$
the sequence of contacts $*$ of successive terms form \textit{a circle}. Thus, \textit{there is an index $k$} such that
$$\text{One} * \text{Being} = \text{One}_k * \text{Being}_k.$$
But now we must appeal to two of our earlier conclusions: The first earlier conclusion (as explained in Section 5.3.3) is that the contacts are interpreted as \textit{ratios} of the sequence of successive parts of the One Being. Thus the ratios of successive parts of the sequence of parts of the One Being
$$ \text{One} > \text{Being} > \text{One}_1 > \text{Being}_1 > \ldots > \text{One}_k > \text{Being}_k > \text{One}_{k+1} > \text{Being}_\textit{{k+1}} > \ldots $$

form a circle. Thus

$$\text{One}/\text{Being}=\text{Onek}/\text{Beingk}.$$
%$$\frac{\text{One}}{\text{Being}} = \frac{\text{One}_k}{\text{Being}_k}.$$
The second earlier conclusion (established in Section 3) is that the sequence of parts of the One Being
$$\text{One} > \text{Being} > \text{One}_1 > \text{Being}_1 > \ldots > \text{One}_k > \text{Being}_k > \text{One}_{k+1} > \text{Being}_{k+1} > \ldots$$
is the sequence of the `remainders' of the anthyphairesis of One to Being. It is notable that the anthyphairetic nature of the Division of the dyad One Being is being used here for the first time. The equality of ratios just established above,
$$\frac{\text{One}}{\text{Being}} = \frac{\text{One}_k}{\text{Being}_k} \quad \text{for some } k,$$
is now seen to be precisely
\begin{enumerate}
    \item [(a)] the dyad One \& Being satisfies the \textbf{Logos Criterion} for the periodicity of the anthyphairesis of the One to the Being.
\end{enumerate}
Part (b) is an immediate consequence of the Logos Criterion: By (a), $\text{One}_k$ is cyclically identified with One; and $\text{One}_k$ is contained in the usual manner, as part, into Being. Hence
\begin{enumerate}
    \item [(b)] the One is contained in itself and in the Being/Other by means of a circle.
\end{enumerate}
This explains the passages 144c2--d4, 138e3--7, and 145b6--e6.

\subsection{The One of the second Hypothesis is described in 155d3--6 as an intelligible Being, whose knowledge is given by Name and Logos 155d8--e1, confirming that the structure of the paradigmatical intelligible Being One is the fact that the One with its part Being forms a dyad satisfying the philosophical analogue of a (full) periodic anthyphairesis}

The description of the intelligible Being in the \emph{Meno} as True Opinion \& Logos, and in the \emph{Sophist} as Name \& Logos has been interpreted as the philosophical analogue of Theaetetus' Logos Criterion for anthyphairetic periodicity in Negrepontis 2012~\cite{Negrepontis2012}, 2018~\cite{Negrepontis2018}, 2024~\cite{Negrepontis2024a}, Negrepontis, Farmaki, Brokou 2024~\cite{Negrepontis2024c}.

The One of the second hypothesis is an intelligible Being:
\begin{quote}
—Ἦν ἄρα τὸ ἓν καὶ ἔστι καὶ ἔσται καὶ ἐγίγνετο καὶ γίγνεται καὶ γενήσεται. 

— Τί μήν; 

— Καὶ εἴη ἄν τι ἐκείνῳ καὶ ἐκείνου, καὶ ἦν καὶ ἔστιν καὶ ἔσται. 

— Πάνυ γε. 155d3--6

“Then the one was and is and will be and was becoming and is becoming and will
become.” 

“Certainly.” 

“And there would be and was and is and will be something which is in relation to
it and belongs to it?” 

“Certainly.”~
\end{quote}
In the Second Hypothesis of the \emph{Parmenides} the use of the Theaetetean term Logos is mostly avoided, since it would be anachronistic to put in the mouth of Parmenides or Zeno the philosophic analogue of the Logos Criterion for the periodicity of anthyphairesis of magnitudes; in fact, the first theory of ratios of magnitudes is due to Theaetetus, as hinted by Aristotle in \emph{Topics} 158b, reconstructed in Negrepontis and Protopapas, 2025~\cite{Negrepontis2025b}. However, Logos does appear in relation to the Knowledge of the One of the Second Hypothesis in the \emph{Parmenides} 155d8--e1
\begin{quote}
Καὶ ὄνομα δὴ καὶ λόγος ἔστιν αὐτῷ, καὶ ὀνομάζεται καὶ λέγεται 155d8--e1, \\
``And it has a Name and Logos, is named and in possession of Logos,"
\end{quote}
and Logos is prominent in the definition of dialectical number Two, as explained in section 6.7. below, both constituting a confirmation of our common interpretation of both \emph{Theaetetus, Meno, Sophist}, with \emph{Parmenides} in terms of periodic anthyphairesis.

The structure of the paradigmatical intelligible Being the One of the Second Hypothesis in the \emph{Parmenides} in terms of the philosophical analogue of periodic anthyphairesis
\begin{quote}
the One in the Second Hypothesis,

a paradigmatical intelligible Being (155d--e),  

forms, with its other (heteron) part Being. the dyad One \& Being 

(a) satisfying the philosophical analogue of infinite anthyphairesis
(142b--143b), and

(b) satisfying in fact the philosophical analogue of periodic anthyphairesis
(144c2--d2, 138a3--7, 148d--149d). 
\end{quote}

\subsection{The dialectical proof of the infinity of the anthyphairesis of the intelligible Being One 144e3--7 answers the question posed in 143a4--9:}

\begin{quote}
``The One itself (τὸ ἓν αὐτὸ), then, split up (kekermatismenon) by the Being part, is many and \textbf{infinite in multitude} (apeira to plethos)."\\
``Clearly."\\
``Then not only (monon) the One \& Being (to hen on) is many,\\
but the One itself (auto to hen) divided (dianenememenon) by the Being part,\\
must be many."\\
``Certainly." 144e3--7
\end{quote}

In passage 144e3--7 Plato returns and answers the question posed in 143a, cf. Section 5.6. It was initially proved, in 142b--143a (Section 3), that

\begin{quotation}
the dyad One \& Being has infinite anthyphairesis.
\end{quotation} 
This was a proof, without any recourse to contacts/logoi, and was thus an \textbf{eristic} proof of infinity; by contrast, we have now show that the One and the Being, are in periodic anthyphairesis by the Logos Criterion formed by the cyclical containmment of the One in the Being (Section 5.8), and the question was then if there could be a proof that the One of the Second Hypothesis is a paradigmatical intelligible Being,
\begin{quotation}
the Platonic Idea `the One itself' possesses an infinite multitude of parts.
\end{quotation}
 The dialectical proof that the One itself, the paradigmatical Platonic Idea, does have an infinite multitude of parts is the following: the anthyphairesis of the two initial parts in the One itself is periodic, and periodicity of anthyphairesis implies infinity. This is then a \textbf{`dialectical'}, intelligible proof of the infinite nature of the One itself, quite different from the `eristic' proof, given in 142b--143a (Section 3).

\textit{Note}. Since the difference of the two proofs of infinity relies on the fact that the division of the One to the Being is anthyphairetic, and since scholars on Plato have completely ignored this nature of the dyad (One Being), it is not surprising that they have failed to understand the philosophically fundamental difference of the two proofs of infinity. Plato does not refer to the original proof of anthyphairetic infinity as eristic nor to the second as dialectical, but his intent is clear and revealed in the \emph{Philebus} 16e4--17a5, where he distinguishes between eristic potential infinite, not involving dialectical number, and dialectical already completed infinite, a clear rebuke to Archytas' proofs of incommensurabilities.

\subsection{Does the proof of the plus one rule in the \emph{Parmenides} 149a7--c3 employ Mathematical Induction?}

There has been a heated discussion on whether Greek Mathematics possessed mathematical induction. (Freudenthal, 1953~\cite{Freudenthal1953}, Mueller, 1981~\cite{Mueller1981}, Fowler, 1987~\cite{Fowler1987}, Unguru 1991~\cite{Unguru1991}, Fowler, 1994~\cite{Fowler1994}, Unguru 1994~\cite{Unguru1994}).

The proof of the Proposition number of terms $=$ hapseis$+1$ given by Plato in the \emph{Parmenides} 148d--149d and explained in Section 5.5, has been proposed as a perfect specimen of ancient proof by mathematical induction. Wedberg 1955~\cite{Wedberg1955}, and also Acerbi, 2000~\cite{Acerbi2000}, have analyzed Plato's \emph{Parmenides} 149a7--c3 about `contacts' (`hapseis'), and have concluded that it is a quite simple but logically impeccable proof by mathematical induction. However, our analysis above has shown that Plato means this to be an argument involving only finite induction, since there are only finitely many `hapseis'.

On the wider question of Mathematical Induction in Greek Mathematics,
let us note that Zeuthen~\cite{Zeuthen1886} p.27--28, Heath~\cite{Heath1926} vol.1,
p.398--401, Freudental, van der Waerden argue in favor of an
inductive proof, Fowler, Unguru, Acerbi (2000~\cite{Acerbi2000} p.72
`notably some among those connected with the property of side and
diagonal numbers, the proof itself is altogether lacking') against the
presence of a valid inductive proof, but the true arguments in the
ancient sources in favor of a valid inductive proof appear to have
been missed; in Negrepontis, Farmaki, Brokou,
2024~\cite{Negrepontis2024c} and with further arguments in
Negrepontis, Farmaki, 2025~\cite{Negrepontis2025a}, we argue that
Proclus in his \emph{Commentary to Plato's Republic} 2,24,16--25,13
and 2,27,1--29,4 2022~\cite{Proclus2022}, together with Theon of
Smyrna 43,5--43,15, 44,18--45,8, 1878~\cite{TheonSmyrnaeus1878},
1979~\cite{TheonSmyrna1979}, and Iamblichus, \emph{Comments to
  Nicomachus Arithmetic} 91,21--92,3,92,23--93,6,
2014~\cite{Jamblique2014} describe a valid proof by mathematical
induction of the basic Pell type property of the side and diameters
numbers.

\section{The definition of dialectical number in an intelligible Being in terms of the Logos Criterion (138a3--7) and the plus one rule (148d--149d), the finiteness of the dialectical number of the parts of the Being, the equality of the dialectical number of the parts of the One and of the dialectical number of the parts of the Being, and the rejection of the eristic numbers (144d4--e3)}

\subsection{The definition of the dialectical numbers in an intelligible Being, in terms of the anthyphairetic periodicity and the plus one rule}

The plus one rule between terms and contacts ratios applied to the sequence of parts/anthyphairetic remainders of an intelligible Being generates the dialectical numbers.
\begin{quotation}
Definition of the dialectic numbers in the One of the Second Hypothesis (\emph{Parmenides} 148d--149d). The dialectic number of a finite or infinite sequence $S$ of successive terms in the sequence of remainders

$$\text{One} > \text{Being} > \text{One}_1 > \text{Being}_1 > \ldots > \text{One}_k > \text{Being}_k > \text{One}_{k+1} > \text{Being}_{k+1} > \ldots $$

is equal to (the number of all different contacts (hapseis), namely, ratios formed from successive terms in $S$) $+ 1$.
\end{quotation}

Since, by anthyphairetic periodicity, a sequence of successive remainders is periodic, the multitude of all different ratios formed from the whole infinite sequence of parts of the One is finite, and thus the whole dialectic number of an intelligible Being is finite, never exceeding the length of the anthyphairetic period plus one.

\subsection{\emph{Parmenides} 144d4--5 on the equality and finiteness of the dialectical number of parts of the One with the dialectical number of parts of the Being}

\subsubsection{The passage 144d4--5 is the following}
\begin{quote}
— Καὶ μὴν τό γε μεριστὸν πολλὴ ἀνάγκη εἶναι τοσαῦτα ὅσαπερ μέρη.\\
— Ἀνάγκη. 144d4--5
\end{quote}

\subsubsection{A usual incorrect rendering of passage 144d4--5}

In the classical translation by H. N. Fowler, 1926~\cite{Fowler1926}, the passage is rendered as
\begin{quote}
`And that which is divided into parts (`to meriston') must certainly be as numerous as its parts (`mere')',
\end{quote}
as if `the parts (mere) were the parts of the `meriston', turning the statement into a trivial tautology; But according to the immediately previous 144d2--4 to 144d4--5 sentence the One is by necessity `memerismenon' (divided), because the One could not possibly be simultaneously present whole in all the `mere'(parts) of the Being part. Hence the sentence 144d4--5 cannot be rendered as Fowler does.

\subsubsection{Our rendering of the passage 144d4--5}

It is clear that the `meriston' of 144d4--5 is the `memerismenon' One of 144d2--4, while the `mere' of 144d4--5 are the `mere' of the Being of 144d2--4. Thus the correct rendering of 14d4--5 is:
\begin{quote}
`The One is divided into as many (`tosa') parts as (`hosaper') the parts of the Being', 14d4--5
\end{quote}

\subsection{Linguistic Note: The expression `tosa\ldots hosa' always indicates the equality of two finite numbers}

We show by examples that in general the expression `tosa\ldots hosa' indicates equality of two finite numbers.

\begin{enumerate}
   \item [{[1]}] Πότερα οὖν ὅσα ἐμοὶ δοκεῖ δεῖν ἀποκρίνεσθαι, 

τοσαῦτά σοι ἀποκρίνωμαι, 

ἢ ὅσα σοί; \textit{Protagoras}, translation \cite{Plato} 334e2--3

Then are my answers so many, 

as many I think they should be, 

or as many as you think they should be?

\item [{[2]}] ΣΩ. Ἦ καὶ ὁπόσα ἂν φῇ τις ἑκάστῳ ὀνόματα εἶναι,  

τοσαῦτα ἔσται καὶ τότε ὁπόταν φῇ; \textit{Cratylus}, translation \cite{Plato} 385d5--6

And as many names anyone says a thing has, 

so many will really be at the time when he says it?

\item [{[3]}] Οἶσθ' οὖν, ἦν δ' ἐγώ, ὅτι καὶ ἀνθρώπων εἴδη 

τοσαῦτα ἀνάγκη τρόπων εἶναι, 

ὅσαπερ καὶ πολιτειῶν; 

ἢ οἴει ἐκ δρυός ποθεν ἢ ἐκ πέτρας τὰς πολιτείας γίγνεσθαι, 

ἀλλ' οὐχὶ ἐκ τῶν ἠθῶν τῶν ἐν ταῖς πόλεσιν, 

ἃ ἂν ὥσπερ ῥέψαντα τἆλλα ἐφελκύσηται; 

Οὐδαμῶς ἔγωγ', ἔφη, ἄλλοθεν ἢ ἐντεῦθεν.

Οὐκοῦν εἰ τὰ τῶν πόλεων πέντε, καὶ αἱ τῶν ἰδιωτῶν κατασκευαὶ τῆς ψυχῆς 

πέντε ἂν εἶεν. \textit{Republic}, translation \cite{Plato} 544d6--e6

“Are you aware, then,” said I, “that there must be 

as many types of character among men 

as there are forms of government? 

Or do you suppose that constitutions spring from the proverbial oak or rock~and
not from the characters~of the citizens, which, as it were, by their momentum
and weight in the scales~draw other things after them?” 

“They could not possibly come from any other source,” he said. 

“Then if the forms of government are five, 

the patterns of individual souls must be five also.”

 ~

 \item [{[4]}] Ἀλλὰ μὴν πλειόνων γε μέτρων ὂν ἢ ἐλαττόνων,

ὅσωνπερ μέτρων, τοσούτων καὶ μερῶν ἂν εἴη·

καὶ οὕτω αὖ οὐκέτι ἓν ἔσται ἀλλὰ τοσαῦτα 

ὅσαπερ καὶ τὰ μέτρα. \textit{Parmenides}, translation \cite{Plato} 140c8--d2

“But whether it have more measures or less,

it will have as many parts as measures 

and thus one will be no longer one, 

but will be as many as are its measures.”

\item [{[5]}] ΞΕ. Πότερον οὖν

 τὸν πολιτικὸν καὶ βασιλέα καὶ δεσπότην καὶ ἔτ' οἰκονόμον θήσομεν 

ὡς ἓν πάντα ταῦτα  προσαγορεύοντες,

ἢ τοσαύτας τέχνας αὐτὰς εἶναι φῶμεν 

ὅσαπερ ὀνόματα ἐρρήθη; \textit{Statesman}, translation \cite{Plato} 258e8--12

Shall we then assume that the statesman, king, master, and householder too, 

for that matter, are all one, to be grouped under one title, 

or shall we say that 

there are so many arts 

as many names?

\item [{[6]}] καὶ τὸ γήινον ἤδη πᾶν ἀνήλωτο γένος, 

πάσας ἑκάστης τῆς ψυχῆς τὰς γενέσεις ἀποδεδωκυίας,

ὅσα  ἦν ἑκάστῃ προσταχθὲν 

τοσαῦτα  εἰς γῆν σπέρματα πεσούσης, \textit{Statesman}, translation \cite{Plato} 272e2--3

and every soul had fulfilled all its births by

so many seeds falling into the earth 

as many were prescribed to each

\item [{[7]}] οἷαί τε ἔνεισι καὶ ὅσαι, καθορᾷ, 

τοιαύτας καὶ τοσαύτας διενοήθη δεῖν καὶ τόδε σχεῖν. 

εἰσὶν δὴ τέτταρες, 

μία μὲν οὐράνιον θεῶν γένος, 

ἄλλη δὲ πτηνὸν καὶ ἀεροπόρον, 

τρίτη δὲ ἔνυδρον εἶδος, 

πεζὸν δὲ καὶ χερσαῖον τέταρτον. \textit{Timaeus}, translation \cite{Plato} 39e8--40a2

According, then, as Reason perceives Forms existing in the Absolute Living
Creature, 

such and so many 

as many exist therein did He deem that this World also should possess. 

And these Forms are four,

—one the heavenly kind of gods; 

another the winged kind which traverses the air; 

thirdly, the class which inhabits the waters; and 

fourthly, that which goes on foot on dry land.~

\item [{[8]}] τὰ τῶν ψυχῶν γένη, σχημάτων τε 

ὅσα ἔμελλεν αὖ σχήσειν οἷά τε καθ' ἕκαστα εἴδη, 

τὸν μυελὸν αὐτὸν τοσαῦτα καὶ τοιαῦτα  διῃρεῖτο σχήματα εὐθὺς 

ἐν τῇ διανομῇ τῇ κατ' ἀρχάς.  \textit{Timaeus}, translation \cite{Plato} 73c4--6 

He straightway divided the marrow itself, 

in His original division, 

into shapes corresponding in their number and their nature 

to the number and the nature of the shapes which should belong 

to the several kinds of Soul.

\item [{[9]}] ΑΘ. Ἐννοεῖς οὖν ὅτι 

νόμων εἴδη τινές φασιν εἶναι τοσαῦτα

ὅσαπερ πολιτειῶν, 

πολιτειῶν δὲ ἄρτι διεληλύθαμεν ὅσα λέγουσιν οἱ πολλοί; 

\textit{Leges}, translation \cite{Plato} 714b3--5

Are you aware that, according to some,

there are so many kinds of laws 

as many  kinds of constitutions? 

And how many constitutions are commonly recognized 

we have recently recounted 

[[in Leges 712c, they were found to be three]].~

\item [{[10]}] ὅσαι εἰσὶ τὰ νῦν ἡμῖν ἑστίαι διανεμηθεῖσαι  τὸν ἀριθμόν,

ταύτας δεῖν ἀεὶ τοσαύτας εἶναι 

καὶ μήτε τι πλείους γίγνεσθαι μήτε τί ποτε ἐλάττους. 

\textit{Leges}, translation \cite{Plato} 740b3--5

the number of hearths, as many are now appointed by us, 

so many must remain, 

and must never become either more or less.

\item [{[11]}] ἐὰν ἄρα ἀδίκως δοκῇ  ὁ ξένος τὸν ἐπιχώριον τύπτειν, 

τῇ μάστιγι τὸν ξένον

ὅσας ἂν αὐτὸς πατάξῃ 

τοσαύτας δόντες, 

τῆς θρασυξενίας παυόντων· 

\textit{Leges}, translation \cite{Plato} 879e2--5

and if the stranger really appears to have beaten the native unjustly, 

they shall give the Stranger 

so many strokes of the scourge 

as many he himself inflicted,

and make him cease from his foreign forwardness~

\item [{[12]}] Οτι δὲ ἔστιν αἴτια, καὶ ὅτι τοσαῦτα τὸν ἀριθμὸν ὅσα φαμέν, δῆλον·

 τοσαῦτα γὰρ τὸν ἀριθμὸν τὸ διὰ τί περιείληφεν·  

Aristotle \textit{Physics}, 1984~\cite{Hardie1984} 198a14--16

 It is clear then that there are causes, and that the number of them is what we
have stated.  The number is the same as that of the things comprehended under
the question

`why' (translation \cite{Hardie1984}).

\item [{[13]}] ἀναγκαῖον ἦν ἂν τοσαῦτα εἶναι τὰ πάντα γράμματα ὅσαπερ τὰ στοιχεῖα, 

μὴ ὄντων γε δύο τῶν αὐτῶν μηδὲ πλειόνων. Aristotle \textit{Metaphysics} translation, 1924~\cite{Aristotle1924}, 1983~\cite{Aristotle1989} 1000a2--4

Hence just as, if the elements of language~were so many in number, the whole of
literature would be so many as those elements—that is, if there were not two
nor more than two of the same (translation \cite{Aristotle1989}).
\end{enumerate}

\noindent In Euclid’s \textit{Elements}, there are 62 occurences of ‘tosauta,,, hosa’; all refer to
a correspondence between two finite families having the same number. A typical
occurrence is in proposition X.6:

\begin{enumerate}
\item [{[14]}] Ὅσαι γάρ  εἰσιν ἐν τῷ Δ μονάδες, 

εἰς τοσαῦτα ἴσα διῃρήσθω τὸ Α
\end{enumerate}

\subsection{The rejection of the eristic numbers 144d5--7}

\begin{quote}
— Οὐκ ἄρα ἀληθῆ ἄρτι ἐλέγομεν λέγοντες ὡς 

πλεῖστα μέρη ἡ οὐσία νενεμημένη εἴη.  144d5--7

Therefore (‘ara’) what we said just now

—that the Being part was divided into ‘pleista’ parts— was not true, 
\end{quote}

We recall that our interpretation of `pleista' parts of the Being part (144c1--2), presented in Section 4.6. above, was that
\begin{quotation}
the parts of the Being part were infinitely,\\
not simply in multitude, but in number,
\end{quotation}
a statement that, which was declared\textit{ true} if ``number" is taken to be ``the eristic number", now, by 144d5--7, is declared false if ``number" is taken to be ``the dialectic number". With this statement the eristic numbers, defined with unequal units, solely by means of the infinity of the philosophical analogue of the anthyphairesis of the One to the Being, are rejected, as we would expect from the \emph{Philebus} 56d4--e6 passage (Section 4.5), in favor of the dialectic numbers.

\subsection{The dialectic numbers satisfy the following three statements\\
   (a) the dialectical number of the parts of the One is finite,\\
    (b) the dialectical number of the parts of the Being is finite,\\
    (c) the dialectical number of the parts of the One is equal to the dialectical number of the parts of the Being, 144d7--e1}

\subsubsection{Proof of (a)}

First of all, let us make sure that, from what we have already explained in Section 5.4, the number of parts of the One part is finite. This is so because this number is computed in terms of `hapseis', by an appeal to the plus one rule: 
\begin{quotation}
the number of parts of the One part $=$ the hapseis starting from the One part $+1$;
\end{quotation}
 but an `hapsis' is identified with the ratio of two successive parts, and by periodicity and the logos criterion (established in Section 5): \textit{the ratios of successive Logoi is finite in multitude.}

According to Section 5, the presence of the One in the Being, true for the One Being of the second hypothesis (144c2--d4) and described in 138a in terms of hapseis forming a circle, implies that
\begin{quotation}
there is a natural number $k$ such that $\text{One}*\text{Being}=\text{One}_k*\text{Being}_k$.
\end{quotation}
  Since we know (from Section 3) that the sequence $$\text{One} > \text{Being} > \text{One}_1 > \text{Being}_1 > \ldots > \text{One}_k > \text{Being}_k > \ldots $$ is the sequence of remainders of the infinite anthyphairesis of the One to the Being, it now follows from the \textit{Logos Criterion}, as noted in Section 2, that \textit{the anthyphairesis of the dyad (One Being) is periodic}. We now appeal to the plus one rule, in order to compute the numbers involved. Since
\begin{quotation}
the ratios of successive parts-remainders is finite in multitude,
\end{quotation}
namely that 

\begin{quotation}
the hapseis in the One Being form a finite set,
\end{quotation}
it follows that
\begin{quotation}
the number of parts of the One part $=$ the hapseis starting from the One part $+1$,
\end{quotation} 
 because this number is computed in terms of the hapseis-logoi, as set in 148--9 (Section 5).

\subsubsection{(c) is stated in 144d7--e1}
\begin{quote}
οὐδὲ γὰρ πλείω τοῦ ἑνὸς νενέμηται,\\ ἀλλ' ἴσα, ὡς ἔοικε, τῷ ἑνί· \\
For it is not divided, you see, into any more parts than one,\\ but, as it seems, into the \textbf{equal} number as the One' 144d7--e1.
\end{quote}

Proof of (c). But because of periodicity in anthyphairesis, it follows immediately that
\begin{quotation}
the hapseis starting from the One part\\
are equal to\\
the hapseis starting from the Being part 
\end{quotation}
(since if $\text{One}*\text{Being}=\text{One}_k*\text{Being}_k$, then $\text{Being}*\text{One}_k=\text{Being}_k*\text{One}_{k+1}$).

Hence, we conclude 
\begin{quotation}
the number of the parts of the One part\\
is equal to\\
the number of the parts of the Being part
\end{quotation}
exactly as claimed in 144d4--5.

\subsubsection{(b) follows from (a) and (c).}

\subsubsection{Strong support for the finiteness of the whole dialectical number is given in the by Plato's \emph{Philebus} 16c10--e1}

The finiteness of the whole number in every Platonic Being
\begin{quote}
 ‘This being the way in which these things are set in order (‘diakekosmemenon’),

we must always (‘aei’) assume (‘themenous’) that 

there is in every~case~(‘hekastote’) one Idea (‘mian idean’) for everything
(‘peri pantos’)

and must look for (‘zetein’) it—for we shall find that it is there—

and if we get a grasp (‘metalabomen’) of this,

we must look next for two, if there be two, and if not, for three 

or some other number (‘arithmon’); 

and again we must treat each of those units in the same way 

(‘kai ton hen ekeinon hekaston palin hosautos’), 

until we can see 

not only (‘monon’) that the original One (‘to kat’ archas hen’) is one and many
and infinite, 

but just how many (‘alla kai hoposa’) it is. 

and we must not apply the idea of infinite (‘ten tou apeirou idean’) to
plurality (‘pros to plethos’)

until (‘prin’) we have a view of its whole number (‘ton arithmon autou panta’,
16d8) 

between (‘metaxu’) infinity and one;’ \textit{Philebus}, translation \cite{Plato} 16c10--e1
\end{quote}

\textit{Comment}. There results a finite number, by the hapseis rule in the \emph{Parmenides} 148d--149d: \\
`arithmos pas'$=$hapseis $+1$. Thus\\
 `arithmos pas' (16d8) is the number $n+1$. \\
 The $n+1$ many equalized units of this number are \\
 `ta panta' (16e2), the equalized units $e_1,e_2,\ldots, e_{n+1}$ of `arithmos pas'\\
 `hen hekaston ton panton' (16e1--2) is each of the equalized units $e_1,e_2,\ldots,e_{n+1}$'\\
  `hen hekaston ton ekeinon' (16d4--5) coincides with `hen hekaston ton panton' `ton ekeinon' (16d4) means `ton panton'\\
  \begin{quotation}
  `kai ton hen ekeinon palin hosautos' (16d4--5)
\end{quotation}  
 refers to the fact that, because of periodicity, each of the units-parts of `arithmos pas' has the same number of parts-units (with the same argument that the number of parts of the One is equal to the number of parts of the Being in 144d4--e3.   
Note this statement has been erroneously taken to mean departure from binary division). The equalization of the units implies the correct concept of presence of the One in the parts of the Platonic Idea (138a, 144d5--e3).

\subsubsection{Strong support for the finiteness of dialectical number is given in Proclus' \emph{Platonic Theology} 4, 94,26--95,4 passage}

\begin{quote}
It is clear that this infinity of multitudes must not be transferred to number

(’Καὶ μάλιστα δῆλον ὡς οὐ χρὴ τὴν ἀπειρίαν ταύτην ἐπὶ τὸ ποσὸν μεταφέρειν),

for how is it possible that number be infinite, 

since this infinity is against the nature of number’;

(πῶς γὰρ άπειρος ἀριθμός, 

αὐτῆς τῆς ἀπειρίας πρὸς τὴν τοῦ ἀριθμοῦ φύσιν διαμαχομένης;)

and how can the parts of the Being be equal to the parts of the One?

(Πῶς δὲ ἴσαι τοῖς τοῦ ὄντος κέρμασιν αἱ μοῖραι τοῦ ἑνός;)

Because in the infinite multitudes equality does not exist

(Ἐν γὰρ τοῖς ἀπείροις τὸ ἴσον οὐκ ἔστιν.)

Proclus, \textit{Platonic Theology}, translated 1816 \cite{ProclusTheology}, translated by the author, 4, 94,26--95,4 translation
by the author.
\end{quote}

\subsection{The paradoxical statement about dialectic numbers: the multitude of the parts of the Being is infinite, while the (dialectic) number of the parts of the Being is finite}

Now the falsity of this claim can only mean that
\begin{quotation}
the parts of the Being part are finite in number.
\end{quotation}

 The combination of the statements obtained in 142b--143b (Section 3) and in 144d4--7 (present Section 6.5) bring us in a rather difficult and paradoxical situation: 
the family of parts of the One has been shown in Section 3 to satisfy an infinity condition, namely that
\begin{quotation}
the parts of the One is infinite in multitude, and,
\end{quotation}

in seeming conflict to this infinity condition, is now claimed, in the present Section 6.5, to satisfy a finiteness condition, namely that 

\begin{quotation}
the parts of the One are finite in dialectical number.
\end{quotation}

\subsection{The description of the generation of dialectic number Two, given in 143c1--d5, in terms of Logos can now be understood as the dialectical and not the eristic number Two}

Right after the generation of the anthyphairetic sequence of parts-remainder, the number Two was defined in 143c1--d5 (in Section 4):
\begin{quotation}
‘Τί οὖν; ἐὰν προελώμεθα αὐτῶν 

εἴτε βούλει τὴν οὐσίαν καὶ τὸ ἕτερον εἴτε τὴν οὐσίαν καὶ τὸ ἓν εἴτε τὸ ἓν καὶ τὸ
ἕτερον,

ἆρ' οὐκ ἐν ἑκάστῃ τῇ προαιρέσει προαιρούμεθά τινε ὣ ὀρθῶς ἔχει καλεῖσθαι
ἀμφοτέρω; 

— Πῶς;

 — Ὧδε· ἔστιν οὐσίαν εἰπεῖν(143c5);

 — Ἔστιν. 

— Καὶ αὖθις εἰπεῖν (143c5) ἕν; 

— Καὶ τοῦτο.

— Ἆρ' οὖν οὐχ ἑκάτερον αὐτοῖν εἴρηται (143c6); 

— Ναί. 

— Τί δ' ὅταν εἴπω (143c7) οὐσία τε καὶ ἕν, ἆρα οὐκ ἀμφοτέρω; 

— Πάνυ γε.

— Οὐκοῦν καὶ ἐὰν οὐσία τε καὶ ἕτερον ἢ ἕτερόν τε καὶ ἕν,

καὶ οὕτω πανταχῶς ἐφ' ἑκάστου ἄμφω λέγω (143c9); 

— Ναί. 

—  Ὣ δ' ἂν ἄμφω ὀρθῶς προσαγορεύησθον, 

ἆρα οἷόν τε ἄμφω μὲν αὐτὼ εἶναι, δύο δὲ μή; 

— Οὐχ οἷόν τε. 

— Ὣ δ' ἂν δύο ἦτον, ἔστι τις μηχανὴ μὴ οὐχ ἑκάτερον αὐτοῖν ἓν εἶναι; 

— Οὐδεμία.

— Τούτων ἄρα ἐπείπερ σύνδυο ἕκαστα συμβαίνει εἶναι, καὶ ἓν ἂν εἴη ἕκαστον. 

— Φαίνεται.’
\end{quotation}

At this early point of the second hypothesis it was not possible to understand the passage 143c1--d5. In support of this view we consider two classical renderings of the passage, one by H.N. Fowler \cite{Plato}, the other by Allen \cite{Allen1997}.

The translation of the passage by H. N. Fowler.

\begin{quote}
‘Well, then, if we make a selection (‘proelometha’) among them,

whether we select (‘proairoumetha’) 

Being (‘ousian’) and the Other (‘heteron’), 

or Being (‘ousian’) and One, 

or One and the Other (‘heteron’),

in each instance we select (‘en hekastei tei proairesei’) two parts 

which may properly be called both (‘amphotero’)?’

‘What do you mean?’ 

‘I will explain.

We can speak of (‘eipein’, 143c5) Being (‘ousian’)?’ 

‘Yes.’

‘And we can also speak of (‘eipein’,143c5) One?’ 

‘Yes, that too.’

‘Then have we not spoken of (‘eiretai,143c6’) each of them?’ 

‘Yes.’

‘And when I speak of Being (‘ousia’) and One, 

do I not speak of (‘eipo’, 143c7) both?’ 

‘Certainly.’

‘And also when I speak of (‘lego’, 143c9) Being (‘ousia’) and Other (‘heteron’),
or Other (‘heteron]) and One, 

in every case I speak of (‘lego’, 143c9) each pair as both?’ 

‘Yes.’

‘If things are correctly called both (‘ampho’), 

can they be both (‘ampho’) without being two (‘duo’)?” ‘They cannot.’

‘And if things are two (‘duo’), must not each of them be one (‘hen’)?’
‘Certainly.’

‘Then since the units of these pairs are together two (‘sunduo’), each
‘(hekaston’) must be individually one (‘hen’).’ 

‘That is clear.
\end{quote}

The reader cannot fail to notice that there is a volley of several forms of `legein' (`eipein'143c5, `eipein'143c5, `eiretai'143c6, `eipo'143c7, `lego'143c9) in short succession, rendered as `speaking of'; it might be attributed improbably to an awkwardness on Plato's part, or is there some other reason?

Allen, 1997~\cite{Allen1997} p.262, realizes the awkwardness and attempts to explain it by the lack of dual case in English!

\begin{quote}
‘Parmenides next proceeds to show that Unity, taken apart from its own being,
implies number (143c--d). 

The exact force of his argument \textit{cannot be reproduced in English. 
}
Greek possesses, as English does not, a dual as well. as a singular and plural;
when Parmenides argues that since it is possible 

to \textit{mention} (143c5) Unity and 

to \textit{mention} (143c5) Being, 

each of two has been \textit{mentioned} (143c6), 

the English ``two" is more explicit than the text, 

which contains only the genitive dual ‘autoin’. 

It is from this feature in the syntax of his language 

that Parmenides goes on to infer that 

both have been \textit{mentioned} (143c7), 

and that since both have been \textit{mentioned} (143c9), 

two have been \textit{mentioned} (143c9).’ 

[emphasis added by present author] .
\end{quote}

On the other hand, the definition of the number two, with units One and Being, has a problem, since, as noted, the two units are unequal, and thus the number two as it stood would be an eristic and not a dialectical number. It was only after 144c2--e3 that these units are \textit{equalized}, and the definition of Two becomes dialectical. 
This \textit{equalization} is possible only because the number is considered as not only the dyad One Being, but \textit{with the hapsis} between them, namely \textit{the ratio/Logos} between them; equalisation is effected by the logos-hapsis of periodicity that joins together the One with the Being. In retrospect, we realize that this is exactly what Plato was trying to convey in 143c1--d5 with the repeated, then unexplained, references to `legein' and Logos.

\subsection{Finiteness of dialectical number in Plotinus, \emph{Enneades} 6, 6, 17, 1--10: Only a finite dialectic number can be generated in a Platonic Being, because there is only a finite number of `hapseis'-'logoi'}

It follows that the generation of eristic numbers attempted in 143d5--144c2, specifically the generation of the number $n$ times $m$, $n \cdot m$, if the numbers $n, m$ are already generated, is not in general valid. In fact if $n$ is number generated, 
there is no guarantee that the number $2n$ can be generated. Indeed, suppose that $e_1 * e_2 * \ldots * e_n$ has been generated; in order that $e_{n+1} * e_{n+2} * \ldots * e_{2n}$ can also be generated, we would need to have that there are further `hapseis', 
since additional units are generated only with additional `hapseis'. But there is absolutely no reason why we would have a new hapsis $*$ at the point $e_n * e_{n+1}$, so as to form $e_1 * e_2 * \ldots * e_n * e_{n+1} * e_{n+2} * \ldots * e_{2n}$.

A strong confirmation of our interpretation that not every number can be represented in a Platonic Idea, but there is a greatest number that can be represented in a Platonic Idea, and this number is determined by the length of the period of the logoi is found in the Plotinus' \emph{Enneades}. Plotinus clearly states, in the \emph{Enneades} 6,6,17,1--11 passage, 
\textit{the finiteness} of number, and identifies its cause: the finite number of true \textit{contacts} (`sunapsanta', `prosapsais')
\begin{quote}
‘But what about the number called infinite (‘apeiros’)? 

For these ‘logoi’ provide a limit (‘peras’). 

And this is correct, if it is going to be a number; 

for the infinite (‘to apeiron’) clashes with number. 

Why, then, do we say ‘The number is infinite (‘apeiros’)’? 

Is it with number as it is when we say a line is infinite (‘apeiron’)?

--but we say a line is infinite (‘apeiron’) 

not because there is any line of this kind 

but because it is possible with the longest line, 

that of the universe for instance, to think (‘epinoesai’) of a longer. 

For having obtained knowledge (‘gnosthentos’) how much (‘hosos’) a number is, 

it is possible to double it in thought (‘tei dianoiai’) 

without connecting (‘sunapsanta’) it to the known number. 

For how you could attach (‘prosapsais’) 

a thought and mental image (‘noema kai phantasma’) 

which is only in you to things that really exist?’ (translation \cite{Plotinus})
\end{quote}

\subsection{The connection of the dialectic numbers with Logos in the \emph{Epinomis} 977b9--e2 confirms the need of Logos, namely the Logos Criterion for the periodic anthyphairesis, nrcessary for the definition of dialectic numbers}

In the \textit{Epinomis}, translation \cite{Plato} 977b--e passage where logos and number are insistently
correlated.

[1] 

Ἔτι δὲ σμικρὸν ἐπανελθόντες πως τοῖς λόγοις ἀναμνησθῶμεν

Moreover, turning back some little way,  

Let us recollect the logoi (τοῖς λόγοις ἀναμνησθῶμεν, 977b9--c1, rendered by Lamb
‘could not tell’) 

[2] 

ὅτι καὶ μάλ' ὀρθῶς ἐνοήσαμεν ὡς, 

εἴπερ ἀριθμὸν ἐκ τῆς ἀνθρωπίνης φύσεως ἐξέλοιμεν, 

οὐκ ἄν ποτέ τι φρόνιμοι γενοίμεθα.

οὐ γὰρ ἂν ἔτι ποτὲ ψυχὴ τούτου τοῦ ζῴου πᾶσαν ἀρετὴν λάβοι σχεδόν,

ὅτου λόγος ἀπείη· 

on how entirely right we were in conceiving that 

if we should deprive human nature of number (ἀριθμὸν, 977c1) 

we should never attain to any knowledge.

For then the soul of that creature which could not posess logos (λόγος, 977c4,
rendered by Lamb ‘clearly tell of’)  

would never any more be able, one may say, to attain virtue (ἀρετὴν) in general;

[3] 

ζῷον δὲ ὅτι 

μὴ γιγνώσκοι δύο καὶ τρία 

μηδὲ περιττὸν μηδὲ ἄρτιον, 

ἀγνοοῖ δὲ τὸ παράπαν ἀριθμόν, 

οὐκ ἄν ποτε διδόναι λόγον ἔχοι 

περὶ ὧν αἰσθήσεις καὶ μνήμας [ἔχοι] μόνον εἴη κεκτημένον,

τὴν δὲ ἄλλην ἀρετήν, ἀνδρείαν καὶ σωφροσύνην, οὐδὲν ἀποκωλύει. 

and the creature that 

did not know two and three, or odd or even, 

and was completely ignorant of number (ἀριθμόν, 977c6), 

could never clearl provide logos) (διδόναι λόγον, 977c6, rendered by Lamb ‘ tell
of things’) 

about which it had only acquired sensations and memories. 

From the attainment of ordinary virtue—courage and temperance—

it is certainly not debarred: 

[4] 

στερόμενος δὲ ἀληθοῦς λόγου σοφὸς οὐκ ἄν ποτε γένοιτο, 

ὅτῳ δὲ σοφία μὴ προσείη, 

πάσης ἀρετῆς τὸ μέγιστον μέρος, 

οὐκ ἂν ἔτι τελέως ἀγαθὸς γενόμενος εὐδαίμων ποτὲ γένοιτο. 

οὕτως ἀριθμὸν μὲν ἀνάγκη πᾶσα ὑποτίθεσθαι· 

but if a man is deprived of true logos (ἀληθοῦς λόγου, 977d2, rendered by Lamb’
true telling’) he can never become wise (σοφὸς), and he who has not the
acquirement of wisdom—the greatest part of virtue as a whole—

can no more achieve the perfect goodness which may make him happy. 

Thus it is absolutely necessary to postulate number (ἀριθμὸν, 977d5); 

[5] 

διότι δὲ τοῦτο ἀνάγκη, λόγος ἔτι πλείων πάντων γίγνοιτ' ἂν τῶν εἰρημένων. 

ἀλλὰ καὶ ὁ νῦν ὀρθῶς ῥηθήσεται, 

ὅτι καὶ τὰ τῶν ἄλλων τεχνῶν λεγόμενα, 

ἃ νυνδὴ διήλθομεν ἐῶντες εἶναι πάσας τὰς τέχνας,

οὐδὲ τούτων ἓν οὐδὲν μένει, 

πάντα δ' ἀπολείπεται τὸ παράπαν,

ὅταν ἀριθμητικήν τις ἀνέλῃ.  

and why this is necessary can be shown by a still fuller logos (λόγος, 977d6,
rendered by Lamb ‘argument’) than any that has been advanced. But here is one
that will be particularly correct—that of the properties of the other arts,
which we recounted just now in granting the existence of all the arts not a
single one can remain, but all of them are utterly defective, when once you
remove numeration (ἀριθμητικήν, 977e2). \textit{Epinomis}, translation \cite{Plato} 977b9--e2

[translated by W.R.M. Lamb. 1925 \cite{Plato}]
%Cambridge, MA, Harvard University Press; London,
%William Heinemann Ltd. 1925, modified by the author]

In this remarkable passage in the \emph{Epinomis} Plato refers to the recollection of the  and correlates dialectical number with logos four times. In view of our analysis, in the \emph{Parmenides} 142b1--144e3 connection between anthyphairetic logos and number, we will have no difficulty in identifying the logos that occurs in the \emph{Epinomis} passage with the anthyphairetic logos criterion for anthyphairetic periodicity, from which true dialectical numbers, with equal (equalized) units are formed, and true knowledge is achieved (as shown in the \emph{Sophist} and \emph{Statesman}). 
The recollection of logoi, as in the \emph{Meno, Phaedo} and \emph{Phaedrus}, 
refers to the repetition of the anthyphairetic logos after a full revolution instituting equalization of units and acquisition of knowledge (episteme) (cf. Negrepontis, 2024~\cite{Negrepontis2024a}).

\section{Plato's Indivisible Line, considered mysterious by modern Platonists, is identified with Plato's Intelligible Being}

\subsection{Indivisible Lines, a mystery for the Platonists}

As was acknowledged by A. E. Taylor, 1960~\cite{Taylor1960} and more recently by M. Rashed, 2013~\cite{Rashed2013b}, p.106--109, Plato's and Xenocrates' indivisible line remains a mystery for the scholars on Plato up to our day. We will outline our interpretation that identifies the Indivisible Line with Plato's intelligible Being. For more details the reader is referred to Negrepontis, 2024~\cite{Negrepontis2024b}.

\subsection{Understanding the nature of the dialectical numbers introduced in an intelligible Being, called by Aristotle (eidetic numbers), and their connection with the kinds (eide) (\emph{Metaphysics} 987b21--22, 1081a21, 1086a5,8, 1088b34, 1090b35) provides a natural reason to call the intelligible Being indivisible after the completion of the period of its anthyphairesis}

The term indivisible line is in analogy to the paradox of the infinite multitude of the parts of the One (of the second hypothesis in the \emph{Parmenides}) vs. the finite dialectical number of parts of the One. In the same way, the dyad One, Being of the paradigmatical intelligible Being is divisible ad infinitum, by infinite anthyphairesis, but the kinds (eide) which are identified with (the units of) the dialectical numbers, are finite. The identification of the dialectical numbers with the kinds is seen in Aristotle' description in the \emph{Metaphysics}, in which it is stated that the Platonic kinds are numbers (\emph{Metaphysics} 987b21--22), while he calls the dialectical numbers as eidetic/kind-like (\emph{Metaphysics} 1081a21, 1086a5,8, 1088b34, 1090b35); 
cf. Simplicius, \emph{Commentary to Aristotle's Physics} 140,6--13, 142, 16–18. Thus, after the completion of a period, the anthyphairetic division continues ad infinitum, but there are no more new kinds, and it is in exactly this sense that we describe an intelligible Being as an \textbf{indivisible line}. Proclus warns us about the symbolic meaning of the term:
\begin{quote}
But Plato, for the sake of concealment (di'epikrupsin), employed mathematical names (onomaton), as screens (parapetasmasin) of the truth of the things, in the same manner as theologists employed fables, and the Pythagoreans symbols (tois sumbolois). \emph{Proclus, Commentary to Timaeus} 2,245,23--246,7, translation Thomas Taylor 1820 \cite{Proclus1820}.
\end{quote}

\subsection{In the description of the Noble Sophistry in the \emph{Sophist} 226b1--231b8, an example of an intelligible Being, by the Name and Logos, it is declared indivisible after reaching the Logos Criterion and anthyphairetic periodicity}

The acquisition of knowledge of the intelligible Being Noble Sophistry in the \emph{Sophist} 226b1--231b8 and the description of the Division of every intelligible Being in the \emph{Phaedrus} 277b5--c6 by Nane and Logos, 
in conjunction with our anthyphairetic interpretation of the method of Name and Logos in the \emph{Sophist} (Negrepontis 2012~\cite{Negrepontis2012}), quickly lead us to the conclusion that an intelligible Being becomes ``indivisible" exactly at the moment that the Logos Criterion establishing anthyphairetic periodicity is instituted. 
Indeed, one step before achieving the Logos, the Eleatic Stranger in the \emph{Sophist} is asking the crucial question:
\begin{quote}
[Stranger] But we must examine further and see whether it is already (hede, 229d5) an \textbf{indivisible} (atomon, 229d5) whole or still admits of division important enough to have a name. 229d5--6
\end{quote}

Similarly, Socrates in the \emph{Phaedrus}:
\begin{quote}
Socrates. Before someone obtains knowledge (eidei) of the truth about each of those he sets logoi (legei, 277b6) [speaks] (or writes, he must be able to define (horizesthai, 277b6) everything; and having defined 
(horisamenos, 277b7), he must obtain knowledge (epistethei, 277b8) by dividing (temnein, 277b7) according to kinds (kat' eide, 277b7) till the \textit{indivisible} (mechri tou atmetou, 277b7, translation \cite{Plato}).
\end{quote}

\subsection{The comparison of the ``mysterious" Aristotle's \emph{Metaphysics} 992a19--22, with the One in the first and in the second hypothesis in the \emph{Parmenides}, naturally suggests that the One of the second hypothesis, the paradigmatical intelligible being, is an indivisible line}

The heretofore mysterious Aristotle's \emph{Metaphysics} 992a19--22 passage, according to which Plato does not regard the geometric point suitable for the foundation of Geometry, but instead he opts for the indivisible line, becomes fully understandable if we compare it with the possibility of each of the two Hypotheses in Plato's \emph{Parmenides} to generate Geometry. 

Thus, the One of the first hypothesis, a partless One, exactly like a geometric point 137c5--6, is declared unfit to generate the two basic geometric constructions, the straight and the circular (\emph{Parmenides} 137d8--e4), and thus unfit to serve as the foundation of Geometry,

While the One of the second hypothesis, the paradigmatical intelligible Being, does generate the two basic geometric constructions, 
the straight and the circular (\emph{Parmenides} 145b3--b5), and thus fit to serve as the foundation of Geometry,

Comparing Plato's \emph{Parmenides} two Hypotheses and Aristotle's \emph{Metaphysics} 992a19--22, we naturally conclude that the \textbf{indivisible line} (in the \emph{Metaphysics} passage) precisely corresponds to Plato's intelligible Being.

On the relation of the indivisible line with the \textit{Republic} "hypothesis-free"
(anhupotheron) consult Negrepontis 2019 \cite{Negrepontis2019}.

\subsection{We conclude that the indivisible line coincides with the intelligible Being.}

\section{Our interpretation of the self-similar Oneness of the One of the second hypothesis in the \emph{Parmenides} 144e1--3, 145a2--4, and of Plato's statement that Not-Being is a Being in the \emph{Sophist} 257d7--258c5, 258d5--e5}

\subsection{The crucial passage on the equalization of the One and the Being \emph{Parmenides} 144e1--3}

In the second hypothesis of the \emph{Parmenides}, the initial dyad of the paradigmatical intelligible Being, the One and the Being, are equalized by means of the periodic anthyphairesis and the resulting introduction of dialectical numbers.
\begin{quote}
οὔτε τὸ ἓν τοῦ ὄντος, 

ἀλλ' ἐξισοῦσθον δύο ὄντε ἀεὶ παρὰ πάντα.

— Παντάπασιν οὕτω φαίνεται. 

‘for Being is not wanting to the One, 

nor the One to Being, 

but being two they are equalized throughout.” 

“That is perfectly clear.” ’ 144e1--3
\end{quote}
\begin{quote}
— Τὸ ἓν ἄρα ὂν 

ἕν τέ ἐστί που καὶ πολλά, καὶ 

ὅλον καὶ μόρια, 

καὶ πεπερασμένον καὶ ἄπειρον πλήθει. 

— Φαίνεται.  

“Then (‘ara’) the One Being (‘to hen on’) is, apparently, 

both one (‘hen’) and many (‘polla’), 

a whole (‘holon’) and parts (‘moria’), 

finite (‘peperasmenon’) and of~infinite~multitude (‘apeiron plethei’).” 

“So it appears.” 145a2--4 (translation in \cite{Plato})
\end{quote}

The One of the second hypothesis is

(i) many, in fact infinite in multitude, by the infinite anthyphairetic division of the initial indefinite dyad One and Being;

(ii) finite in number, by the anthyphairetic periodicity (finiteness of hapseis/logoi) and the plus one rule: number of units$=$ hapseis$+1$; and,

(iii) One, because, although One and Being is as an indefinite dyad, the equalization of the dialectical number of the One and of the Being, in fact the equalization of the dialectical number of the One and of the dialectical number with any anthyphairetic remainder $\text{One}_k$ or $\text{Being}_k$ for all indices $k$, Being the intelligible Being One into One, in fact it turns it into a \textit{self-similar One}, 
in which each of the infinitely many parts is equalized to the whole. By exactly the same argument, every part of the One generated by the anthyphairesis of the dyad One, Being is equalized to the One; indeed the dialectical number of the sequence of all the parts of $\text{One}_k$, or of $\text{Being}_k$, is equal to the length of the period of the anthyphairesis of $\text{One}_k$ to $\text{Being}_k$, or of $\text{Being}_k$ to $\text{One}_{k+1}$, respectively, for any $k=1,2,\ldots,$ and by anthyphairetic periodicity, 
this dialectical number is equal to the dialectical number of the sequence of all part οὔτε γὰρ τὸ ὂν τοῦ ἑνὸς ἀπολείπεται

\subsection{The statement of ``not-Being is" in Plato' \emph{Sophist} 257d7--258c5, 258d5--e5}

\begin{quote}
ΞΕ. Ἔστι τῷ καλῷ τι θατέρου μόριον ἀντιτιθέμενον;

ΘΕΑΙ. Ἔστιν.

ΞΕ. Τοῦτ' οὖν ἀνώνυμον ἐροῦμεν ἤ τιν' ἔχον ἐπωνυμίαν;

ΘΕΑΙ. Ἔχον· ὃ γὰρ μὴ καλὸν ἑκάστοτε φθεγγόμεθα,

τοῦτο οὐκ ἄλλου τινὸς ἕτερόν ἐστιν ἢ τῆς τοῦ καλοῦ φύσεως.257d7--11

Stranger Is there a part of the other which is opposed to the beautiful?

Theaetetus There is.

Stranger Shall we say that this is nameless or that it has a name?

Theaetetus That it has one; for that which in each case we call not-beautiful is
surely the other of the nature of the beautiful and of nothing else.
\end{quote}
The Beautiful and the not-Beautiful are other (heteron) to each other. two parts that are other (heteron) to each other In the second hypothesis of the \emph{Parmenides} the corresponding dyad of others is the dyad One and Being. The difference between the \emph{Parmenides} dyad One, Being and the \emph{Sophist} dyad Beautiful, not-Beautiful is just in name.

\begin{quote}
ΞΕ. Ἴθι νυν τόδε μοι λέγε.

ΘΕΑΙ. Τὸ ποῖον;

ΞΕ. Ἄλλο τι τῶν ὄντων τινὸς ἑνὸς γένους ἀφορισθὲν

καὶ πρός τι τῶν ὄντων αὖ πάλιν ἀντιτεθὲν 

οὕτω συμβέβηκεν εἶναι τὸ μὴ καλόν;

ΘΕΑΙ. Οὕτως.

ΞΕ. Ὄντος δὴ πρὸς ὂν ἀντίθεσις, ὡς ἔοικ', εἶναί τις συμβαίνει τὸ μὴ καλόν.

ΘΕΑΙ. ᾿Ορθότατα.257d12--e8

Stranger Now, then, tell me something more.

Theaetetus What?

Stranger Does it not result from this that the not-beautiful is a distinct part
of some one class of being and also, again, opposed to some class of being?

Theaetetus Yes.

Stranger Then, apparently, it follows that the not-beautiful is a contrast of
being with being.

Theaetetus Quite right.
\end{quote}
Thus, the Beautiful is an intelligible Being, it possesses an heteron, the not-Beautiful, and the dyad Beautiful and not-Beautiful is a dyad in the philosophic analogue of periodic anthyphairesis. Because of periodic anthyphairesis, the definition of the dialectical numbers, and the criterion of equalization of parts in terms of the dialectical number of their parts, the Beautiful and the not-Beautiful are \textbf{equalized} and both are intelligible Beings.

\begin{quote}
ΞΕ. Τί οὖν; κατὰ τοῦτον τὸν λόγον ἆρα 

μᾶλλον μὲν τὸ καλὸν ἡμῖν ἐστι τῶν ὄντων, 

ἧττον δὲ τὸ μὴ καλόν;

ΘΕΑΙ. Οὐδέν.

ΞΕ. Ὁμοίως ἄρα τὸ μὴ μέγα καὶ τὸ μέγα αὐτὸ εἶναι λεκτέον;

ΘΕΑΙ. Ὁμοίως.

ΞΕ. Οὐκοῦν καὶ τὸ μὴ δίκαιον τῷ δικαίῳ 

κατὰ ταὐτὰ θετέον πρὸς τὸ μηδέν τι μᾶλλον εἶναι θάτερον θατέρου;

ΘΕΑΙ. Τί μήν; 257e9--258a6

Stranger Then, in that case, according to logos, 

the beautiful is more and the not-beautiful less one of the beings?

Theaetetus Not at all.

Stranger Hence the not-great must be said to be no less truly than the great?

Theaetetus No less truly.

Stranger And so we must recognize the same relation between the just and the
not-just, in so far as neither has any more being than the other?

Theaetetus Of course.
\end{quote}

Because of Logos, namely because of the Logos Criterion and the resulting anthyphairetic periodicity, the two parts of the initial dyad of the intelligible Being, the Beautiful and the Not-Beautiful are \textit{equally intelligible Beings}

\begin{quote}
ΞΕ. Καὶ τἆλλα δὴ ταύτῃ λέξομεν, 

ἐπείπερ ἡ θατέρου φύσις ἐφάνη τῶν ὄντων οὖσα, 

ἐκείνης δὲ οὔσης ἀνάγκη δὴ καὶ τὰ μόρια αὐτῆς μηδενὸς ἧττον ὄντα τιθέναι.

ΘΕΑΙ. Πῶς γὰρ οὔ;

ΞΕ. Οὐκοῦν, ὡς ἔοικεν, 

ἡ τῆς θατέρου μορίου φύσεως καὶ τῆς τοῦ ὄντος πρὸς ἄλληλα ἀντικειμένων ἀντίθεσις
οὐδὲν ἧττον, εἰ θέμις εἰπεῖν, αὐτοῦ τοῦ ὄντος οὐσία ἐστίν, 

οὐκ ἐναντίον ἐκείνῳ σημαίνουσα 

ἀλλὰ τοσοῦτον μόνον, ἕτερον ἐκείνου.

ΘΕΑΙ. Σαφέστατά γε. 258a7--b4

Stranger And we shall, then, by Logos (lexomen), the same holds of other things,
since the nature of the other is proved to possess real being; and if it has
being, we must necessarily ascribe being in no less degree to its parts also.

Theaetetus Of course.

Stranger Then, as it seems, the opposition of the nature of a part of the other,
and of the nature of being, when they are opposed to one another, is no less
truly existence than is being itself, if it is not wrong for me to say so, for
it signifies not the opposite of being, but only the other of being, and
nothing more.

Theaetetus That is perfectly clear.
\end{quote}

But the equalization is not restricted only in the equalization of the not-Beautiful to the Beautiful, but it extends to every part generated from the infinite periodic anthyphairetic division. Thus every part is equalized to the Beautiful, and therefore constitutes an intelligible Being.

\begin{quote}
ΞΕ. Τίν' οὖν αὐτὴν προσείπωμεν;

ΘΕΑΙ. Δῆλον ὅτι τὸ μὴ ὄν, ὃ διὰ τὸν σοφιστὴν ἐζητοῦμεν, αὐτό ἐστι τοῦτο.258b5--7

Stranger Then what shall we call this?

Theaetetus Evidently this is precisely not-being, which we were looking for
because of the sophist.
\end{quote}

According to the Name \& Logos of the \emph{Sophist} (Negrepontis 2012 \cite{Negrepontis2012}, 2018 \cite{Negrepontis2018}), the Sophist is a Not-Being with respect to its opposite/heteron 
Demologue, and since the anthyphairesis of the dyad Demologue, Sophist is clearly periodic, it forms an intelligible Being.

\begin{quote}
ΞΕ. Πότερον οὖν, ὥσπερ εἶπες, ἔστιν οὐδενὸς τῶν ἄλλων οὐσίας ἐλλειπόμενον,

καὶ δεῖ θαρροῦντα ἤδη λέγειν ὅτι 

τὸ μὴ ὂν βεβαίως ἐστὶ τὴν αὑτοῦ φύσιν ἔχον, 

ὥσπερ τὸ μέγα ἦν μέγα καὶ τὸ καλὸν ἦν καλὸν 

καὶ τὸ μὴ μέγα <μὴ μέγα> καὶ τὸ μὴ καλὸν <μὴ καλόν>, 

οὕτω δὲ καὶ τὸ μὴ ὂν κατὰ ταὐτὸν ἦν τε καὶ ἔστι μὴ ὄν, 

ἐνάριθμον τῶν πολλῶν ὄντων εἶδος ἕν;

ἤ τινα ἔτι πρὸς αὐτό, ὦ Θεαίτητε, ἀπιστίαν ἔχομεν;

ΘΕΑΙ. Οὐδεμίαν. 258b8--c5

Stranger And is this, as you were saying, as fully endowed with being as
anything else, and shall we henceforth say with confidence that not-being has
an assured existence and a nature of its own?

Just as we found that the great was great and the beautiful was beautiful, the
not-great was not-great and the not-beautiful was not-beautiful, shall we in
the same way say that not-being was and is not-being, to be counted as one
class among the many classes of being? Or have we, Theaetetus, any remaining
distrust about the matter?

Theaetetus None whatever.
\end{quote}

The not-Being is a Being and the dialectical number of the not-Beings is finite.

\begin{quote}
ΞΕ. Ἡμεῖς δέ γε

οὐ μόνον τὰ μὴ ὄντα ὡς ἔστιν ἀπεδείξαμεν,

ἀλλὰ καὶ τὸ εἶδος ὃ τυγχάνει ὂν τοῦ μὴ ὄντος ἀπεφηνάμεθα·

τὴν γὰρ θατέρου φύσιν

ἀποδείξαντες

οὖσάν τε καὶ κατακεκερματισμένην ἐπὶ πάντα τὰ ὄντα πρὸς ἄλληλα,

τὸ πρὸς τὸ ὂν ἕκαστον μόριον αὐτῆς ἀντιτιθέμενον ἐτολμήσαμεν εἰπεῖν

ὡς αὐτὸ τοῦτό ἐστιν ὄντως τὸ μὴ ὄν.

ΘΕΑΙ. Καὶ παντάπασί γε, ὦ ξένε, ἀληθέστατά μοι δοκοῦμεν εἰρηκέναι.258d5--e5

Stranger But we have not only pointed out that things which are not-Beings  are
in fact Beings, but we have even shown the form of Being that not-Being is; 

for we have proved that the nature of the other is a Being and is distributed in
small bits

throughout all Beings in their relations to one another, and we have ventured to
say that each part of the other which is opposed to the Being, exactly the
not-Being is really a Being.

Theaetetus And certainly, Stranger, I think that what we have said is perfectly
true. 
\end{quote}

[Plato. \emph{Plato in Twelve Volumes}, Vol. 12 translated by Harold N. Fowler. 
Cambridge, MA, Harvard University Press; London, William Heinemann Ltd. 1921 \cite{Plato}, with modifications by the author]

\subsection{Our interpretation of the statement the not-Being is a Being in the \emph{Sophist}}

According to our interpretation of the dialogues \emph{Theaetetus, Meno, Sophist, Parmenides}, an intelligible Being is the philosophical analogue of a geometric dyad in periodic anthyphairesis.

In the \emph{Sophist}, the dyad defining an intelligible Being receives a standard name: Being, Not-Being, as in the cases Beautiful, Not-Beautiful, or Just, Not-Just, and so on. In other words, in the dyad One, Being of the Second Hypothesis of the \emph{Parmenides}, the One corresponds to the Being and the Being corresponds to the Not-Being.

For exactly the same reason that the One and the Being are equalized in the \emph{Parmenides}, the Being and the Not-Being are \textbf{equalized} in the \emph{Sophist}, 
and the Not-Being is as good an intelligible Being as the Being.

Τhe paradigmatical intelligible Being One in the Second Hypothesis in the \emph{Parmenides}, possessing the part Being, 
forming the dyad One \& Being satisfying the philosophical analogue of periodic anthyphairesis has its exact analogue in the intelligible Being Being in the \emph{Sophist}, 
possessing the part not-Being forming the dyad Being \& not-Being satisfying the philosophical analogue of periodic anthyphairesis.

The theory of intelligible Being in the second hypothesis of the \emph{Parmenides} in terms of the dyad One \& Being is fully equivalent to the theory of the 
intelligible Beings in the \emph{Sophist} in terms of the dyad Being \& not-Being. The divisions in the \emph{Sophist} and \emph{Statesman} are periodic anthyphairetic, as in the second hypothesis of the \emph{Parmenides}.

A recent paper by S. Meister, 2025~\cite{Meister2025} provides a clear description of the difficulties faced by the Platonists in obtaining a satisfactory interpretation of the statement ``the not-Being is a Being", in the absence of the structure of an intelligible Being as the philosophical analogue of periodic anthyphairesis.

\section{The compresence of almost contradictory properties, Infinite \& Finite, One \& Many, In itself \& in the Other, Motion \& Rest, in Plato's intelligible Being; and the use of these properties to separate the intelligibles from the sensibles}

\subsection{Infinite and Finite 145a2--3}

The One of the second hypothesis is \textbf{Infinite}, because the anthyphairesis of the One to the Being is infinite, thus generating an infinite multitude of parts; and is \textbf{Finite}, by anthyphairetic periodicity.

\subsection{One and Many 145a2--3}

This has been discussed in Section 8.

\subsection{The One is in itself and in the Other}
\begin{quote}
τὸ ἓν ἀνάγκη αὐτό τε ἐν ἑαυτῷ εἶναι καὶ ἐν ἑτέρῳ.\\
— Ἀνάγκη. \emph{Parmenides} 145b6--e6 \\
and thus the One must be both in itself and in the Other."\\
``It must."
\end{quote}
\textit{Proof} (after \emph{Parmenides} 145b6--e6, 148d--149d, 138a3--7). By anthyphairetic periodicity there is an index $k$ such that $\text{One}/\text{Being}=\text{One}_k/\text{Being}_k$. Then

[1] the One is in itself,\\
since the $\text{One}_k$, which is the same as the One, is a part of the One; and,

[2] the One is in the Other\\
since the $\text{One}_k$, which is the same as the One, is a part of the Being, and is Other to the One.

\subsection{The One is at rest and in motion 145e7--146a8}
\begin{quote}
Οὕτω δὴ πεφυκὸς τὸ ἓν ἆρ' οὐκ ἀνάγκη καὶ κινεῖσθαι καὶ ἑστάναι; 

— Πῇ;

This being its nature, must not the one be both in motion and at rest?” 

“How is that?”
\end{quote}
\textit{Proof}
\begin{enumerate}
    \item [1)] — Ἕστηκε μέν που, εἴπερ αὐτὸ ἐν ἑαυτῷ ἐστιν·

ἐν γὰρ ἑνὶ ὂν καὶ ἐκ τούτου μὴ μεταβαῖνον ἐν τῷ αὐτῷ ἂν εἴη, ἐν ἑαυτῷ. 

— Ἔστι γάρ. 

“It is at rest, no doubt, if it is in itself; for being in one, and not passing
out from this, it is in the same, namely in itself.” 

“It is.”

    \item [2)] — Τὸ δέ γε ἐν τῷ αὐτῷ ἀεὶ ὂν ἑστὸς δήπου ἀνάγκη ἀεὶ εἶναι. 

— Πάνυ γε. 

“But that which is always in the same, must always be at rest.” 

“Certainly.” 

— Τί δέ; τὸ ἐν ἑτέρῳ ἀεὶ ὂν οὐ τὸ ἐναντίον ἀνάγκη μηδέποτ' ἐν ταὐτῷ εἶναι,

μηδέποτε δὲ ὂν ἐν τῷ αὐτῷ μηδὲ ἑστάναι, μὴ ἑστὸς δὲ κινεῖσθαι;

— Οὕτως. 

— Ἀνάγκη ἄρα τὸ ἕν, αὐτό τε ἐν ἑαυτῷ ἀεὶ ὂν καὶ ἐν ἑτέρῳ, ἀεὶ κινεῖσθαί τε καὶ
ἑστάναι. 

— Φαίνεται. 145e7--146a8

“Well, then, must not, on the contrary, that which is always in other be never
in the same, 

and being never in the same be not at rest, and being not at rest be in motion?”

“True.” 

“Then the one, being always in itself and in other, must always be in motion and
at rest.” 

“That is clear.”

\end{enumerate}

\subsection{The final and most seemingly contradictory final Proposition for the One of the second hypothesis (\emph{Parmenides} 151c3--155d3)}

\subsubsection{The final Proposition for the One of the second hypothesis of \emph{Parmenides} 151c3--155d3}
\begin{quote}
Ἆρ' οὖν καὶ χρόνου μετέχει τὸ ἕν, καὶ ἐστί τε καὶ γίγνεται

νεώτερόν τε καὶ πρεσβύτερον αὐτό τε ἑαυτοῦ καὶ τῶν ἄλλων,

καὶ οὔτε νεώτερον οὔτε πρεσβύτερον οὔτε ἑαυτοῦ οὔτε τῶν ἄλλων,

χρόνου μετέχον; Plato’s Parmenides 151e3--6  

“And does the one partake of time and if it partakes of time,

is it and does it become younger and older than itself and other things, 

and neither younger nor older than itself and the others?”~
\end{quote}

The proof of this proposition for the One of the second hypothesis of
the \emph{Parmenides} occupies the passage 151c3--155d3. The
proposition appears to contain several contradictions at the same
time; scholars of Plato have had difficulty in making sense of this passage.

\subsubsection{The interpretation of \textit{Parmenides} 151c3--155d3 by Meinwald, 1991~\cite{Meinwald1991}, 1992~\cite{Meinwald1992}, Peterson, 1996~\cite{Peterson1996}}

Meinwald, 1991~\cite{Meinwald1991}, 1992~\cite{Meinwald1992}, proposed an interpretation of the second hypothesis in the
Parmenides by distinguishing two types of predication, 

predication in relation to itself or tree predication, and 

predication in relation to others or ordinary predication: 
\begin{quotation}
‘In the \textit{Sophist}, \textit{Statesman}, and \textit{Philebus} Plato devotes a great deal of attention
to such trees’… 

‘the sort of genus-species tree familiar from the Linnaean classification
system’, such as ‘dividing Animal into Vertebrate and Invertebrate, dividing
Vertebrate in turn into Mammal and so on, and continuing with such divisions
through Feline and Cat, to produce at last such infimae species as Persian Cat.
..

In such a tree, a kind A appears 

either directly below or far below another kind B… 

In any such case, B can be truly predicated of A… 

(I will sometimes use the phrase ``tree predication" for this type.)….. 

Predication of a subject in relation to itself takes some explaining because we
ourselves are not in the habit of making such predications. 

Predication in relation to the others is much easier to understand, because this
is the category into which Plato would put our own common or garden
predications. (Thus I will sometimes call these ``ordinary" or ``everyday"
predications.).’ Meinwald, 1992~\cite{Meinwald1992} p.378--380:
\end{quotation}

Peterson, 1996~\cite{Peterson1996} p.178
\begin{quotation}
For example, 155b4--c4 establishes this apparent contradiction (discussed
on Meinwald, 1991~\cite{Meinwald1991} p.116): 

``the one does not become younger than the others (pros ta alia) and 

the one does become younger than the others (pros ta alla)." 

Naturally a full understanding of its truth requires understanding the proof
which establishes it, of which I forego discussion. The point here is that the
result is not a genuine contradiction. 

Plato makes clear that the result means that the one does not become younger in
one respect and does become younger in another respect. And Plato tells us the
respects (I paraphrase): 

``The one does not become younger than the others (pros ta alla) in that their
ages always differ from one another by the same amount; and 

the one does become younger than the others (pros ta alla) in that the ratio of
the ages of older to younger gets smaller over time." 

A parallel example with a subject different from the one would be that when a
child has his first birthday, a parent may be about thirty times the child's
age. Twenty years later, that parent will be only about twice as old as the
child. 

The parent becomes younger in relation to the child. But of course in a more
familiar way the parent does not become younger than the child.

Plato explicitly supplies the specific qualifications for this particular
apparent

contradiction about the one: he says 

``they always differ from each other by an equal number" and 

``[they] must continually differ from each other by a different portion." 

But often the qualifications needed to make sense of the results within a
section are not explicit. There is work left for the reader. p.178]]

The problem with these attempts to interpret the almost contradictory behavior of an intelligible Being is that they do not take into account its structure as a philosophical analogue of a dyad in periodic anthyphairesis.
\end{quotation}

\subsubsection{The interpretation of \emph{Parmenides} 151c3--155d3 by Proclus, in his \emph{Commentary to Plato's Parmenides} 1225,30--1227, 31 provides an essential and lucid explanation of this seemingly contradictory passage in terms of the cyclical, periodical nature of time and of motion}

\begin{quote}
Βούλεται μὲν κατασκευάσαι διὰ τούτων τε καὶ τῶν ἑξῆς ὅτι

πᾶν τὸ χρόνου μετέχον 

πρεσβύτερόν ἐστι καὶ ἰσήλικον ἑαυτῷ· 

τοῦτο δὲ βουλόμενος ἀναγκαίως προαποδείκνυσιν ὅτι 

αὐτό τι ἑαυτοῦ πρεσβύτερόν ἐστιν· 

ᾗ δὲ αὐτὸ ἑαυτοῦ πρεσβύτερον, δῆλον ὅτι καὶ ἑαυτοῦ νεώτερόν ἐστιν·

πρὸς γὰρ ἑαυτὸ καλεῖται νεώτερον ἅμα καὶ πρεσβύτερον. 

His purpose in this and the following passage is to establish that 

everything that partakes in time is older than, and of equal age with, itself: 

and since he wishes this, he necessarily begins by demonstrating that 

something is older than itself; 

and in so far as it is older than itself, it is plain that 

 it is also younger than itself; 

for it is in relation to itself that it is called both younger and older. 1225,30--37
\end{quote}

\begin{quote}
Δόξειε δ' ἂν ἀπορώτατος εἶναι, καὶ, ἵν' εἴπω, 

σοφιστικός πως οὗτος ὁ λόγος· 

πῶς γὰρ αὐτό τι ἑαυτοῦ καὶ πρεσβύτερον ἅμα καὶ νεώτερόν ἐστιν; 

Now this argument might seem to be problematic in the extreme, 

one might even say sophistic; 

for how could something be simultaneously older and younger than itself? 
\end{quote}

\begin{quote}
Οὔκουν δὴ ὅ γε Σωκράτης, \\
ὁ ἑαυτοῦ πρεσβύτερος γεγονὼς, καὶ νεώτερος  αὑτοῦ ἐστι·\\
τὸ μέν γε πρεσβυτικὸν αὐτῷ πάρεστιν, οἴχεται δὲ ἡ νεότης. \\
Surely the Socrates who has become older than himself \\
is not also younger than himself; \\
at any rate being older is present to him, while being young is gone from him.\\
1225,37--1226,

\medskip

[1] Ἔνιοι μὲν οὖν ἐνέδωσαν τῷ λόγῳ \\
καὶ οὐδὲ τοῦτο ἀπέσχοντο εἰπεῖν, τι σοφιζομένῳ ἔοικεν ὁ Πλάτων, \\
πάνυ προχείρως τὴν ἑαυτῶν ἄγνοιαν ἐπὶ τὸν λόγον μετάγοντες· \\
So some commentators have given up in face of this argument, and\\
have not scrupled even to say that Plato is appearing here to indulge in sophistry,\\
all too readily transferring their own ignorance to the argument. 1226,2--6

\medskip

[2] οἱ δὲ καὶ ἀντισχεῖν πρὸς αὐτὸν πειραθέντες \\
μαλακώτερον προέστησαν τῆς ἀληθείας· \\
τὸ γὰρ αὐτὸ, φασὶν, ἅμα νεώτερον καὶ πρεσβύτερον·\\
κατὰ μὲν τὸν μέλλοντα χρόνον νεώτερον, οὔπω γὰρ ἐκεῖνον προσείληφε· \\
κατὰ δὲ τὸν παρελθόντα πρεσβύτερον, ὃν ἤδη βεβίωκεν. 

Ἀλλὰ τοῦτο \\
οὐκ ἦν ἑαυτοῦ νεώτερον ἅμα γίγνεσθαι καὶ πρεσβύτερον, \\
ἀλλὰ νεώτερον μὲν ἄλλου, πρεσβύτερον δὲ ἄλλου· \\
ταχὺ ἄρα ἡμῖν ἀσθενὴς ἐφάνη ὁ τούτου λόγος. \\
Others, again, in trying to stand up to these critics, \\
have championed the truth rather too weakly \\
for they say the same thing is at the same time younger and older; \\
in respect of future time it is younger, for it has not yet attained to that; \\
in respect of past time, on the other hand, it is older, since it has already\\
lived through that.

But this is \\
not what it is to become both simultaneously younger and older than oneself, \\
but rather younger than one thing and older than another; \\
so that the argument of this commentator is quickly revealed to us as being weak.
1226,6--15 

\medskip

[3] Οἱ δὲ ἀνάπαλιν ἔφασαν \\
τὸ πρεσβύτερον ἑαυτοῦ πᾶν εἶναι καὶ νεώτερον,\\
πρεσβύτερον μὲν τὸ νῦν ὂν, νεώτερον δὲ ὃ ἦν ἔμπροσθεν,\\
τὸ γοῦν νῦν πρεσβύτερον τοῦ πρότερον ὄντος νεωτέρου λέγεσθαι πρεσβύτερον· \\
οὐδ' οὗτοι συνέντες ὅπως φησὶν ὁ Πλάτων, \\
τὸ μὲν ὂν πρὸς τὸ ὂν λέγεσθαι, τὸ δὲ ἐσόμενον πρὸς τὸ ἐσόμενον,\\
ἀεὶ τῶν πρὸς τὸ τόδε· \\
γεγονὸς οὖν τοῦ γεγονότος νεωτέρου οὐ δυνατὸν φάναι τὸ \\
νῦν ὂν πρεσβύτερον· \\
τοῦτο γάρ ἐστιν ἐπαλλάττειν τοὺς χρόνους \\
καὶ μὴ φυλάττειν ὃν αὐτὸς εἶπε κανόνα περὶ τῶν πρός τι πάντων. 

Another set, again, declared that everything is both older and younger than itself,\\
what is now existent being older, and \\
what was before younger, \\
and that \\
what is now older can be said to be older than what was formerly younger; \\
but these too fail to understand the sense of Plato's statement (141b), \\
which puts what is against what is, what will be against what will be, and what\\
has been against what has been, these in all cases being relative expressions;\\
so then, it is not possible to say that \\
that which now is older has become older than that which had become younger; \\
for this is to mix up times and not to preserve the rule which he himself laid\\
down regarding all relative expressions. 1226,15--2

\medskip

%%%
[4] Πάλιν οὖν ἡμῖν ἐπὶ τὸν ἡμέτερον καθηγεμόνα τρεπτέον, \\
καὶ τὴν ἐκείνου παράδοσιν εἰς μέσον ἀκτέον, \\
φῶς ἀνάπτουσαν εἰς πάντα τὸν προκείμενον λόγον.\\
Διττὸν δὴ τὸ χρόνου μετέχον ἐστί\\
τὸ μὲν οἷον κατ’ εὐθεῖαν ὁδεῦον, καὶ ἀρχόμενόν τε ἀπό τινος καὶ εἰς ἄλλο καταλῆγον·\\
τὸ δὲ κατὰ κύκλον περιπορευόμενον, καὶ ἀπὸ τοῦ αὐτοῦ πρὸς τὸ αὐτὸ τὴν κίνησιν ἔχον, \\
ᾧ καὶ ἀρχὴ καὶ πέρας ἐστὶ ταὐτὸν καὶ ἡ κίνησις ἀκατάληκτος,\\
ἑκάστου τῶν ἐν αὐτῇ καὶ ἀρχῆς καὶ πέρατος ὄντος, καὶ οὐδὲν ἧττον ἀρχῆς καὶ πέρατος. \\
Τὸ δὴ κυκλικῶς ἐνεργοῦν μετέχει τοῦ χρόνου περιοδικῶς, \\
καὶ ἐπειδὴ τὸ αὐτὸ καὶ πέρας αὐτῷ τῆς κινήσεώς ἐστι καὶ ἀρχὴ,\\
καθόσον μὲν ἀφίσταται τῆς ἀρχῆς, πρεσβύτερον γίγνεται,\\
καθόσον δὲ ἐπὶ τὸ πέρας ἀφικνεῖται, νεώτερον γίγνεται·\\
γιγνόμενον γὰρ ἔγγιον τοῦ πέρατος ἐγγύτερον γίγνεται τῆς οἰκείας ἀρχῆς· \\
τὸ δὲ τῆς οἰκείας ἀρχῆς ἐγγυτέρω γιγνόμενον νεώτερον γίγνεται·\\
τὸ ἄρα ἐπὶ τὸ πέρας ἀφικνούμενον κυκλικῶς \\
νεώτερον γίγνεται τὸ αὐτὸ κατὰ τὸ αὐτὸ καὶ πρεσβύτερον γιγνόμενον·\\
τὸ γὰρ τῷ ἑαυτοῦ πέρατι συνεγγίζον ἐπὶ τὸ πρεσβύτερον πρόεισιν.\\
Ω μὲν οὖν ἀρχὴ ἄλλο καὶ τὸ πέρας, \\
τούτῳ καὶ τὸ νεώτερον ἕτερον ἢ τὸ πρεσβύτερον·\\
ᾧ δὲ ταὐτὸν ἀρχὴ καὶ πέρας, \\
οὐδὲν μᾶλλον νεώτερόν ἐστι τὸ νεώτερον ἢ πρεσβύτερον,\\
ἀλλ’ ὡς ὁ Πλάτων φησὶν, ἅμα νεώτερον ἑαυτοῦ καὶ πρεσβύτερον γίγνεται.\\
Τὸ μὲν γὰρ κατ’ εὐθεῖαν κινούμενον οὐκ ἔχει, \\
διότι πέρας καὶ ἀρχὴ διαφέρετον ἐπ’ αὐτοῦ,\\
τὸ δὲ κατὰ κύκλον ἔχει τὸ νεώτερον ἑαυτοῦ\\
τῷ πρὸς τὸ αὐτὸ καὶ ὡς ἀρχὴν καὶ ὡς πέρας γίγνεσθαι τὴν κίνησιν.\\
Πᾶν ἄρα τὸ χρόνου μετέχον, \\
εἰ πρεσβύτερον ἑαυτοῦ γίγνεται, καὶ νεώτερον ἑαυτοῦ γίγνεται·\\
τοιοῦτον δὲ τὸ κυκλικῶς κινούμενον.\\
Ὅθεν εἰκότως ἐθορυβοῦντο καὶ οἱ ἀρχαῖοι\\
μή πη σοφισματώδης οὗτος ὁ λόγος ἐστὶν,\\
εἰς τὰ κατ’ εὐθεῖαν κινούμενα βλέποντες,\\
δέον διελέσθαι καὶ θεωρῆσαι \\
τίσι μὲν τὸ αὐτὸ καὶ ἀρχὴ καὶ πέρας,\\
τίσι δὲ ἕτερον,\\
καὶ ὅτι νῦν περὶ τῶν θείων ψυχῶν ὁ λόγος,\\
αἳ καὶ χρόνου μετέχουσι περιοδικὸν ἔχουσαι\\
τὸν χρόνον τῆς θείας κινήσεως, ὥσπερ καὶ τὰ ὀχήματα τὰ ἐξημμένα αὐτῶν.\\
Ἀλλ’ οὗτος μὲν ὁ λόγος τοῦ καθηγεμόνος ἡμῶν.

\medskip

We must therefore turn once again to our Master, \\
and bring to bear upon the problem his discussion,\\
which illuminates the whole preceding argument.\\
That which partakes in time is of two sorts;\\
the one which, as it were, proceeds in a straight line,\\
beginning from one point and ending in another;\\
the other which travels round in a circle,\\
and pursues its motion from the same point to the same point,\\
so as to have a beginning and an end which are the same\\
and a motion which is unceasing,\\
since each point of its progress is both beginning and end,\\
and is no less a beginning than an end.\\
That, then, which enjoys cyclical activity partakes in time by circuits,\\
and since the same point is for it both an end and a beginning of motion;\\
in so far as it departs from a beginning, it becomes older,\\
whereas in so far as it arrives at an end, it becomes younger;\\
for as it comes to be nearer to its end,\\
it comes to be nearer to its own beginning;\\
and that which comes to be nearer to its own beginning becomes younger;\\
so then, that which arrives cyclically at its end becomes younger,\\
while at the same time and by the same process\\
also becoming older; \\
for that which draws near to its own end proceeds towards being older;\\
so for that which has a beginning different from its end,\\
becoming younger is different from becoming older;\\
but for that of which its beginning is the same as its end,\\
its youngness is no younger than it is older,\\
but as Plato says, \\
“It becomes simultaneously younger than itself and older.”\\
For that which is moved in a straight line \\
does not have this characteristic,\\
because end and beginning are different in its case,\\
whereas that which moves in a circle \\
has the quality of being younger than itself \\
by reason of the fact that its motion comes about\\
in relation to the same point as both beginning and as end.\\
So then, everything that partakes in time,\\
if it becomes older than itself, also becomes younger than itself;\\
and what is of this sort is that which moves in a circle.\\
For this reason the ancients\\
were understandably disturbed \\
lest this argument be in some way sophistical,\\
since they were looking to things\\
that moved in a straight line,\\
whereas they should have made the distinction\\
and considered what things have the characteristic\\
of having their beginning and end the same,\\
and what have them as different,\\
and they should have considered that\\
the subject of discussion now is divine souls,\\
which partake in time in the respect\\
that they have a periodic time of their proper motion,\\
as indeed do the vehicles which are dependent upon them.\\
This, then, is the argument of our Master. 1226,26–1227,31\\
%%%

\medskip

[[Translation by G.M.~Morrow and J.M.~Dillon 1987 \cite{Proclus1987}]]
\end{quote}

\subsubsection{Our interpretation of the \textit{Parmenides} 151c3--155d3 completes Proclus’ interpretation and rejects Meinwald’s interpretation}

Our interpretation of the \textit{Parmenides} agrees and completes the one given by Proclus. The commentary by Proclus clearly builds its interpretation in terms of the periodicity of motion, and we complete this interpretation by suggesting that the periodic motion is generated by the periodic anthyphairesis of the intelligible Being. In fact, as we have seen above, the One of the second hypothesis is in motion when it is in the other, and at rest when it is in itself, and we have interpreted the statesment to be “in itself” and “in the other” in terms of periodic anthyphairesis; hence the difficult and seemingly contradictory passage is interpreted by means of periodic anthyphairesis.

According to our interpretation the tree-like division in the divisions in the \textit{Sophist} and the \textit{Statesman} is an abbreviated form of the anthyphairetic Division in the \textit{Parmenides} 142d9--143a3 (One participates in Being), in Section 3. This Division is not of the form of Linnean classification and does not end with an infimae species, as with Aristotle’s and as Meinwald believes, but is infinite. Participation, related to the infinity of the anthyphairetic division, is not sufficient to explain the statements in the second hypothesis in the \textit{Parmenides}, but must be complemented with the presence of the One in Being (\textit{Parmenides} 144c2--d4) and the equalization of the One with the Being, via periodicity of the anthyphairesis.
Meinwald’s approach is not related to (anthyphairetic) periodicity.

\subsection{The use of the Compresence properties for separating the Intelligible Beings from the sensibles}
The almost contradictory compresence of these properties of an intelligible Being are useful for separating the intelligible Beings from the sensibles. This separation originated by Zeno, as described in Plato’s \textit{Parmenides} 127e–128a, analyzed in detail in Negrepontis, 2023~\cite{Negrepontis2023}, whose general argument, has the following structure:
\begin{enumerate}
    \item[Step 1] Suppose that the sensible is an intelligible Beng.
    \item[Step 2] An intelligible Being satisfies one of the Compresence properties
    \item[Step 3] Then a sensible entity must satisfy this Compresence property.
    \item[Step 4] But from our experience, we know that this is impossible, a contradiction.
    \item[Step 5] Hence the intelligible Beings are different from the sensible entities.
\end{enumerate}
This argument was adopted by Plato, as seen in the \textit{Sophist} 250b, where the intelligible Beings are separated from the sensibles because an intelligible Being is simultaneously in motion and at rest, something “most impossible” for the sensibles.

\section{The Third Man Argument in the Parmenides 132a1--b2, Vlastos’ interpretation, and its resolution: Vlastos’ Non-Identity principle in the Third Man Argument is false, because the self-similar Oneness of an intelligible Being prevents the infinite regress}
The complete description of the One of the second hypothesis, obtained in terms of periodic anthyphairesis, is powerful, and so naturally there are several consequences, some of which we have already seen. A further consequence is the resolution of the Third Man Argument: Vlastos’ reconstruction of the Third Man Argument fails, exactly because the non-Identity principle does not hold, on account of the self-similar Oneness of an intelligible Being.

\subsection{The Third Man Argument in the \textit{Parmenides} 132a1--b2}
\begin{quotation}
Οἶμαί σε ἐκ τοῦ τοιοῦδε ἓν ἕκαστον εἶδος οἴεσθαι εἶναι·\\
ὅταν πόλλ' ἄττα μεγάλα σοι δόξῃ εἶναι, \\
μία τις ἴσως δοκεῖ ἰδέα ἡ αὐτὴ εἶναι ἐπὶ πάντα ἰδόντι, \\
ὅθεν ἓν τὸ μέγα ἡγῇ εἶναι.

Ἀληθῆ λέγεις, φάναι.\\

Τί δ' αὐτὸ τὸ μέγα καὶ τἆλλα τὰ μεγάλα, \\
ἐὰν ὡσαύτως τῇ ψυχῇ ἐπὶ πάντα ἴδῃς, \\
οὐχὶ ἕν τι αὖ μέγα φανεῖται, ᾧ ταῦτα πάντα μεγάλα φαίνεσθαι;

Ἔοικεν.\\

Ἄλλο ἄρα εἶδος μεγέθους ἀναφανήσεται, \\
παρ' αὐτό τε τὸ μέγεθος γεγονὸς καὶ τὰ μετέχοντα αὐτοῦ· \\
καὶ ἐπὶ τούτοις αὖ πᾶσιν ἕτερον, ᾧ ταῦτα πάντα μεγάλα ἔσται· \\
καὶ οὐκέτι δὴ ἓν ἕκαστόν σοι τῶν εἰδῶν ἔσται, \\
ἀλλὰ ἄπειρα τὸ πλῆθος.
\end{quotation}
\begin{quotation}
 “I fancy your reason for believing that each idea is one is something like this; \\
when there are many things (polla) which seem to you to be great (megala) [[$S$]],\\
you may think, as you look at them all, that there is one and the same idea\\
(mia… idea he aute) [[$\Sigma_1$]] over all of them (epi panta), and hence you think\\
the great is one (hen to mega).\\

“That is true,” he said.\\

“But if with your mind's eye you regard the great itself (auto to mega) and the\\
other great things (talla ta megala) [[$S$, $\Sigma_1$]] in the same way, will not\\
another great (hen ti au mega) [[$\Sigma_2$]] appear beyond, by reason of which (hoi)\\
all these (tauta panta) must appear to be great (megala)?”\\

“So it seems.”\\

“That is, another idea of greatness (allo eidos megethous) [[$\Sigma_2$]]  will appear,\\
in addition to (par’) greatness itself (auto to megethos) and the objects which\\
participate in it (ta metechonta autou) [[$S$, $\Sigma_1$]];\\

and another (heteron) [[L3]] again on all these (epi toutois… pasin) [[$S$, $\Sigma_1$, $\Sigma_2$]], \\

by reason of which (hoi) all these (tauta panta) are great (megala); \\

and each of your ideas will no longer be one, but will be infinite in multitude\\
(apeira to plethos) [[$\Sigma_1,\Sigma_2,\Sigma_3,\dots,\Sigma_n,\dots $]].”
\end{quotation}

[[Plato. Plato in Twelve Volumes, Vol. 9 translated by Harold N. Fowler. Cambridge, MA, Harvard University Press; London, William Heinemann Ltd.\ 1925 \cite{Plato} With modifications by the author; symbols added by the author]]

\subsection{Vlastos’reconstruction of the Third Man Argument}
Gregory Vlastos reconstructs the Platonic Third Man Argument (TMA) in the \textit{Parmenides} 132a1--b2 by introducing three hypotheses:

\textit{One Over Many (OOM):} For any collection of entities possessing property F, there is an intelligible Being of F-ness, such that every element of the collection has property F by participating in an intelligible Being of F-ness.

\textit{Self-Predication (SP):} Every intelligible Being of F-ness is itself an entity possessing property F.

\textit{Non-Identity (NI):} If an intelligible Being participates in an intelligible Being then it is not identical to it. No form participates of itself.

\subsection{Rickless’ description of Vlastos’ reconstruction of the Third Man Argument}
Rickless’ description, 2020~\cite{Rickless2020}, based on Vlastos’ reonstruction
of the Third Man’s Argument (Vlastos, 1954~\cite{Vlastos1954}, 1955~\cite{Vlastos1955},1956~\cite{Vlastos1956},
1969~\cite{Vlastos1969}), is the following:

\begin{quotation}
‘Parmenides generates the infinite regress as follows.
\begin{enumerate}
    \item Consider a plurality of large things, $A$, $B$, and $C$. By One-over-Many, there is a form of largeness (call it “$L_1$”) by virtue of partaking [participates] of which $A$, $B$, and $C$ are large. By Self-Predication, $L_1$ is large.
    \item So there is now a new plurality of large things, $A$, $B$, $C$, and $L_1$. Thus, by One-over-Many, there is a form of largeness (call it “$L_2$”), by virtue of partaking[participating] of which, $A$, $B$, $C$, and $L_1$ are large. Hence $L_1$ partakes of $L_2$. At this point, Parmenides assumes something like the following Non-Identity assumption:
    (Non-Identity) No form is identical to anything that partakes [participates] of it.
    
    \dots

    From the fact that $L_1$ partakes [participates] of $L_2$, Non-Identity entails that $L_2$ is numerically distinct from $L_1$. Thus, there must be at least two forms of largeness, $L_1$ and $L_2$. But this is not all. By Self-Predication, $L_2$ is large.
    \item So there is now a new plurality of large things, $A$, $B$, $C$, $L_1$, and $L_2$. Thus, by One-over-Many, there is a form of largeness (call it “$L_3$”) by virtue of partaking [participating] of which $A$, $B$, $C$, $L_1$, and $L_2$ are large. Hence $L_1$ and $L_2$ both partake [participate] of $L_3$. But then, by Non-Identity, $L_3$ is numerically distinct from both $L_1$ and $L_2$.
    \item Thus, there must be at least three forms of largeness, $L_1$, $L_2$, and $L_3$.’ [Rickless]
\end{enumerate}
\end{quotation}

\subsection{Some arguments in justification of the Non-identity principle and some evidence that the (strong) non-Identity principle is disturbing}

\subsubsection{ Although Plato does not explicitly state the non-Identity principle, nevertheless Vlastos, 1954~\cite{Vlastos1954} p.325, regards it as obvious:}
    \begin{quotation}
    ‘In the many modern discussions of the Argument I can find no explicit statement that this Nonidentity Assumption, or an equivalent one, is strictly required in just this way. This may be because the role of such an assumption at this point strikes critics more nimble-witted than myself as so obvious that they feel it an insult to their reader's intelligence to put it into words or symbols. However, there are times when the drudgery of saying the obvious is rewarded, and this is one of them.’
    \end{quotation}

Th. Scaltsas, 1992~\cite{Scaltsas1992} mentions some arguments in favor of accepting the non-Identity principle:

    \begin{quotation}
    [p.221] ‘We can find justification for the Non-Identity thesis in Aristotle's system, in which, as in Plato's Theory, the generation of things requires the existence of forms. Form is not generated, according to Aristotle.  If it were, the generation of anything would require an infinite regress of generations. So, if not generated, form must be passed down from cause to product. Hence, the cause must be different from the product (NI), and must possess the form it passes down to the product (SP).’
    \end{quotation}

    \begin{quotation}
    [p.225] ‘Understanding the Theory in terms of the biological model, a Form F is the source of F-ness, and things that participate in the Form inherit their F-ness from the Form. The source from which the f-ness is inherited must be different from the thing that inherits it. Otherwise, no explanation is offered since a thing can inherit nothing from itself. Hence, the Non-Identity between an entity and the Form it partakes of.’
    \end{quotation}

    \subsubsection{The weak and strong non-identity}
    On the other hand, an indication that the strong non-identity principle is disturbing is the early differentiation of the strong from the weak principle
\begin{itemize}
    \item Fine, 1993 \cite{Fine1993}
    \begin{quotation}
    As Vlastos noted long ago, it is important to distinguish between weak and strong non-identity. According to weak non-identity, sensibles are F not in virtue of themselves (they are F by being suitably related to a form of F). According to strong non-identity, nothing is F in virtue of itself; not even a form of F can be F in virtue of itself. Strong non-identity implies weak non-identity, but not conversely. NI expresses strong nonidentity, and P-TMA requires strong non-identity. G. Fine,1993 \cite{Fine1993} p.207
    \end{quotation}
    \begin{quotation}
    NI in the Accurate One over Many Argument says that what is predicated of a plurality of things is different from the things of which it is predicated and, as we have seen, in the Accurate One over Many Argument predications extend beyond sensibles. NI is therefore the claim that nothing is F by being predicated of, i.e. in virtue of, itself. A[lexander]-TMA thus involves the same non-identity assumption as P-TMA and the Resemblance Regress. p.216
    \end{quotation}
    \begin{quotation}
    In the last chapter I distinguished between weak and strong non-identity. According to weak non-identity, sensible Fs are F in virtue of something distinct from themselves; according to strong non-identity (= NI), nothing is F in virtue of itself. Now, in some phases of his career, Plato seems to believe that there is a form for every property-name; so he is sometimes committed to weak non-identity. For if there is a form of F, then sensible Fs are F by being suitably related to it, and so in virtue of something distinct from themselves. Plato could, however, accept weak non-identity without being committed to strong non-identity. Yet strong non-identity—NI—is the non-identity assumption needed for the regress. p.225
    \end{quotation}
\end{itemize}

\subsection{Some interpretations of the Third Man Argument}
Various modifications of Vlastos’ original reconstruction have been suggested. A common feature of these modifications is that, with the exception of Rickless, explicitly or implicitly assume the Non-identity principle Here is a sample of them.

W. Sellars, 1955~\cite{Sellars1955} argued that non-identity is not consistent with self-predication, as formulated by Vlastos, and proposed a modification, in particular formulating NI as follows:
    \[\text{(NI) If x is F, then x is not identical with the F-ness by virtue of which it is F. (Sellars, 1955~\cite{Sellars1955} p.418)}\]

S. M. Cohen, 1971 \cite{Cohen1971} has proposed a modification of the argument so that self-predication and non-identity assumptions are built in the hypotheses, but not made explicit.

H. Teloh and D.J. Louzecky, 1972~\cite{Teloh1972}, and H. Teloh, 1982~\cite{Teloh1982} p.158--167, argue that an infinite regress follows from a single premiss that does not require the notion of Self-Predication, of course keeps the principle of non-identity. All that is needed is
    \begin{quotation}
    (T): `If a number of things are F, there is a single Form in virtue of which we apprehend these things as F, and these things . . . are not identical with this Form', which requires only the notions of predication and non-identity (Teloh and Louzecky, 1972~\cite{Teloh1972} p. 87).
    \end{quotation}
    In Meinwald, 1992~\cite{Meinwald1992} she applied her interpretation of
    Plato’s Parmenides, to the Third Man argument. Frances, 1996~\cite{Frances1996}
    describes Meinwald’s rejection of the Third Man Argument, finding
    fault with the SP hypothesis, on the basis of her distinction of X
    is F pros heauto and pros alla.

\begin{itemize}
   \item  Frances, 1996~\cite{Frances1996}
    \begin{quotation}
    According to Meinwald the problem lies in the self-predicative claim that Fness is F. Montblanc and Venus may be large, but the Form the Large certainly is not large, at least not in the same way Montblanc and Venus are large. So
    either ‘the Large is large’ is obviously false and Plato has made a rather gross error,
    or the Form the Large is being predicated of the Large in a way different from the way it is predicated of Montblanc and Venus.
    In support of this latter option Meinwald’s two forms of predication play the essential role. Recall that to say a is F pros heauto is to say that the F is at least part of the nature of a; the F is definitionally true of a.
    Meinwald thinks it is pretty clear that the Large is large pros heauto:
    it is true that Largeness is definitionally true of the Form the Large.
    On the other hand, to say that a is F pros ta alla is to say that a displays some feature that conforms to the nature of the F. Since Montblanc, Venus, and other ordinary large objects are large pros ta alla but not large pros heauto, and the two kinds of predication are different, we conclude that largeness is being predicated of the Large in a way different from the way it is being predicated of Montblanc and Venus.

    In light of these different forms of predication Meinwald would formulate the true intuitions behind the third man argument as follows (cf. 1991 \cite{Meinwald1991} 155--157):
    \begin{enumerate}
        \item[1´.] a and b are F pros ta alla.
        \item[2´.] If there is a plurality P each of which is X pros ta alla,
        then there is exactly one Form that
        \begin{enumerate}
            \item[(a)] is distinct from each of P,
            \item[(b)] is such that each of P is X pros ta alla in virtue of it, and
            \item[(c)] is X pros heauto.
        \end{enumerate}
        \item[3´.] By (1´) and (2´) we infer that
        there is exactly one Form—call it Fness—that
        \begin{enumerate}
            \item[(a)] is distinct from a and b,
            \item[(b)] is such that a and b are F pros ta alla in virtue of it, and
            \item[(c)] is F pros heauto.
        \end{enumerate}
    \end{enumerate}
    She claims that the proper formulations of the intuitions that drive the argument fail to lead to regress or inconsistency.
    As before, (1´) and (3´c) entail that there is a plurality Q consisting of a, b, and Fness
    each of which is F—but this time not all of the members of Q are F pros ta alla.
    In particular, Fness is not F pros ta alla—it is F pros heauto.
    For example, the Large is not large and Redness is not red pros ta alla; Forms are neither colored nor have size.
    So we cannot reuse (2´) for plurality Q to derive the regress or inconsistency.
    What is essential to Meinwald’s solution to the third man argument is that
    \begin{enumerate}
        \item[(a)] the members of plurality P, i.e., the participants of Fness, are F in a way different from the way Fness is F,
        \item[(b)] the members of P all are F pros ta alla, and
        \item[(c)] the crucial claim $\lceil\text{Fness is F}\rceil$ is a pros heauto claim.
    \end{enumerate}
    Meinwald claims that this is the end of the problem of the third man:
    ‘Plato’s metaphysics can say goodbye to the Third Man’ (1991 \cite{Meinwald1991} 157).
    Surely Meinwald is correct in saying that in order for ‘The Large is large’ to present any problems it must be the same kind of claim as ‘Montblanc is large’, i.e., a pros ta alla claim. It is clear that the Large is large pros heauto and is not large pros ta alla.
    And since, on Meinwald’s view, claims of the form ‘the F is F pros heauto’ are truisms,
    it is equally obvious that the Large is large pros heauto.
    The same holds for the Forms Justice, the Cat, etc.

    Meinwald’s solution appears to succeed in that it seems to block the third man argument for some forms, and it is reasonable to hold that it was the solution Plato was indicating in the Parmenides, for the following reason.
    Consider first that it appears that part of Plato’s objective in the long second part of the dialogue was to introduce and illustrate the pros ta alla/pros heauto predication distinction.
    Furthermore, according to Plato this part of the dialogue was supposed to provide the intellectual exercise necessary in order for one to avoid the traps—such as that of the third man—that Parmenides revealed in the first part of the dialogue and that are associated with an immature theory of Forms and participation.
    Thus, it is reasonable to conclude that the predication distinction is both the central element of the exercise and the heart of the solution of the third man problem. By arguing that Plato had this predication distinction in mind for this purpose Meinwald has certainly offered an appealing answer to the age-old question of how Plato thought the problem of the third man was to be solved by the paradoxical dialectic.
\end{quotation}

    \item Petersom, 1996 \cite{Peterson1996} p.170
    \begin{quotation}
      Meinwald's insistence that the qualifying phrases are to be supplied
      systematically---that the results of the negative sections need supplement by pros heauto and
      the results of the positive sections need supplement by pros ta alla---leads
      to a new understanding of apparent contradictions
    between results deduced in different sections.
    For example,
    the first of the eight sections [First Hypothesis],
    from the hypothesis that the one is, gets the result that the one is not many,
    while the second section [Second Hypothesis],
    from the same hypothesis, gets the result that the one is many
    If we supply pros heauto to the negative result from the first section,
    and supply pros ta alla to the positive result from the second section,
    the combined results are:
    The one is not many pros heauto and the one is many pros ta alla
    which does not have the form of a contradiction.
    \end{quotation}
    \item Variants of her interpretation of the Third Man argument are in F.~J.~Pelletier, and E.~N.~Zalta, 2000~\cite{Pelletier2000}.
    \item Rickless, 2020~\cite{Rickless2020}
    \begin{quotation}
    There is another way of answering the three central interpretive questions,
    one on which Parmenides’ criticisms as well as the Deductions come out as serious and valid. (This is the interpretation defended in Rickless, 2007~\cite{Rickless2007}, and one aspect of which is defended, though on different grounds, in Gill, 2014~\cite{Gill2014}.)
    What Parmenides’ criticisms reveal is that, whether combined with the Pie Model conception of partaking or with Paradigmatism,
    Plato’s middle period theory of forms is internally inconsistent.
    It turns out that there are three principles the abandonment of which would eliminate all inconsistencies apart from the Greatest Difficulty:
    Careful logical analysis of the second part of the dialogue then reveals that the Deductions establish
    not only that the forms posited by the middle period theory exist,
    but also that Purity-F, Uniqueness, and No Causation by Contraries
    are all false.
    It is then reasonable to suppose that Plato meant the reader to recognize that the proper way to save the forms is by abandoning these three basic assumptions.
    And, importantly, this can be done without abandoning the most important principles at the heart of the middle period theory, namely One-over-Many and Separation.
    \end{quotation}
    \begin{quotation}
    [[Definitions of concepts in Rickless’ interpretation, cf Rickles, 2007~\cite{Rickless2007}, Abbreviations:
    --P (Purity) For any property F that admits a contrary (con-F), the F is not con-F.
    --NCC (No Causation by Contraries) For any property F that admits a contrary (con-F),
    whatever makes something be (or become) F cannot itself be con-F.
    --U (Uniqueness) For any property F, there is exactly one form of F-ness.
    --OM (One-over-Many) For any property F and any plurality of F things, there is a form of F-ness by virtue of partaking of which each member of the plurality is F.
    --SP (Self-Predication) For any property F, the F is F. xiii
    --NSP (Non-Self-Partaking) No form partakes of itself. xii, xiii]]
    \end{quotation}
\end{itemize}

\subsection{The falsity of the Non-Identity principle and of the Third Man Argument}
From our analysis of the One of the second hypothesis in the \textit{Parmenides}, given in sections 5,6, and 8, above, we conclude that
\begin{quotation}
Vlastos’ reconstruction of the Third Man’s Argument is faulty, in that he assumes the Non-Identity principle.
\end{quotation}

The intelligible Being Large, denoted by $\Sigma_1$, has the same structure as the paradigmatical intelligible Being, the One of the second hypothesis, hence for convenience we identify $\Sigma_1$ with the One of the second hypothesis, and we denote by $S$ the family of all sensible entities, the large things, that participate in $\Sigma_1$.

\noindent \textit{Step 1 of the Vlastos reconstruction of TMA works well:}

\begin{enumerate}
    \item[Step 1.1.] Indeed, every member of $S$ is large by virtue of participating in the intelligible Being, $\Sigma_1$. Thus, the One Over Many (OOM) hypothesis is satisfied.
    \item[Step 1.2.] We note that
    \begin{enumerate}
        \item[(i)] Every sensible entity participates in an intelligible Being \textit{Parmenides} 130b1--131a3
        
A necessary step for answering the Question set in 134c--d is the unqualified statement in 130b1--131a3, 
\begin{quotation}
every sensible entity participates in an intelligible Being.
\end{quotation}

    \item[(ii)] A true Opinion of an intelligible Being is a finite initial segment of the full infinite anthyphairesis of the intelligible Being.

        One good way to learn this is from the \textit{Meno}, as explained in Negrepontis 2024~\cite{Negrepontis2024a}. In the \textit{Meno} the intelligible Being is the diameter in relation to the side of a square, whose sequence of quotients is the infinite sequence $[1,2,2,2,\dots]$. A true Opinion of the diameter to the side is given by a finite initial segment $[1,2,2,\dots,2 \text{(n times)}]$, describing the anthyphairesis of the $n$-th dyad of side and diameter numbers.
        
        \item[(iii)] The knowledge of a sensible entity participating in an intelligible Being is a true Opinion of the intelligible Being

The best Platonic passage establishing that
\begin{quotation}
The knowledge of a sensible entity participating in an intelligible Being\\
is a true Opinion of the intelligible Being
\end{quotation}        

is \textit{Republic} 475d--480a.
    \end{enumerate}
   
    Thus,
    \begin{itemize}
        \item a sensible entity that is large participates in the intelligible Being Large $\Sigma_1$, and is fully known is known by a finite segment of infinite anthyphairesis of the dyad $\langle\text{One, Being}\rangle$; and,
        \item the intelligible Being $\Sigma_1=\text{One}$ is clearly not identical to any member $s$ of the family $S$; and the Non-identity (NI) hypothesis indeed holds for Step 1 of the TMA.
    \end{itemize}
    \item[Step 1.3.] By the hypothesis of Self-Predication (SP), the intelligible Being $\Sigma_1=\text{One}$ is itself large. Thus Step 1 of Vlastos’ reconstruction of TMA indeed holds true.
\end{enumerate}

\noindent \textit{However, Step 2 of Vlastos’ reconstruction breaks down.}

\begin{enumerate}
    \item[Step 2.1.] We are to find a Platonic Idea $\Sigma_2$, satisfying the One Over Many (OOM) hypothesis, for the family consisting of $S$ together with $\Sigma_1$ (\textit{Parmenides} 131a6--11).
    Indeed this is possible. Following the notation above, we set $\Sigma_2= \text{One}_k$, where $k$ is such that $\text{One}_k/\text{Being}_k=\text{One}/\text{Being}$.
    Here we use in an essential way the anthyphairetic periodicity of the One of the second hypothesis. Since the anthyphairesis $\text{Anth}(\text{One}_k, \text{Being}_k)$ of $\Sigma_2$ is equal to the anthyphairesis $\text{Anth}(\text{One}, \text{Being})$ of $\Sigma_1$, it follows that every member of $S$ participates in $\Sigma_2$, and $\Sigma_1$ is large because it is identical to $\Sigma_2$.
    \item[Step 2.2.] The Non-Identity (NI) hypothesis claims that $\Sigma_2$ is not identical with $\Sigma_1$.
    But this claim is against the equalization obtained by means of the anthyphairetic periodicity and the introduction of the dialectical numbers.
    We conclude that Step 2.2 is false.
\end{enumerate}

Plato simply states TMA in the Parmenides Introduction as a problem, wishing perhaps to make clear the fundamental difference between a sensible participating in an intelligible Being and an intelligible Being participating in an intelligible Being; in the first case participation is from an inferior to a superior entity, in the second case participation is between two equalized Beings, equalization being a consequence of anthyphairetic periodicity.

It follows then that Vlastos’ Non-Identity principle:
\[\text{If an Idea } \Sigma \text{ participates in an idea } \Sigma’, \text{ then } \Sigma \text{ is not identical to } \Sigma’.,\]
on which his analysis of the Third Man argument is based, is false, and with it Vlastos’ interpretation of the Third Man Argument breaks down; there is no regression to infinity, the intelligible Being is one and unique.

Notes
\begin{itemize}
    \item Vlastos’ reconstruction of TMA is inductive; our analysis implies that while Vlastos’ induction is true for the first step, the argument breaks down already at the second step. A similar failure of an induction at the second step in a mathematical example, claimed in W. Hurewicz and H. Wallman, 1941~\cite{Hurewicz1941}. was noted by V. V. Fedorchuk and J. van Mill, 2000~\cite{Fedorchuk2000}.
    \item Our comment on Rickless interpretation is the following:
    \begin{itemize}
        \item On the one hand, the negation of P [Purity] and the negation of NCC [No causation by Contraries] are indeed shown in the second hypothesis of the \textit{Parmenides}, e.g. by the equalization of the One and the Being, and equivalently in the Sophist by the equalization of Being and not Being.
        \item On the other hand, the U [Uniqueness] holds true for intelligible Beings, also by the equalization of the parts One and Being in the Second Hypothesis of the \textit{Parmenides}, equivalently the equalization of Being and not Being in the Sophist, and more generally by the equalization of all parts of the intelligible Being, namely the self-similar Oneness of intelligible Beings.
        \item On the contrary, Rickless mistakenly accepts the validity of the Third Man Argument and the denial of the Uniqueness U
    \end{itemize}
\end{itemize}

\subsection{The resolution of the Third Man Argument by means of the equalization of the One with each of its parts}
We set
\[\Sigma_1=\text{One}> \Sigma_2=\text{One}_k> \Sigma_3=\text{One}_{2k+1}> \Sigma_4=\text{One}_{3k+2},\dots> \Sigma_n=\text{One}_{nk+n-1}>\dots,\]
and we see that the sequence ($\Sigma_n$) constitutes an infinite multitude of intelligible Beings, such that $\Sigma_{n+1}$ is equalized to $\Sigma_n$ for all $n$.

Here let us note that one might think that it would be simpler to take
\[\Sigma_1=\text{One}> \Sigma_2=\text{Being}> \Sigma_3=\text{One}_1> \Sigma_4=\text{Being}_2>\dots> \Sigma_{2n-1}=\text{One}_n>\dots> \Sigma_{2n}= \text{Being}_n>\dots,\]
because then we would have not only that $\Sigma_{n+1}$ is equalized to $\Sigma_n$, but also that $\Sigma_n$ participates in $\Sigma_{n+1}$; but then there would be some problem with the participation of all the sensibles $S$ in $\Sigma_{n+1}$.
Instead, with our choice, we also have that $\text{Anth}(\text{One}_{nk+n-1}, \text{Being}_{nk+n-1})=\text{Anth}(\text{One}, \text{Being})$ for all $n$, hence every sensible $s$ participates in every $\Sigma_n$, for all $n$.
so that the intelligible Being $\Sigma_{n+1}$ is One Over both the sensibles $S$ and the intelligible being $\Sigma_n$ for all $n$.

\renewcommand{\refname}{Bibliography}
% \nocite{*} 
% \bibliographystyle{plain}
% \bibliography{negr251206}

%\printbibliography

\bigskip

Author's address:

\medskip

Stelios Negrepontis, 

Department of Mathematics,

University of Athens,

Athens, Greece

snegrep@gmail.com

\end{document}